\documentclass[11pt]{article}
\usepackage{amsmath, amssymb,anysize,indentfirst}
\def\qed{\hbox to 0pt{}\hfill$\rlap{$\sqcap$}\sqcup$}
\newtheorem{lemma}{Lemma}[subsection]
\newtheorem{theorem}[lemma]{Theorem}
\newtheorem{remark}[lemma]{Remark}
\newtheorem{proposition}[lemma]{Proposition}

\newtheorem{corollary}[lemma]{Corollary}
 \textheight
220mm \textwidth 145mm \voffset=-6mm
\begin{document}
\title{\textbf{Derivations for the even parts of
modular Lie superalgebras $W$ and $S$ of Cartan type}}
\author{Wende  Liu$^{1,2}$ and Yongzheng Zhang$^{1}$ \\
\\
\textit{$^{1}$Department of Mathematics,} \\
\textit{Harbin Normal University, Harbin 150080,   China}\\
\textit{Email: wendeliu@sohu.com}\\
\\
\textit{$^{2}$Department of Mathematics,} \\
\textit{Northeast Normal University, Changchun 130024, China}\\
\textit{Email: zhyz@nenu.edu.cn}}
\date{ }
\maketitle
\begin{quotation}
\small\noindent Let $\mathbb{F}$ be the  underlying base field of
characteristic $p>3 $ and denote by $\mathcal{W}$ and
$\mathcal{S}$ the even parts of  the finite-dimensional
generalized Witt Lie superalgebra $W$ and  the special Lie
superalgebra $S,$  respectively. We first give the generator sets
of the Lie algebras $\mathcal{W}$ and $\mathcal{S}.$  Using
certain properties of the canonical tori of $\mathcal{W}$ and
$\mathcal{S},$ we then determine the derivation algebra of
$\mathcal{W}$ and the derivation space of $\mathcal{S} $ to
$\mathcal{W},$ where $\mathcal{W}$ is viewed as $\mathcal{S}
$-module by means of the adjoint representation. As a result, we
describe explicitly the derivation algebra of $\mathcal{S}.$
Furthermore, we prove that the outer derivation algebras of
$\mathcal{W}$ and $\mathcal{S} $ are abelian Lie algebras or
metabelian Lie algebras with explicit structure. In particular, we
give the dimension formulae of the derivation algebras and outer
derivation algebras of $\mathcal{W}$ and $\mathcal{S}.$ Thus we
may make a comparison between the  even parts of
 the (outer) superderivation algebras
 of $W$ and $S$  and  the (outer) derivation algebras of the
 even parts of $W$ and $S,$ respectively.\\



\end{quotation}

\noindent {\Large\textbf{0. \ \  Introduction}}
\newline

\noindent This paper considers finite-dimensional Lie algebras and
Lie superalgebras over
 a field of prime characteristic.
  As is well known, the theory of Lie superalgebras over
  a field of characteristic zero has
 obtained plentiful fruits (see [4--5, 14]). But that is
 not the situation for modular Lie superalgebras, for example, the classification
 of finite-dimensional simple modular Lie superalgebras has not been completed. We sketch the
 recent development for modular Lie superalgebras. As far as we
 know, [6] may be the earliest paper on modular Lie
 superalgebras, in which the $p$-structure and $2p$-structure for
 modular Lie superalgebras
 (analogous to $p$-mappings for modular Lie algebras) are introduced.
 The (restricted) enveloping algebras for
  modular Lie superalgebras
 are studied in [12]. The reader is also referred to a paper on Frobenius extensions
 and restricted modular Lie superalgebras (see [3]).  In [19]
  the four families of finite-dimensional
 Cartan-type modular Lie superalgebras
 $X(m,n;\underline{t})$
are constructed and the simplicity and
 restrictiveness are studied, where $X=W,S,H,$ or $K$
 (These notations and other concepts
 mentioned in this introduction will be further explained in
 Section 1).
  These modular Lie superalgebras are analogous to both
 finite-dimensional modular Lie algebras of Cartan type and finite-dimensional
 Lie superalgebras of Cartan type over a field of characteristic
 zero (see [4, 15]). In a recent paper [9] (see also [10]),
 we consider a new family of
 finite-dimensional simple modular Lie superalgebras of Cartan type, which are
  analogous to neither the finite-dimensional
 modular Lie algebras of Cartan type, nor the finite-dimensional
  Lie superalgebras of Cartan type over a field of characteristic
 zero. It is therefore conceivable that the classification of
 finite-dimensional simple modular Lie superalgebras should not be
 trivial. For more information on modular Lie superalgebras, the
 reader is referred to [7--12, 19--23].

Given any Lie superalgebra $\frak g=\frak
g_{\overline{0}}\oplus\frak g_{\overline{1}} ,$ then the
superderivation  algebra of $\frak g$ is also a Lie superalgebra,
denoted by $\mathrm{Der}(\frak
g)=\mathrm{Der}_{\overline{0}}(\frak  g)\oplus
\mathrm{Der}_{\overline{1}}( \frak g).$ Since $\frak
g_{\overline{0}}$ is a Lie algebra, one can consider  the
derivation algebra $\mathrm{Der}(\frak{g}_{\overline{0}} ).$ A
question naturally arises: What are the relations between the
derivation algebra of the Lie algebra $\frak g_{\overline{0}}$ and
the superderivation algebra of $\frak g$? More precisely, if we
denote still by
 $\mathrm{Der}_{\overline{0}}(\frak {g}) $ the Lie algebra consisting
 of  restrictions  of the even
 superderivations   to
 $\frak{g}_{\overline{0}},$
 can we assert that
$\mathrm{Der}(\frak g_{\overline{0}})=
\mathrm{Der}_{\overline{0}}(\frak g )$? Equivalently, can we
assert that every derivation of the Lie algebra $\frak
g_{\overline{0}},$ the even part of $\frak{g},$ may   extend  to a
superderivation of the Lie superalgebra $\frak g$? In this paper,
as direct consequences of our results, these questions will be
answered  for the Lie superalgebras $W$ and $S$ of Cartan type
over a field  of prime characteristic.

Let $\mathbb{F}$ be the  underlying base field of characteristic
$p>3. $ Let $\mathcal{W}$ and $\mathcal{S} $ denote the even parts
of the Lie superalgebras $W$ and $S ,$  respectively. In this
paper we shall study the derivation algebras and the outer
derivation algebras of these two Lie algebras in a  systematic way
and one of the main purposes of this paper is to lay the
foundations for future studies on modular Lie superalgebras of
Cartan type. Let $\mathcal{L}$ denote  $\mathcal{W}$ or
$\mathcal{S},$ and $\mathrm{Der}(\mathcal{L},\mathcal{W})$ the
space of derivations of $\mathcal{L}$ to $\mathcal{W}, $ where
$\mathcal{W}  $ is viewed as $\mathcal{L}$-module by means of
adjoint representation. By our notation,
$\mathrm{Der}(\mathcal{W})=\mathrm{Der}(\mathcal{W},\mathcal{W})$
is the derivation algebra of $\mathcal{W}.$ Let
$\mathcal{L}=\oplus_{i\geq -1}\mathcal{L}_{i}$ be the natural
$\mathbb{Z}$-gradation of $\mathcal{L}.$ Just as in the case of
Lie superlagebras  $W$ and $S$ of Cartan type, one may prove that
any derivation in $ \mathrm{Der}(\mathcal{L},\mathcal{W})$ can be
reduced (by modulo a suitable inner derivation) to be a derivation
vanishing on $\mathcal{L}_{-1}.$ However, in contrast to the case
of Lie
 superalgebras $W$ and $S$ of
Cartan type,
one cannot obtain directly that a derivation
vanishing on $\mathcal{L}_{-1}$ must be zero, since $\mathcal{L}$  is neither transitive
 nor admissibly graded in general. Thus it is conceivable that the top
$\mathcal{L}_{-1}\oplus \mathcal{L}_{0}$  will play an important
role in the studies on the derivations. Indeed, just as we shall
see, any two derivations of nonnegative $\mathbb{Z}$-degree in
$\mathrm{Der}(\mathcal{L},\mathcal{W})$ which coincide on the top
$\mathcal{L}_{-1}\oplus\mathcal{L}_{0}$ differ with only an inner
derivation (see Corollaries 3.1.5 and 4.1.4), where
$\mathcal{L}=\mathcal{W}$ or $\mathcal{S}.$ This motivates us  to
consider whether  every derivation of nonnegative $\mathbb{Z}$-degree
in $\mathrm{Der}(\mathcal{L},\mathcal{W})$ can be reduced to be a
derivation vanishing on the top. For
$\mathrm{Der}(\mathcal{W})=\mathrm{Der}(\mathcal{W},\mathcal{W})$,
  using a known result which tells us that any homogeneous derivation
  of nonnegative $\mathbb Z$-degree of  a centerless $\mathbb Z$-graded Lie algebra
   may be reduced to be vanishing on the given torus contained in
   the component of degree zero
    (see [15, Proposition 8.4, p. 193]), we can obtain the desired result  by a brief
   argument (see Section 3.2). But the same work on the derivation space
    $ \mathrm{Der}(\mathcal{S},\mathcal{W})$ is more difficult, since we cannot
   apply the
    known result mentioned above in this case.
    This observation leads to the study on the exterior algebra $\Lambda(n)$ and the
    canonical torus of $\mathcal{S}$ (see Corollary 2.1.5),
    where $\Lambda(n)$ is viewed as a module of the canonical torus of $\mathcal{S}$
    in the obvious way. Then it is proved that any homogeneous derivation of nonnegative
    $\mathbb{Z}$-degree in $ \mathrm{Der}(\mathcal{S},\mathcal{W})$ can be reduced to be vanishing on
    the  canonical torus of $\mathcal{S}.$  As a result, we can show that
    the homogeneous derivations of nonnegative $\mathbb{Z}$-degree in
    $ \mathrm{Der}(\mathcal{L},\mathcal{W})$ are all inner. By giving
     the generator set of $\mathcal{L},$ we can compute the homogeneous derivations
     of negative $\mathbb{Z}$-degree. Finally,  the derivation algebras and
     the outer derivation algebras of $\mathcal{W}$ and $\mathcal{S}$
     are determined completely; in particular, we give the dimension formulae of the derivation
     algebras of $\mathcal{W}$ and $\mathcal{S}$.

      The original motivation for this paper comes from the encouragement of
      the anonymous referee for the
      paper [8].  Our work  is  motivated  by the
 results and methods on Lie algebras and Lie superalgebras  [2, 15, 17,
 23]
 and based on certain results in  [15, 23] on modular Lie algebras and Lie
 superalgebras of
 Cartan type. Certain results of this paper are closely parallel to those obtained by
 Celousov [2].  In particular, we use many ideas from [2,
 15] and benefit much from reading [13, 18]. For more information
 on derivations of modular Lie superalgebras of cartan type, the
 reader is referred to [8, 10--11, 17, 23].

The paper is organized as follows. In Section 1, we review the necessary notions
concerning Lie
algebras and Lie superalgebras and the notions of modular Lie superalgebras $W$ and
$S$ of Cartan type. In Section 2, we study certain subalgebras of the even parts
 $\mathcal{W}$ and $\mathcal{S} $  and give the generator sets of  $\mathcal{W}$ and
 $\mathcal{S}.$ We establish also some
 technical lemmas concerning the canonical tori of   $W$ and $S, $
 which will be used throughout this paper. In Section 3, we first study
 the derivations
 vanishing on the top $\mathcal{W}_{-1}\oplus \mathcal{W}_{0}.$ Then we characterize
 the homogeneous derivations of nonnegative $\mathbb{Z}$-degree and
 negative $\mathbb{Z}$-degree, respectively. As a result, the derivation algebra of
 $\mathcal{W}$ is determined. In Section 4, the derivation space
 $\mathrm{Der}(\mathcal{S}, \mathcal{W}) $ and the derivation algebra
 $\mathrm{Der}(\mathcal{S})$ are determined. In Section 5, by using the results obtained
 in Sections 3 and 4, the outer derivation algebras of $\mathcal{W}$ and $\mathcal{S} $
 are described explicitly and the dimension formulae of derivation algebras are given.

\section{Preliminary}

\subsection{Basic notion}

\noindent Let $\mathbb{F}$ be an arbitrary field in this subsection  and
$\mathbb{Z}_{2}=\{\overline{0}, \overline{1}\}$ be the field of
two elements. In this paper, all vector spaces,  linear mappings,
tensor products are over the underlying base field  $\mathbb{F}.$

Recall that a  \textit{vector superspace} is a
$\mathbb{Z}_{2}$-graded vector space $V=V_{\overline{0}}\oplus
V_{\overline{1}}. $ We denote by $\mathrm{p}(a)=\theta$ the
\textit{parity of a homogeneous element} $a\in V_{\theta},
\theta\in \mathbb{Z}_{2}.$ A subspace $U$ of a vector superspace
$V$ is by definition   $\mathbb{Z}_{2}$-graded; that is, $U=U\cap
V_{\overline{0}} \oplus U\cap V_{\overline{1}}.$

We assume
throughout that if $\mathrm{p}(x)$ occurs in an expression,
 then $x$ is assumed to be $\mathbb{Z}_2$-homogeneous.

A \textit{superalgebra} is a vector superspace $\mathcal{A}=\mathcal{A}_{\overline{0}}\oplus
\mathcal{A}_{\overline{1}}$ endowed with an algebra structure such that
$\mathcal{A}_{\theta}\mathcal{A}_{\mu}\subset \mathcal{A}_{\theta+\mu} $
for all
$\theta,\mu\in  \mathbb{Z}_{2}.$

A \textit{Lie superalgebra}  is a superalgebra satisfying the
super-anticommutativity and super-Jacobi identity (see [4 , 14]). Let
$\frak g=\frak g_{\overline{0}}\oplus \frak g_{\overline{1}} $ be
a Lie superalgebra. Then the even part $\frak g_{\overline{0}}$ is
a Lie algebra. Note that in the case $\mathrm{char} \mathbb{F}=2,$
a Lie superalgebra is a $\mathbb{Z}_{2}$-graded Lie algebra. Thus
one usually adopts the convention that $\mathrm{char}
\mathbb{F}=p>2 $ in the modular case.

Let $V=V_{\overline{0}}\oplus V_{\overline{1}} $ be a vector
superspace. The algebra $\mathrm{End}_{\mathbb{F}}(V)$ consisting
of the $\mathbb{F}$-linear mappings
 of $V$ into itself becomes an associative superalgebra if one
defines
$$\mathrm{End}_{\mathbb{F}} (V)_{\theta }
:=\{A\in \mathrm{End}_{\mathbb{F}} (V)\mid A (V_{\mu})\subset
V_{\theta +\mu}, \mu\in \mathbb{Z}_{2}\}
$$
for  $\theta\in \mathbb{Z}_{2}.$ On the
vector superspace $\mathrm{End}_{\mathbb{F}} (V)
=\mathrm{End}_{\mathbb{F}} (V)_{\overline{0}}\oplus
\mathrm{End}_{\mathbb{F}} (V)_{\overline{1}}$ we define a new
multiplication $[\ ,\ ]$ by
$$ [A,B]=AB-(-1)^{\mathrm{p}(A)\mathrm{p}(B)}BA
\quad \mbox{for  }\ A, B\in \mathrm{End}_{\mathbb{F}} (V).$$ This
algebra endowed with the new multiplication will be denoted by
 $\mathrm{pl}(V)=\mathrm{pl}_{\overline{0}}(V)\oplus\mathrm{pl}_{\overline{1}}(V);$
it is  a Lie superalgebra and is said to be the \textit{general
linear Lie superalgebra}.

Suppose that $\mathcal{A}=\mathcal{A}_{\overline{0}}\oplus
\mathcal{A}_{\overline{1}}$ is  superalgebra. Let
$$\mathrm{Der} _{\theta}(\mathcal{A}):
= \{ D\in\mathrm{pl}_{\theta}(\mathcal{A})\mid D(xy)
=D(x)y+(-1)^{\theta \mathrm{p}(x) }xD(y)\  \mbox{for  } x,y\in\mathcal{A}\} $$
for all $\theta \in \mathbb{Z}_{2}.$  Define
$$\mathrm{Der} (\mathcal{A}):=
 \mathrm{Der} _{ \overline{0}}(\mathcal{A})
 \oplus\mathrm{Der} _{ \overline{1}}(\mathcal{A}),$$
 then it is easy to see that
 $\mathrm{Der} (\mathcal{A})$ is a subalgebra of
  $\mathrm{pl}(\mathcal{A}),$ which is called the
  \textit{superderivation algebra of} $\mathcal{A}$.
 If $\theta=\overline{0}$ (resp. $\theta=\overline{1}$),  the elements in
$\mathrm{Der} _{ \theta}(\mathcal{A})$ are called \textit{even
superderivations} (resp. \textit{odd superderivation})\textit{ of} $\mathcal{A}$ and
the elements in $\mathrm{Der}(\mathcal{A})$ are called
\textit{superderivations of} $\mathcal{A}.$ For more details on
superderivations for Lie superalgebra, the reader is referred to
[14].

Let $\frak{g}$ be a  Lie algebra and $V$   an $\frak{g}$-module. A
linear mapping $D:\frak{g}\rightarrow V$ is called a
\textit{derivation of} $\frak{g}$  to $V$ if $D(xy)=x\cdot
D(y)-y\cdot D(x)$ for all $x,y\in \frak{g}.$ A derivation
$D:\frak{g}\rightarrow V$ is called \textit{inner} if there is
$v\in V $ such that $D(x)=x\cdot v$ for all $x\in \frak{g}.$
Following [15, p. 13], denote by $\mathrm{Der}(\frak{g},V)$ the
\textit{space of derivations} of $\frak{g}$  to $V.$  Then
$\mathrm{Der}(\frak{g},V)$ is an $\frak{g}$-submodule of
$\mathrm{Hom}_{\mathbb{F}}(\frak{g},V).$ Assume in addition that
$\frak{g}$ and $V$ are finite-dimensional and that
$\frak{g}=\oplus_{r\in \mathbb{Z}} \frak{g}_{r}$ is
$\mathbb{Z}$-graded and $V=\oplus_{r\in\mathbb{ Z}} V_{r}$ is a
$\mathbb{Z}$-graded $\frak{g}$-module. Then
$\mathrm{Der}(\frak{g},V)=\oplus_{r\in \mathbb{Z}}
\mathrm{Der}_{r}(\frak{g},V)$ is $\mathbb{Z}$-graded
$\frak{g}$-module  by setting
$$\mathrm{Der}_{r}(\frak{g},V):=\{D\in
\mathrm{Der}(\frak{g},V)\mid D(\frak{g}_{i})\subset V _{r+i} \
\mbox{for all} \ i \in\mathbb{Z}\}.$$ In the case of $V=\frak{g},$
the \textit{derivation algebra} $\mathrm{Der}(\frak{g})$ coincides
with $\mathrm{Der}(\frak{g},\frak{g}) $ and
$\mathrm{Der}(\frak{g}) =\oplus _{r\in \mathbb{Z}}
 \mathrm {Der}_{r}(\frak{g})$ is a $\mathbb{Z}$-graded Lie algebra.

 For a Lie superalgebra
  $\frak g=\frak g_{\overline{0}}\oplus \frak g_{\overline{1}},$
   the restriction  of an  even superderivation
   $ D\in \mathrm{Der}_{\overline{0}}(\frak g)$ to  $\frak g_{\overline{0}}$ is a
     derivation of  the Lie algebra $\frak g_{\overline{0}}.$
     For convenience, we shall  write
   still
   $ \mathrm{Der}_{\overline{0}}(\frak g)$  for the set
    $\{D|_ {\frak g_{\overline{0}}}\mid D\in \mathrm{Der}_{\overline{0}}(\frak g)\}.$
    By using this notation, it is easy to see that
    $\mathrm{Der}_{\overline{0}}(\frak g)\subset  \mathrm{Der} (\frak g_{\overline{0}})$ is a subalgebra
    $\mathrm{Der} (\frak g_{\overline{0}}).  $ As mentioned in the introduction,
    as one of the main results, we shall prove in this paper that
    the converse inclusion is also valid for
    $\frak g=W,S,$ the \textit{generalized Witt superalgebra} and
     the \textit{special superalgebra of Cartan type over a field
      of finite characteristic} (for a
     definition, see Section 2.2).

    If $\frak g=\oplus_{-r\leq i\leq s}\frak g_i$ is a $\mathbb
    Z$-graded Lie  algebra, then $\oplus_{-r\leq i\leq 0}\frak
    g_i$ is called the \textit{top of} $\frak g$ (with respect to the
    gradation).

\subsection{Modular Lie superalgebras $W$ and $S$}
\noindent In this subsection we review the notions of  modular Lie
superalgebras $W$ and $S$ of Cartan type and their gradation
structures.

In the following sections,  $\mathbb {F}$ denotes a field of
characteristic $p>3.$ In addition to the standard notation
$\mathbb {Z},$ we write  $\mathbb {N}$ for the set of positive
integers, and ${\mathbb {N}}_0$ for the set of nonnegative
integers.
Henceforth, we will let  $m,n$  denote fixed integers in
$\mathbb {N}\setminus \{1,2\} $ without notice. For $\alpha = (
\alpha _1,\ldots ,\alpha _m ) \in \mathbb {N}_0^m,$ we put
$|\alpha| =\sum_{i=1}^m\alpha _i.$ Let $\frak{A}(m)$ denote the
\textit{divided power algebra over} $\mathbb{F}$ with an ${\mathbb
F}$-basis $\{ x^{( \alpha ) }\mid \alpha \in \mathbb{N}_0^m \} .$
For $\varepsilon _i=( \delta_{i1}, \ldots ,\delta _{im}) ,$ we
abbreviate $x^{( \varepsilon _i)}$ to $x_i,$ $i=1,\ldots ,m.$ Let
$\Lambda (n)$ be the \textit{exterior  superalgebra over} $\mathbb{F}$ in
$n$ variables $x_{m+1},\ldots ,x_{m+n}.$ Denote the tensor product
by $\frak{A}(m,n) =\frak{A}(m)\otimes_{\mathbb{F}} \Lambda(n).$
Obviously, $\frak{A}(m,n)$ is an associative superalgebra with a
$\mathbb{Z}_2$-gradation induced by the trivial
$\mathbb{Z}_2$-gradation of $\frak{A}(m)$ and the natural
$\mathbb{Z}_2$-gradation of $\Lambda (n).$ Moreover,
$\frak{A}(m,n)$ is super-commutative.

For $g\in \frak{A}(m) ,f\in \Lambda(n),$ we write $gf $  for $
g\otimes f.$   The following formulas hold in $\frak{A}(m,n):$
$$
x^{(\alpha) }x^{(\beta) }=\binom{\alpha +\beta }{ \alpha}
 x^{( \alpha +\beta)}\quad\mbox{for}\  \alpha ,\beta \in
{\mathbb N}_0^m;
$$

$$
x_kx_l=-x_lx_k\quad\mbox{for}\ k,l=m+1,\ldots ,m+n;
$$
$$
x^{( \alpha ) }x_k=x_kx^{( \alpha ) }\quad\mbox{for}\ \alpha \in
\mathbb{N}_0^m, k=m+1,\ldots ,m+n,
$$
where $\binom{ \alpha +\beta} {\alpha}:=\prod_{i=1}^m\binom{
\alpha _i+\beta _i}{ \alpha _i}.$

Put $Y_0:=\{ 1,2,\ldots ,m \} ,$  $Y_1:=\left\{m+1,\ldots
,m+n\right\} $ and $Y:=Y_0\cup Y_1.$ For convenience, we adopt the
notation $r':=r+m$ for $r=1,\ldots,n.$ Thus, $Y_1:=\left\{ 1', 2',
\ldots ,n'\right\}.$ Set
$$\mathbb{B}_k:=\left\{ \langle i_1,i_2,\ldots
,i_k\rangle |m+1\leq i_1<i_2<\cdots <i_k\leq m+n\right\} $$ and $
\mathbb{B}:=\mathbb{B}(n)=\bigcup\limits_{k=0}^n\mathbb{B}_k,$
where $\mathbb{B}_0:=\emptyset.$ For $u=\langle i_1,i_2,\ldots
,i_k\rangle \in \mathbb{B}_k,$ set $ |u| :=k,| \emptyset | :=0$,
$x^\emptyset :=1,$ and $x^u:=x_{i_1}x_{i_2}\ldots x_{i_k};$
 we use  also $u$ to stand for the set $\{i_1,i_2,\ldots
,i_k\}.$ Clearly,
$\left\{ x^{\left( \alpha \right) }x^u\mid\alpha \in
\mathbb{N}_0^m,u\in \mathbb{B}  \right\} $ constitutes an
$\mathbb{F}$-basis of $\frak{A} \left( m,n\right) .$

Let $D_1,D_2,\ldots ,D_{m+n}$ be the linear transformations of
$\frak{A} \left( m,n\right) $ such that
\[
D_r ( x^{\left( \alpha \right) }x^u ) =\left\{
\begin{array}{l}
x^{\left( \alpha -\varepsilon _r\right) }x^u,\quad \quad \quad
\quad r\in Y_0
\\
x^{\left( \alpha \right) }\cdot \partial x^u/\partial x_r,\quad
\quad r\in Y_{1.}
\end{array}
\right.
\]
Then $D_1,D_2,\ldots ,D_{m+n}$ are superderivations of the
superalgebra $\frak{A} \left( m,n\right).$ Let
\[
W\left(m,n\right) =\Big\{ \sum\limits_{r\in Y}f_r D_r \mid f_r\in
\frak{A} \left( m,n\right) ,r\in Y\Big\} .
\]
Then $W\left( m,n\right) $ is a Lie superalgebra, which is
contained in ${\rm Der} ( \frak{A} \left(m,n\right)). $

Obviously, ${\rm p}( D_i) =\tau (i) ,$ where
\[
\tau \left( i\right) :=\left\{
\begin{array}{l}
\overline{0},\quad \quad i\in Y_0 \\
\overline{1},\quad \quad i\in Y_{1.}
\end{array}
\right.
\]
One may verify that
\[[fD,gE]=fD(g)E-(-1)^{{\rm
p}(fD){\rm p}(gE)}gE(f)D+(-1)^{{\rm p}(D){\rm
p}(g)}fg[D,E]\eqno(1.2.1)
\]
for $ f,g\in \frak{A} (m,n) ,$ $ D,E\in {\rm Der}\ \frak{A} (m,n
).$ In particular, the following formula holds in $W\left(
m,n\right):$
\begin{equation*}
\left[ fD_r,gD_s\right] =fD_r\left( g\right) D_s-\left( -1\right)
^{{\rm p}(fD_r) {\rm{p}}(gD_s) }gD_s\left( f\right) D_r
\end{equation*}
for $f,g\in \frak{A} \left( m,n\right) ,r,s\in Y.$

Let
\[
\underline{t}:=\left( t_{1},t_{2},\ldots ,t_m\right) \in
\mathbb{N}^m,\quad \pi:=\left( \pi _1,\pi _2,\ldots ,\pi _m\right)
\]
where $\pi _i:=p^{t_i}-1,i\in Y_0.$ Let
$\mathbb{A}:=\mathbb{A}\left( m;\underline{t}\right) =\left\{
\alpha \in \mathbb{N}_0^m\mid\alpha _i\leq \pi _i,i=1,2,\ldots
,m\right\}.$ Then
\[
\frak{A} \left( m,n;\underline{t}\right) :={\rm
span}_{\mathbb{F}}\left\{ x^{\left( \alpha \right) }x^u \mid
\alpha \in \mathbb{A},u\in \mathbb{B} \right\}
\]
is a finite-dimensional subalgebra of $\frak{A} \left(m,n\right) $
with   a natural ${\mathbb Z}$-gradation $\frak{A}
\left(m,n;\underline{t}\right)=\bigoplus _{r=0}^{\xi}
\frak{A}(m,n;\underline{t})_{r}$ by putting
        $$\frak{A}(m,n;\underline{t})_{r}:=
        {\rm span} _{\mathbb F}\{ x^{(\alpha)} x^{u}\mid |\alpha|+|u|=r \},
        \quad \xi:=|\pi|+n.$$
Set
\[
W\left( m,n;\underline{t}\right):=\Big\{ \sum\limits_{r\in Y}
f_rD_r \mid f_r\in \frak{A} \left( m,n;\underline{t}\right) ,r\in
Y\Big\} .
\]
Then $W\left( m,n;\underline{t}\right) $ is a subalgebra of
$W\left( m,n\right).$ In particular, it is a finite-dimensional
simple Lie superalgebra (see [19]). Obviously, $W (
m,n;\underline{t} ) $ is a free $ \frak{A} \left(
m,n;\underline{t}\right) $-module with $ \frak{A} \left(
m,n;\underline{t}\right) $-basis $ \{ D_r\mid r\in Y \}.$ We note
that $W ( m,n;\underline{t} ) $ possesses a \textit{standard}
$\mathbb F$-\textit{basis} $\{x^{(\alpha)}x^{u}D_r\mid \alpha\in
\mathbb{A}, u\in \mathbb{B}, r\in Y\}.$

Let $r,s\in Y$ and $D_{rs}:\frak{A} (
m,n;\underline{t})\rightarrow W(m,n;\underline{t})$ be the linear
mapping  such that
\begin{equation*}
 D_{rs}(f)=(-1) ^{ \tau(r) \tau (s)}D_{r}(f)D_{s}-(-1)
^{(\tau(r)+ \tau (s)){\rm p}(f)}D_{s}(f)D_{r}\quad\mbox{for}\ f\in
\frak{A} ( m,n;\underline{t}).
\end{equation*}
Then the following equation holds:
\begin{equation*}
[D_{k},D_{rs}(f)]=(-1) ^{\tau (k)\tau(r)
}D_{rs}(D_{k}(f))\quad\mbox{for}\ k,r,s \in Y; \ f\in \frak{A} (
m,n;\underline{t}).\eqno(1.2.2)
\end{equation*}
Put $$ S(m,n;\underline{t}):={\rm span} _{\mathbb F}\{
D_{rs}(f)\mid  r,s\in
Y; \ f\in \frak{A} (m,n;\underline{t})\}.$$
Then $ S(m,n;\underline{t})$ is a finite-dimensional simple
Lie superalgebra (see [19]).

Let
$ \mathrm {div}:W(m,n;\underline{t})\rightarrow \frak
{A}(m,n;\underline{t})$
be the \textit{divergence}  such that
$$
\mathrm {div}\big ( \sum_{r\in Y} f_{r}D_{r}\big )=\sum_{r\in Y }
(-1) ^{\tau(r) \mathrm{p}(f_{r})} D_{r}(f_{r}).
$$
Direct
computation shows that $\mathrm{div}$ is superderivation of
$W(m,n;\underline{t})$ to $\frak {A}(m,n;\underline{t})$ (see [5]
or [19]); that is,
 $${\mathrm {div}}[D,E]=D{\mathrm
{div}}(E)-(-1)^{\mathrm{p}(D)\mathrm{p}(E)}E{\mathrm {div}}(D)\quad \mbox{for all}\
D,E\in W(m,n;\underline{t}).
$$
Following [19], put
$$
\overline{ S }(m,n;\underline{t}):=\{D\in W(m,n;\underline{t})\mid
\mathrm{div}(D)=0\}.
$$
Then $
\overline{ S }(m,n;\underline{t})$
is a subalgebra of
 $W(m,n;\underline{t})$
  and $  S (m,n;\underline{t}) $ is a subalgebra of
  $\overline{ S }(m,n;\underline{t}).$

In the following sections, $W(m,n;\underline{t}),$
$S(m,n;\underline{t}),$ $\overline{S}(m,n;\underline{t}),$ and $\frak{A} (m,n;\underline{t})$ will be
denoted  by $W,S, \overline{S},$ and $\frak{A},$ respectively. In addition,
the the even parts of $W,$ $S,$ and $\overline{S} $ will be denoted by $\mathcal{W},$
$\mathcal{S}, $ and $\overline{\mathcal{S}}, $ respectively; that is,
$\mathcal{W}:=W_{
\overline{0} },$   $\mathcal{S}:=S_{\overline{0}}, $
 and $ \overline{\mathcal{S}}:=\overline{S}_{\overline{0}}.$

We  note also that for the sake simplicity, as mentioned above we
assume throughout that $\mathrm{char} \mathbb{F}>3$ and that the
parameters $m,n>2 $ although sometimes a weaker hypothesis is
sufficient.

\section{Subalgebras and generator sets of $\mathcal{W}$ and $\mathcal{S}$}
\noindent In this section, we present certain results on   some
 subalgebras of $\mathcal{W}$ and then give the generator
sets of $\mathcal{W}$ and $\mathcal{S}.$ These results will be
frequently used in the sequel.

\subsection {Subalgebras }

\noindent In this subsection, we deal with certain subalgebras of
$\mathcal{W},$ which are   important   for   future studies in
this paper. In particular, we will study the property of the canonical
torus subalgebras  of $\mathcal{W}$ and  give a reduction
proposition (Proposition 2.1.6) for derivations of $\mathcal{L}$
to $\mathcal{W},$ where $\mathcal{L}$ is a $\mathbb{Z}$-graded subalgebra of
$\mathcal{W}$ satisfying $\mathcal{L}_{-1}=\mathcal{W}_{-1}.$

Let
$$\mathcal{G}:=\mathrm{span}_{\mathbb{F}}\{x^uD_r
\mid r\in Y, u\in \mathbb{B} , \mathrm{p}(x^uD_r)=\overline{0}\}.
$$ Then $\mathcal{G}=C_\mathcal{W}(\mathcal{W}_{-1}),$ the centralizer of
$\mathcal{W}_{-1}$ in $\mathcal{W}.$ Therefore, $\mathcal{G}$ is a
$\mathbb{Z}$- graded subalgebra of $\mathcal{W}.$ Put
$\mathcal{G}_{i}:=\mathcal{G}\cap \mathcal{W}_{i} $ and
$$
E(\mathcal{G}):=\oplus _{r\in \mathbb{Z}}\mathcal{G}_{2r},\quad
O(\mathcal{G}):=\oplus _{r\in \mathbb{Z}}\mathcal{G}_{2r+1}.
$$
Since  $[O(\mathcal{G}), O(\mathcal{G})]=0,$ $O(\mathcal{G})$ is
an ideal of $\mathcal{G}.$ It is easily seen that
$$
\mathcal{G}/ O(\mathcal{G})\simeq E(\mathcal{G})\simeq
W(n)_{\overline{0}},
$$
where $W(n)$ is the simple Lie superalgebra of Cartan type (see
[4, p. 57]) and $W(n)_{\overline{0}}$ is the even part of $W(n).$

Let $\frak{g}=\oplus_{q=-r}^s\frak{g}_{q}$ be a
$\mathbb{Z}$-graded Lie
 algebra. Recall that $\frak{g}$ is called \textit{transitive }(with
respect to the $\mathbb{Z}$-gradation) provided that
$\{x\in\frak{g}_{q}\mid [x,\frak{g}_{-1}]=0\}=0$ for all $q\in
\mathbb{N}_{0}.$ We say that $\frak{g}$ is   \textit{admissibly
graded} if $C_{\frak{g}}(\frak{g}_{-1})=\frak{g}_{-r}.$

By the remarks above, $\mathcal{W}$ is neither transitive nor
admissibly graded. In particular, $\mathcal{W}$ is not a simple
Lie algebra. Noticing that $\mathcal{S}_{-1}=\mathcal{W}_{-1} $
and $\mathcal{S}\cap \mathcal{G}_{0}\neq 0,$ we have the same
conclusion for $\mathcal{S}.$

We shall use frequently the following simple fact. The proof is
straightforward and therefore is omitted.

\begin{lemma} 
 Let $\mathcal{L}$ be a
$\mathbb{Z}$-graded subalgebra of $\mathcal{W}$ such that $
\mathcal{L}_{-1}=\mathcal{W}_{-1}.$ Suppose that $\phi\in
\mathrm{Der}(\mathcal{L},\mathcal{W})$ and
$\phi(\mathcal{L}_{-1})=0.$ If $E$ is an element of $\mathcal{L}$
such that $[E, \mathcal{W}_{-1}]\subset \ker(\phi),$ then
$\phi(E)\in \mathcal{G}.$\qed
\end{lemma}

Put $\Gamma_{r}:=x_rD_r$ for $r\in Y$  and $\Gamma=\sum_{r\in
Y}\Gamma_{r}.$ Then $\Gamma$ is called the \textit{degree derivation of}
$\frak{A}(m,n;\underline{t}).$  Correspondingly, $\mathrm{ad}
\Gamma $ is called the \textit{degree derivation of} $\mathcal{W}.$  Put
$\Gamma':=\sum_{k\in Y_1}\Gamma_{k}$ and $\Gamma'':=\sum_{i\in
Y_0}\Gamma_{i}.$

In this paper we adopt the following notation. Let $P$ denote a proposition.
Put $\delta _{P}:=1$ if $P$ is true and $\delta _{P}:=0,$ otherwise.
We need the following computational lemma.

 \begin{lemma}
  Let $i \in Y_0,$   $k\in Y_1.$
$\alpha\in \mathbb{A}$ and $u\in \mathbb{B}.$ Then the following
statements hold.

 $\mathrm{(i)}$\ \ \ \ \
$\Gamma_{i}(x^{(\alpha)}x^u)=\alpha_{i}x^{(\alpha)}x^u.$

$\mathrm{(ii)}$\ \ \ \ $\Gamma_k(x^{(\alpha)}x^u)=\delta_{k\in
u}x^{(\alpha)}x^u;$ in particular, $\Gamma_k^2= \Gamma_k.$

$\mathrm{(iii)}$\ \ \ For $f\in\frak{A} _{r}, r\in \mathbb{N}_0,$
$\Gamma(f)=rf.$

$\mathrm{(iv) }$ \ \   Let $r\in \mathbb{Z} $ and $D\in
\mathcal{W}_{r}.$
 Then $ [\Gamma, D]=rD.$

 $\mathrm{(v)}$\ \ \  \   $\Gamma_{i}^p=\Gamma_{i}.$

 $\mathrm{(vi)}$\ \ \ Let $D\in \mathcal{G}_{r}.$ If $r$ is even,
 then $[\Gamma',D]=rD;$ if $ r$ is odd, then $[\Gamma',D]=(r+1)D.$

 $\mathrm{(vii)}$\ \ Every standard basis element of $W$ is an eigenvector
  of $\mathrm{ad}\Gamma _r$ for all
 $ r\in Y.$

 $\mathrm{(viii)}$ $[\Gamma'',x^{(\alpha)}x^{u}D_{i}]
 =(|\alpha|-1)x^{(\alpha)}x^{u}D_{i}$ for all
 $\alpha\in \mathbb{A}$,$u\in \mathbb{B} $ and $i\in Y_{0}.$
\qed
\end{lemma}

Let $(\frak g, [p])$ be a \textit{restricted Lie algebra}. Following
[16, p. 119] and [15, p. 79], define $x\in \frak g$ to be
$p$-\textit{semisimple} if $x\in \sum _{r\in
\mathbb{N}}\mathbb{F}x^{[p]^{r}}.$ If $x^{[p]}=x$ (and hence is
$p$-semisimple) we say that $x$ is \textit{toral}. An abelian restricted
subalgebra $\frak T$ of $\frak g$ is called a \textit{torus} if every
element in $\frak T$ is $p$-semisimple.

If $ \frak g $ is a (not necessarily restricted) Lie algebra with
trivial center, then $ \frak g $ may be identified with a
subalgebra of the restricted Lie algebra $\mathrm{Der}(\frak g).$
Following [1, p.97], we say that a subalgebra $\frak T$ is a
\textit{torus} of $\frak g$ if $\frak T$ is a torus of
$\mathrm{Der}(\frak {g}). $

Let us consider $\mathcal{W}$ and $\mathcal{ S}.$ Note that
$\mathcal{W}$ and $\mathcal{ S}$ are all centerless. Put
$\mathcal{T}:=\sum _{r\in Y}\mathbb {F}\cdot \Gamma _{r}$ and
$\mathcal{T}_{\mathcal{S}}:=\mathcal{S}\cap \mathcal{T}.$ Then by
Lemma 2.1.2, $\mathcal{T}$ and $\mathcal{T}_{\mathcal{S}}$ are
tori of $\mathcal{W}$ and  $\mathcal{S},$ which are called
\textit{canonical  tori}.

  The following lemma will be heavily used in this paper. It is
  essentially a generalization of [15, Proposition 8.2, p. 192].

\begin{lemma}

Let $ V$ be a vector space over $\mathbb{F}$ and $v_1, v_2,
\ldots, v_k\in V.$  Suppose that $A_{i} \in
\mathrm{End}_{\mathbb{F}}V$ is \textit{generalized invertible};
that is, there is $B_{i} \in
 \mathrm{End}_{\mathbb{F}}V$ such that $$A_{i}B_{i} A_{i}= A_{i},  \
 1\leq i\leq k.\eqno(2.1.1)$$
 If the following conditions are satisfied, then there exists $v\in
V$ such that $A_{i}(v)=v_{i}$ for $1\leq i\leq k.$

$\mathrm{(i)}$\ \  $A_{1},A_{2},\ldots,A_{k } $ commute mutually;

$\mathrm{(ii)}$ $A_{i}(v_{j})=A_{j}(v_{i}),\ 1\leq i,j \leq k$;

 $\mathrm{(iii)}$
$$A_{i}B_{i}(v_{i})=v_{i},\ 1\leq i\leq k,\eqno(2.1.2)$$
$$A_{i} B_{j}=B_{j} A_{i}
\ \mbox{whenever}\ i\neq j.\eqno(2.1.3)$$
\end{lemma}

\noindent \textit{Proof.}  We use induction on $k.$ When $k=1,$
  putting $v:=B_{1}(v_{1}),$ one obtains from (2.1.2) that
 $A_{1}(v)=A_{1}B_{1}(v_{1})=v_{1}.$

Let $k\geq 2$ and assume that the  conclusion holds for  $k-1.$
Then there exists $w\in V$ such that $A_{i}(w)=v_{i},\
i=1,\ldots,k-1.$ Set $v:=w+B_{k}(v_{k}-A_{k}(w)).$ Then, for
$i=1,\ldots,k-1,$ we obtain from (2.1.3), (i) and (ii) that

  \begin{eqnarray*}
  A_{i}(v)&=& A_{i}(w) +A_{i}B_{k}(v_{k}-A_{k}(w))\\
  &=&v_{i}+B_{k}A_{i}(v_{k}-A_{k}(w))\\
  &=&v_{i}+B_{k}(A_{i}(v_{k})-A_{k}A_{i}(w))\\
  &=&v_{i}+B_{k}(A_{k}(v_{i})-A_{k}(v_{i}))\\
  &=&v_{i}.
 \end{eqnarray*}
 Moreover, by (2.1.1) and (2.1.2) we have
    \begin{eqnarray*}
  A_{k}(v)&=& A_{k}\Big(w+B_{k}\big(v_{k}-A_{k}(w)\big)\Big)\\
  &=&A_{k}(w)+A_{k}B_{k}\big(v_{k}-A_{k}(w)\big)\\
  &=&A_{k}(w)+A_{k}B_{k}(v_{k})-A_{k}B_{k}A_{k}
  (w)\\
  &=&A_{k}(w)+v_{k}-A_{k}(w)\\
  &=& v_{k}.
 \end{eqnarray*}
  The proof is complete.\qed
 \newline

Using Lemma 2.1.3, we give a   result analogous to [15, Proposition 8.2,
p.192], which will be used to prove the main result in this
subsection (Proposition 2.1.6).

   \begin{corollary} 

  Let $r\leq m $ be an integer and
   $\ f_{1},\ldots,f_{r}\in \frak{A}.$
  Suppose that

$  \mathrm{(a)}$  $D_{i}(f_{j})=D_{j}(f_{i}),\ 1\leq i,j \leq r;$

  $\mathrm{(b)}$  $D^{p^{t_{i}}-1}_{i}(f_{i})=0,\ 1\leq i\leq r.$\\
  Then there is $f\in \frak{A}$
  such that $D_{i}(f)=f_{i}$ for $ 1\leq i\leq r.$
  \end{corollary}

  \noindent \textit{Proof.}
  Define
  $B_{i}:\frak{A}\rightarrow \frak{A}$
  to be a linear transformation such that
  $$
  B_{i}(x^{(\alpha)}x^{u})= x^{(\alpha+\varepsilon_{i})}x^{u}
     \quad\mbox{for}\ i\in Y_0$$
  where $x^{(\alpha+\varepsilon_{i})}:=0$
  for $\alpha+\varepsilon_{i}
  \not\in A(m;\underline{t}).$ We check the conditions of Lemma 2.1.3.  Evidently,
  $D_{1},\ldots,D_{r}$
  satisfy the condition (i). (b) shows that $ B_{1},\ldots, B_{r}$
  satisfy (2.1.1) and that (iii) holds.
  By (a),
  $\ f_{1},\ldots,f_{r}$
  satisfy the condition (ii) in Lemma 2.1.3. The proof is complete.
 \qed\newline

  Obviously, $ \Lambda(n) $ is an invariant  subspace of $\frak A $ under
   $\Gamma _{k}$ for all $k\in Y_{1}.$   Then $\Gamma
   _{k}\in \mathrm{End}_{\mathbb{F}}( \Lambda(n)).$ The following
   corollary will be used in Section 4  to consider how a derivation $\phi$ of
   $\mathcal{S}$ to $\mathcal{W}$ is determined by the action of $\phi$
   on $\mathcal{S}_{0} $ (see Lemma 4.2.4).

  \begin{corollary}

  Let $r\leq m$ be a positive integer and
   $\ f_{1},\ldots,f_{r}\in \Lambda(n).$
  Suppose that
  $\Gamma_{q_{i}}=x_{q_{i}}D_{q_{i}},\ q_{i}\in Y_{1},\ 1\leq i\leq r,$
  satisfy the following conditions:

  $\mathrm{(a)} $  $\Gamma_{q_{i}}(f_j)=\Gamma_{q_{j}}(f_i),\ 1\leq i \leq r;$

  $\mathrm{(b)}$   $\Gamma_{q_{i}}(f_i)=f_i,\  1\leq i\leq r$.\\
  Then there exists $f\in \Lambda(n) $
  such that $\Gamma_{q_{i}}(f)=f_i$ for $ i=1,2,\ldots,r.$
  \end{corollary}

  \noindent \textit{Proof.}  Set $B_i=\mbox{id}_{\Lambda(n)}.$ We check
  the conditions
  of Lemma 2.1.3. Since $\Gamma^2_{q_{i}}=\Gamma_{q_{i}} $ by Lemma
  2.1.2(ii), (2.1.1) holds.
  The remaining conditions may be easily checked.
  \qed\newline

  Let $V=\oplus _{r\in \mathbb{Z}}V_{r}$ be a $\mathbb{Z}$-graded
  vector space. For a $\mathbb{Z}$-homogeneous element $v\in V$,
  we denote by $\mathrm{zd}(v)$ the $\mathbb{Z}$-degree of $v.$

   For the sake of convenience, in the sequel we  usually
   write
   $i,j$  for the elements in $Y_0,$ and write  $k,l$ for the elements in
   $Y_1.$
  View $\mathcal{W}$ as an $\mathcal{L}$-module by   the adjoint
  representation. Recall the definition of a derivation of a Lie
  algebra to its module (see Section 1.1).
  Let $\mathcal{L}$
  be a $\mathbb{Z}$-subalgebra of $\mathcal{W}$ satisfying
  $\mathcal{L}_{-1}=\mathcal{W}_{-1}.$  Just as [15, Proposition 8.3, p. 193],
  the following proposition
  reduces   derivations of $\mathrm{Der} (\mathcal{L},
  \mathcal{W})$ to    the  derivations vanishing on
  $\mathcal{L}_{-1}.$

  \begin{proposition} 

  Let $\mathcal{L}$
  be a $\mathbb{Z}$-subalgebra of $\mathcal{W}$ such that
  $\mathcal{L}_{-1}=\mathcal{W}_{-1}.$
  If $\phi \in \mathrm{Der}_{t}(\mathcal{L},\mathcal{W})$
  where $t:=\mathrm{zd}(\phi)\geq 0,$
  then there exists $E\in \mathcal{W}_{t}$
  such that
  $$\big(\phi -\mathrm{ad}E\big)(\mathcal{L}_{-1})=0.$$
  \end{proposition}

  \noindent  \textit{Proof.}
  Suppose that $\phi(D_{i})=\sum_{r\in Y}f_{ri}D_{r}$
  for $i\in Y_{0},$ where $\ f_{ri}\in \frak{A}.$
  Applying $\phi$
  to the identity $[D_i,D_j]=0$ for $ i,j\in Y_{0},$
  we have
  $[\phi\big(D_i\big),D_j]+[D_i,\phi(D_j)]=0.$
 Then
  $$\Big[\sum_{r\in Y}f_{ri}D_{r},D_j\Big]+\Big[D_i,\sum_{r\in Y}f_{rj}D_{r}\Big]=0.$$
 Consequently,
  $$\sum_{r\in Y}D_{i}(f_{rj})D_{r}-\sum_{r\in Y}D_{j}(f_{ri})D_{r}=0.$$
  Since $\{D_r\mid r\in Y\}$
  is a free basis  of $\frak A$-module $\mathcal{W},$
  we have
  $$D_{i}(f_{rj})=D_{j}(f_{ri})\quad \mbox{for all}\ i,j\in Y_0,\ r\in Y.$$
  This implies that for fixed $r\in Y ,$ $ \{f_{r1},\ldots,f_{rm}\}$
  satisfies the condition (a) in Corollary  2.1.4
   (see also [15, Proposition 8.2, p. 192]).
  On the other hand, since
  $\big(\mathrm{ad}D_i\big)^{\pi_{i}+1}=0$ for  $i\in Y_0,$
  it follows from [15, Lemma 8.1, p. 191]
  that $\big(\mathrm{ad}D_i\big)^{\pi_{i}}\big(\phi(D_i)\big)=0.$
  Consequently, $D_i^{\pi_{i}}(f_{ri})=0$ for $r\in Y$ and $i\in Y_0;$
  that is, $ \{f_{r1},\ldots,f_{rm}\}$
  fulfills the condition (b) in Corollary 2.1.4.
  Therefore, there is  $g_r\in \frak{A}$
  such that
  $$D_i(g_r)=f_{ri} \quad\mbox{for all}\  i\in Y_0.$$
  Let $E''=-\sum_{r\in Y}g_rD_r.$
  Then for $i\in Y_0,$ we have
  $$[E'',D_i]=\sum_{r\in Y}D_i(g_r)D_r=\sum_{r\in Y}f_{ri}D_r
  =\phi(D_i).\eqno(2.1.4)$$
  Note that $\mathrm{zd}(\phi)= t.$
  (2.1.4) shows  that
  $$[E''_t,D_i]=\phi(D_i) \quad\mbox{for}\ i\in Y_0.\eqno(2.1.5)$$
  Let $E'=E''_{t}.$
  Then (2.1.5) implies that $(\phi-\mathrm{ad}E')\big|_{\mathcal{L}_{-1}}=0.$
  Let $E=E'_{\overline{0}}$ (the even part of $E'$). Then $E\in
  \mathcal{W}_{t}$  and
  $(\phi-\mathrm{ad}E)\big|_{\mathcal{L}_{-1}}=0.$
 \qed

   \begin{remark}

   {\rm Since $\mathcal{L}_{-1}$ is an
    abelian Lie algebra, one may consider the derivation space
    $\mathrm{Der}(\mathcal{L}_{-1},\mathcal{W}).$
    Then, under the conditions of Proposition 2.1.6,
    $\phi | _{\mathcal{L}_{-1}} \in \mathrm{Der}(\mathcal{L}_{-1},\mathcal{W})$
    is  inner.}
    \end{remark}

  \begin{remark}

   {\rm In Proposition 2.1.6 the
  element $E$ does not necessarily lie in $\mathrm{Nor}_{\mathcal{W}}(\mathcal{L}).$
  In addition,  in contrast to the  cases of Cartan type Lie algebras and
   Lie superalgebras (see [15, Proposition 8.3, p. 193]
   and \cite[Lemma 5]{ZZ}),
    we cannot prove directly that
   $\phi =(\mathrm{ad}E)\big|_\mathcal{L}, $ since $\mathcal{L}$ is
   not necessarily
    admissibly graded, as noted at the beginning of this
    subsection.}
    \end{remark}

  \subsection{Generator sets}

  \noindent As is well known, the generators play an
  important role in the study on derivations of an $\mathbb{F}$-algebra. More precisely, a
  derivation is completely determined by its action on the generator set. When
  we consider the derivation algebra of a given algebra, it is
  natural to hope finding a `good' generator set, which is
  convenient for computing derivations.

  In this subsection, we study mainly the generator sets of
   $\mathcal{W}$ and
  $\mathcal{S} $ for the purpose mentioned above. We know that
  $W(m,n;\underline{t})$ is a free
   $\mathbb{Z}$-graded module of the
  associative super-commutative algebra
  $\frak A(m,n;\underline{t}) $ with basis $\{D_{r}\mid r\in Y\} $
  consisting of special derivations of $W(m,n;\underline{t}).$ The `coefficients' in
  $\frak {A} (m,n;\underline{t})$
  are tensor products of divided power series and exterior
  products. Thus it is desirable to find a generator set which
  consists of certain elements of low $\mathbb{Z}$-degree and certain elements
  of high $\mathbb{Z}$-degree but of simple form, such as the elements in
    Lie algebras $W(m;\underline{t})$ or
  $S(m;\underline{t}) $ of Cartan type.

  The results obtained in this subsection
  will be used frequently in the sequel.

  Recall that $W=W(m,n;\underline{t})$ and $\mathcal{W}$ denotes the even part of
  $W.$ Recall also our  notations that
  $\underline{t}:=(t_{1},\ldots,t_{m})\in \mathbb{N}$ and
    $\pi_{i}:=p^{t_{i}}-1$ for $i\in Y_{0}.$
    Put
  $$\mathcal{M}:=
  \{x^{(q\varepsilon_{i})}D_{j}  \mid  0\leq q\leq \pi_{i},\ i,j\in Y_0\},$$
  $$\mathcal{N}:=\{x_i x_kD_{l}  \mid  i\in Y_0,\ k,l\in Y_1\},$$
    $$\mathcal{P}:=\{x_k x_lD_{i}  \mid  i\in Y_0,\ k,l\in Y_1\}.$$

   For $u,v\in \mathbb{B} $ with $u\cap v=\emptyset,$
    define $u+v$ to be $w\in \mathbb{\mathbb{ B}}$ such that
   and $ w = u \cup v.$ If in addition $\max u< \min v,$ then
   denote $u \dot{+}v=w$.

 \begin{proposition}

 $\mathcal{W}$ is generated by $\mathcal{M}\cup \mathcal{N}\cup \mathcal{P}.$
 \end{proposition}

   \noindent  \textit{Proof.}  Denote by $\mathcal{X}$
   the subalgebra of $\mathcal{W}$
   generated by $\mathcal{M}\cup \mathcal{N}\cup \mathcal{P}.$ We
    make the following preparatory  remarks.

    (i)  Assert that
   $$x^{(\alpha)}D_i\in \mathcal{X}\quad \mbox{for all}\
   \alpha \in \mathbb{A}\  \mbox{and all}\ i\in Y_0.\eqno(2.2.1)$$
   In fact, it can be easily proved by induction on $r$ that for
   $i\in Y_0, $
   $$x^{(\pi_1\varepsilon_1+\cdots+\pi_r\varepsilon_r)}D_i\in
   \mathcal{X}\quad \mbox{for all}\ r\leq m.$$
   Therefore, $x^{(\pi)}D_i\in
   \mathcal{X}. $ The assertion follows, since $D_j\in \mathcal{X}$ for all $j\in
   Y_0.$

   (ii) Assert that
   $$x^{u}D_i\in \mathcal{X}\quad \mbox{for}\ u\in   \mathbb{B},
   \ |u|   \mbox{  even}; \ i\in Y_0. \eqno(2.2.2)$$
   We are going to prove the assertion by induction on $|u|.$
   For $|u|=2$, it is clear that $x^uD_i\in\mathcal{ P}\subseteq \mathcal{X}.$
   Assume that $|u|>2.$
   Then write $u=v\dot{+}w,$
   where $v,w\in \mathbb{B}$ such that
   $|w|=2,|v|=|u|-2.$

   Note that $x^wD_i\in \mathcal{P}\subset  \mathcal{X} $ and
    $x^vD_i\in  \mathcal{X} $ by the inductive hypothesis.
   We have
  $$x^uD_i=[x^vD_i,[x^wD_i,x^{(2\varepsilon_i)}D_i]]\in \mathcal{X},$$
as asserted.

   (iii) Assert that
   $$x^{u}D_k\in \mathcal{X}\quad \mbox{for}\ u\in \mathbb{B},\
      |u|  \mbox{ odd};\
   k\in Y_1.\eqno(2.2.3)$$
   We may suppose that $|u|\geq 3.$
   Write $u=v+w $
   such that $k\not \in v$ and $|v|=|u|-1,|w|=1.$
   From (2.2.2), we have  $x^vD_1\in \mathcal{X}.$
   Noticing that $x_1x^wD_k\in \mathcal{N},$
   we obtain that
   $$x^uD_k=\lambda[x^vD_1,x_1x^wD_k]\in \mathcal{X} \quad\mbox{where}\
   \lambda=1\ \mbox{or}-1.$$
   Now we are ready to show  that  $\mathcal{X}=\mathcal{W}.$
First claim that
 $$x^{(\alpha)}x^uD_i\in \mathcal{X}
 \quad\mbox{for}\ u\in \mathbb{B}, |u|  \mbox{  even;}
 \ i\in Y_0. \eqno(2.2.4)$$
 Without loss of generality, one may assume  that $|u|\geq 2.$
 Take $k\in u.$ Using (2.2.2) we have $x^uD_i\in \mathcal{X}.$
 Since $x^{(\alpha)}x_kD_k=[x^{(\alpha)}D_1,x_1x_kD_k]\in \mathcal{X},$ we have
 $$x^{(\alpha)}x^uD_i=[x^{(\alpha)}x_kD_k,x^uD_i]\in \mathcal{X},$$
 as desired.

 It remains to prove that
 $$x^{(\alpha)}x^uD_k\in \mathcal{X} \quad\mbox{for}\ u\in \mathbb{B}, \
  |u|  \mbox{ odd;}
 \ k\in Y_1.\eqno(2.2.5)$$
 If $|u|=1,$
 then $x_1x^uD_k\in \mathcal{N}\subseteq \mathcal{X}.$
 It follows from (2.2.1) that $x^{(\alpha)}x^uD_k=[x^{(\alpha)}D_1,x_1x^uD_k]
 \in \mathcal{X}.$
 Assume that $|u|\geq 3.$
 Choose  $r\in  u \backslash k.$
 Since
 $x^{(\alpha)}x_rD_r=[x^{(\alpha)}D_1,x_1x_rD_r]\in \mathcal{X},$
it follows from (2.2.3) that
 $$x^{(\alpha)}x^uD_k=[x^{(\alpha)}x_rD_r,x^uD_k]\in \mathcal{X}.$$

 By (2.2.4) and (2.2.5),  we have $\mathcal{X}=\mathcal{W},$
 completing the proof.
\qed\newline

 We note that in the  three sets above only $\mathcal{M}$  contains
    elements of $\mathbb{Z}$-degree higher than 1. Clearly,
    $\mathcal{M}$ is contained in $W(m;\underline{t}),$ the
    generalized Witt  algebra (see [18]). As we shall see in the following
    sections, this allows us to compute efficiently.

  Recall the torus $\mathcal {T} _{\mathcal{S}}$ of
  $\mathcal{S} $ (see Section 2.1). Clearly,
 $$ \{x_{r}D_{r}-x_{s}D_{s}\mid  \tau(r)=\tau(s);r,s\in Y\}
   \cup \{x_{r}D_{r}+x_{s}D_{s}\mid \tau(r)\neq\tau(s);r,s\in Y\} $$
   is an $\mathbb{F}$-basis of $\mathcal{T}_{\mathcal{S}}$
   consisting of toral elements.

  We first give the following simple fact:

  \begin{lemma}
  $ \mathcal{S}_{0} $ is spanned by
   $
   \mathcal{T}_{\mathcal{S}}
   \cup  \{x_{r}D_{s}\mid  \tau(r)=\tau(s),r\neq s; r,s\in Y\}. $
   \end{lemma}

    \noindent \textit{Proof.}  It is straightforward.
   \qed\newline

  Put $$\mathcal{Q}:=\{D_{ij}(x^{(r\varepsilon_{j})})
  \mid i,j\in Y_0,r\in \mathbb{N}_0\};$$
  $$\mathcal{R}:=\{D_{il}(x^{(2\varepsilon_{i})}x_{k})
  \mid i\in Y_0,k,l \in Y_1\}\cup \{D_{ij}(x_ix^v)\mid
  i,j\in Y_0,v\in \mathbb{B}_2\}.$$

     We now consider the generator set of $\mathcal{S}.$

  \begin{proposition}
   $\mathcal{S} $ is generated by $\mathcal{Q}\cup \mathcal{R}\cup \mathcal{S}_0.$
  \end{proposition}

  \noindent \textit{Proof.}  Let $\mathcal{X}$ be the subalgebra
  generated by $\mathcal{Q}\cup \mathcal{R}\cup \mathcal{S}_0.$
  We proceed in several steps to show that
  $\mathcal{S}=\mathcal{X}.$

  (i) Let $i,j\in Y_0$
   with $i\neq j,$ and let $k,l\in Y_{1}.$
  Then
  $$D_{jl}(x^{(2\varepsilon_{i})}x_k)=
  [D_{ij}(x^{(3\varepsilon_{i})}),D_{il}(x_ix_k) ]
  +3\delta_{kl}D_{ij}(x^{(3\varepsilon_{i})})\in \mathcal{X}\eqno(2.2.6)$$
  and
  $$D_{il}(x^{(\varepsilon_{i}+\varepsilon_{j})}x_k)
  =[D_{il}(x^{(2\varepsilon_{i})}x_k),D_{ij}(x^{(2\varepsilon_{j})}) ]
  \in \mathcal{X}.\eqno(2.2.7)$$
  In addition,  for any $r\in  Y_0\setminus \{i,j\},$
  we have
  $D_{rj}(x^{(\varepsilon_{i}+2\varepsilon_{j})})
  =[D_{rj}(x^{(3\varepsilon_{j})}),D_{ji}(x^{(2\varepsilon_{i})}) ]
  \in \mathcal{X}. $ Therefore,
  $$D_{rl}(x^{(\varepsilon_{i}+\varepsilon_{j})}x_k)
  =[D_{ir}(x^{(2\varepsilon_{i})}),D_{jl}(x^{(2\varepsilon_{j})}x_k) ]
  -\delta_{kl}D_{rj}(x^{(\varepsilon_{i}+2\varepsilon_{j})})
  \in \mathcal{X}.\eqno(2.2.8)$$
  Combining (2.2.6)--(2.2.8) and the assumption that
  $D_{il}(x^{(2\varepsilon_{i})}x_{k})
   \in \mathcal{R}$ for all
  $i\in Y_0$ and  $ k,l \in Y_1,$
  we obtain that
  $$D_{rl}(x^{(\varepsilon_{i}+\varepsilon_{j})}x_k)\in \mathcal{X} \quad
  \mbox{for all}\ r,i,j\in Y_{0}\
  \mbox{and all}\ \ k,l\in Y_1.\eqno(2.2.9)
  $$

  (ii) Let $k,l\in Y_1 $
  with  $k\neq l.$ For all $r\in Y_1 $ and   all $i\in Y_0,$
 it follows that
 $$D_{kr}(x_ix_kx_l)=[D_{ik}(x^{(2\varepsilon_{i})}x_k),D_{ir}(x_ix_l)  ]
 -D_{ir}(x^{(2\varepsilon_{i})}x_l)
 +2\delta_{rl}D_{ik}(x^{(2\varepsilon_{i})}x_k)\in \mathcal{X}.\eqno(2.2.10) $$

 (iii) Let $u\in \mathbb {B}_3.$
  Write  $u=v\dotplus w$ with  $|v|=2,|w|=1.$ Given $i\in Y_0$ and $k\in Y_1,$
  take $j\in Y_0\setminus i $ and
   $l\in v \setminus k,$
  Remember  the assumption  $D_{ij}(x_jx^v)\in \mathcal{R}.$ We
  have

  $$D_{ik}(x^u)=-[D_{ij}(x_jx^v),D_{kl}(x^wx_l)]\in  \mathcal{X}.\eqno(2.2.11)$$

  Combining (2.2.9)--(2.2.11) and the assumption that  $D_{ij}(x_ix^v)\in \mathcal{R} $
 for all $i,j\in Y_0$ and  $v\in \mathbb{ B}_{2},$ we obtain that
 $$\mathcal{S}_1\subset \mathcal{X}.\eqno(2.2.12)$$

 (iv) Assert that $D_{ij}(x^{(\alpha)})\in  \mathcal{X} $ for all $i,j\in Y_0
 $  and $\alpha\in \mathbb{A}.$
 For the purpose, we first show that
  $D_{12}(x^{(\pi_{1}\varepsilon_{1}+\pi_{2}\varepsilon_{2})})\in  \mathcal{X}.$
 It is easy to verify the following equations:
 $$D_{12}(x^{(2\varepsilon_{1}+\varepsilon_{2})})=
 [D_{12}(x^{(3\varepsilon_{1})}),D_{12}(x^{(2\varepsilon_{2})})]\in  \mathcal{X},$$
 $$D_{12}(x_1x^{(\pi_2\varepsilon_{2})})=
 [D_{12}(x^{(\pi_2\varepsilon_{2})}),D_{12}(x^{(2\varepsilon_{1}+\varepsilon_{2})})]\in  \mathcal{X},$$
 $$D_{12}(x^{(\pi_1\varepsilon_{1}+(\pi_2-1)\varepsilon_{2})})=
 [D_{12}(x_1x^{(\pi_2\varepsilon_{2})}),D_{12}(x^{(\pi_1\varepsilon_{1})})]\in  \mathcal{X},$$
 and for $k\in Y_1,$
 $$D_{2k}(x^{(\pi_1\varepsilon_{1}+(\pi_2-1)\varepsilon_{2})}x_k)=
 -[D_{12}(x^{(\pi_1\varepsilon_{1}+(\pi_2-1)\varepsilon_{2})}),
 D_{1k}(x^{(2\varepsilon_{1})}x_k)]\in  \mathcal{X},$$
 $$D_{2k}(x^{((\pi_1-1)\varepsilon_{1}+(\pi_2-1)\varepsilon_{2})}x_k)=
 [D_{1},
 D_{2k}(x^{(\pi_1\varepsilon_{1}+(\pi_2-1)\varepsilon_{2})}x_k)]\in  \mathcal{X},$$
 $$D_{2k}(x^{((\pi_1-1)\varepsilon_{1}+\pi_2\varepsilon_{2})}x_k)=-\frac{1}{2}
 [D_{2k}(x^{((\pi_1-1)\varepsilon_{1}+(\pi_2-1)\varepsilon_{2})}x_k),
 D_{2k}(x^{(2\varepsilon_{2})}x_k)]\in  \mathcal{X}.$$
 Similarly, we also have
  $D_{1k}(x^{(\pi_1\varepsilon_{1}+(\pi_2-1)\varepsilon_{2})}x_k)\in  \mathcal{X}. $
 Therefore,
 \begin{eqnarray*}
  D_{12}(x^{(\pi_1\varepsilon_{1}+\pi_2\varepsilon_{2})})
   &=&
 [D_{1k}(x^{(2\varepsilon_{1})}x_k),
 D_{12}(x^{((\pi_1-1)\varepsilon_{1}+\pi_2\varepsilon_{2})})]\\
  &&+ D_{1k}(x^{(\pi_1\varepsilon_{1}+(\pi_2-1)\varepsilon_{2})}x_k)
 -2D_{2k}(x^{((\pi_1-1)\varepsilon_{1}+\pi_2\varepsilon_{2})}x_k)\in  \mathcal{X}.
\end{eqnarray*}

 Now, proceeding by induction on $q$ one may easily prove that
 $$D_{q,q+1}(x^{(\pi_1\varepsilon_{1}+\cdots+\pi_{q+1}\varepsilon_{q+1})})
 \in  \mathcal{X}  $$
 (cf. [23, Theorem 1]). As a result, $D_{m-1,m}(x^{(\pi)})\in \mathcal{X}.$
 It follows that  $D_{i j}(x^{(\pi)})\in \mathcal{X} $  for all $i,j\in Y_0.$
Hence the assertion holds.

 (v) Suppose that $|u|$
 is odd. We propose to prove that
 $$D_{i k}(x^{u})\in  \mathcal{X}\quad\mbox{for all}\
   i\in Y_0,k\in Y_1.\eqno(2.2.13)$$
 For $|u|=3,$
  from (2.2.12) it is easily seen that (2.2.13) holds. Use induction on
   $|u|\geq 3$ to prove (2.2.13).
  Assume  that  $|u|\geq 5 $
 and write $u=v\dotplus w $
 where $|v|=|u|-2$ and $|w|=2.$
 Let $l=\min  v.$
 By inductive hypothesis, $D_{l i}(x^{v})\in  \mathcal{X}$ for all $  i\in Y_0.$
 Then for $  i\in Y_0 $ and $ k\in Y_1,$
 $$D_{i k}(x^{u})=[D_{l i}(x^{v}),D_{i k}(x_ix_lx^{w})]
 =[D_{l i}(x^{v}),[D_{l i}(x_lx^{w}), D_{i k}(x^{(2\varepsilon_i)}x_l)]]
  \in  \mathcal{X}.$$
 This proves (2.2.13).

 (vi) Suppose that $|u| $ is even.
 Assert that
 $$D_{k l}(x^{u})\in  \mathcal{X} \quad\mbox{for all}\  k,l\in Y_1.\eqno(2.2.14)$$
We first note that if $v\in \mathbb{ B}$
 and $|v|$
 is odd, then
$$D_{i  k}(x_ix^{v})\in  \mathcal{X}\quad\mbox{for all}\ i\in Y_0,k\in Y_1.\eqno(2.2.15)$$
 In fact, this follows from (2.2.13) and the following equation
 $$D_{i k}(x_ix^{v})=
 [D_{q i}(x^{v} ),D_{i k}(x^{(2\varepsilon_i)}x_q) ]\in \mathcal{X}
 $$
 where $q=\min  v .$

 One may write $u=v\dotplus w$ where $|v|=|u|-1,|w|=1.$
 Using (2.2.15) we know that for $k,l\in Y_1,$
 $$D_{kl}(x^{u})=
 [D_{1 k}(x_1 x^{v}),D_{1 l}(x_1 x^w) ]-
 D_{1 l}(x_1D_k(x^{v})x^w)+D_{k 1}(x_1x^vD_l(x^{w}))
 \in  \mathcal{X}.$$
 Then (2.2.14) follows.

 (vii) Let us complete the proof of this proposition. For the purpose, first,
 we show that if $|u|$ is odd where $ u\in
 \mathbb{ B},$ then
 $$D_{i k}(x^{(\pi)}x^{u}) \in  \mathcal{X}
 \quad\mbox{for all}\ i\in Y_0, k\in Y_1 .\eqno(2.2.16)$$
For $|u|=1,$
 let $j\in Y_0\setminus i.$
 Then
 $$D_{i k}(x^{(\pi)}x^{u})=[ D_{i j}(x^{(\pi)}),D_{k j}(x^{(2\varepsilon_j)}x^u)  ]
 \in  \mathcal{X}.$$
 Suppose that $|u|\geq 3.$
 To prove (2.2.16), we note that
 $$D_{j k}(x^{(3\varepsilon_j)}x_l)=-\frac{1}{3}
 [D_{j k}(x^{(2\varepsilon_j)}x_l),D_{j l}(x^{(2\varepsilon_j)}x_l)]
 \in \mathcal{X}  \quad\mbox{for all}\  j\in  Y_0, k,l\in Y_1
   \ \mbox{with }\  k\neq l.\eqno(2.2.17)$$
 Given $k\in Y_1,$
 since $|u|\geq 3,$
 one may choose  $l\in  u \setminus k.$
 Without loss of generality, assume that $l=\min   u.$
 Then we obtain from (2.2.13) and (2.2.17) that
 $$D_{k j}(x^{(2\varepsilon_j)}x^u)=
 [D_{l j}(x^{u}),D_{j k}(x^{(3\varepsilon_j)}x_l)]
 \in \mathcal{X}   \quad\mbox{for all}\  j\in  Y_0, k\in Y_1.
   \eqno(2.2.18)$$
 Take $j\in Y_0\setminus i.$
 Then (iv) and (2.2.18) ensure  that
 $$D_{i k}(x^{(\pi)}x^u)=-[D_{i j}(x^{(\pi)}),D_{k j}(x^{(2\varepsilon_j)}x^u)]
 \in \mathcal{X},$$
 proving (2.2.16).

 Next, we show that for $u\in \mathbb{B},$
 $$D_{ij}(x^{(\pi)}x^u)\in \mathcal{X}
 \quad\mbox{for}\ i,j\in Y_0,\ |u|\ \mbox{even};\eqno(2.2.19)$$
$$D_{k l}(x^{(\pi)}x^u)\in \mathcal{X}
 \quad \mbox{for}\  k,l\in Y_1,\  |u|\ \mbox{even}.\eqno(2.2.20)$$
 One may write  $u=v\dotplus w$ where $ |v|=|u|-1,|w|=1.$
 Let $l:=\min v .$
 Using  (2.2.16), we have
 $$D_{i j}(x^{(\pi)}x^u)=-[D_{i l}(x^{(\pi)}x^v),D_{j i}(x_ix_lx^w)]
 \in \mathcal{X},\ |u|\ \mbox{even};$$
that is, (2.2.19) holds. Similarly,  we  have also
 $$D_{kl}(x^{(\pi)}x^u)=-[D_{ik}(x^{(\pi)}x^v),D_{i l}(x^{(2\varepsilon_i)}x^w)]
 \in \mathcal{X},\ |u|\ \mbox{even},$$
 proving (2.2.20).

  By (2.2.16), (2.2.19) and (2.2.20), it is easily seen that
  $\mathcal{X}=\mathcal{S }.$  The proof is complete.
  \qed

 \begin{remark}

 {\rm In contrast to the case of Lie superalgebra
  $S,$
   the elements in $\mathcal{R}$ cannot be generated by $\mathcal{Q}$ and
   $ \mathcal{S}_0 $ (cf. [23, Theorem 1]). The
   reason
   is that there are not odd   elements in $ \mathcal{S}_0.$ On
   the other hand , in contrast to the case $\mathcal{W},$ the
   elements of $\mathbb{Z}$-degree zero of $\mathcal{S}_0$ cannot be generated
   by $\mathcal{Q}\cup \mathcal{R} .$ }
   \end{remark}

 \section{Derivations of $\mathcal{W}$}

 \noindent In this section, we shall describe explicitly
 the derivation algebra of $\mathcal{W}.$ To do that, we proceed in two steps. First,
 we show in Section 3.1 that every derivation of nonnegative $\mathbb{Z}$-degree
 vanishing on the top of
$\mathcal{W}$ is necessarily inner. Next,  we show in Section 3.2 that
every derivation of nonnegative $\mathbb{Z}$-degree can be reduced to be vanishing on
the top and therefore, every derivation of nonnegative $\mathbb{Z}$-degree
is necessarily inner. Finally, using the generator set, we compute the derivations
of negative  $\mathbb{Z}$-degree and then
 formulate the derivation algebra $\mathrm{Der}(\mathcal{W}).$

\subsection{The top of $\mathcal{W}$}

 \noindent In this subsection, we discuss the influence of the top
 $\mathcal{W}_{-1}\oplus \mathcal{W}_{0}$ over the derivations of $\mathcal{W}.$
 We shall see that the homogeneous derivations of nonnegative
 $\mathbb{Z}$-degree of $\mathcal{W}$
 which vanish on the top $\mathcal{W}_{-1}\oplus \mathcal{W}_{0}$
 turn out to be inner. In other words, if   two homogenous derivations
  of nonnegative $\mathbb{Z}$-degree of $\mathcal{W}$ coincide on the top
 $\mathcal{W}_{-1}\oplus \mathcal{W}_{0}$ then they are congruent modulo
 $\mathrm{ad}\mathcal{W}.$  This observation allows us to focus attention to
 $\mathrm{Der}(\mathcal{W}_{-1}\oplus \mathcal{W}_{0},\mathcal{W})  $ in next subsection.

\begin{lemma}

Let $\phi\in \mathrm{Der}_t\mathcal{W},t\in \mathbb{Z}.$
 Suppose that $\phi(\mathcal{W}_{-1}\oplus \mathcal{W}_{0})=0$
 and $t\not\equiv 0 \ \pmod{ p}.$
 Then $\phi =0.$
 \end{lemma}

 \noindent \textit{Proof.}
 Consider the degree derivation $\Gamma:=\sum_{i\in Y}x_iD_i.$
 By Lemma 2.1.2(iv), we have
 $$[\Gamma,D]=rD\quad \mbox{for all} \ D\in \mathcal{W}_r,\ r\in
 \mathbb{Z}.
 $$
 Then noting that $\Gamma\in \mathcal{W}_0.$
 We have
 $$[\Gamma,\phi(D)]=r\phi(D).\eqno(3.1.1)$$
 Since $\mathrm{zd}(\phi(D))=r+t,$ again by Lemma 2.1.2(iv), we have
 $$[\Gamma,\phi(D)]=(r+t)\phi(D).\eqno(3.1.2)$$
 It follows from (3.1.1) and (3.1.2) that $t\phi(D)=0.$
 Since $t\not\equiv 0\pmod{p},$
 we have $\phi(D)=0$ for all
 $D\in \mathcal{W}_r,\ r\in
 \mathbb{Z}.$
 This proves that $\phi=0.$
\qed\newline

Recall that $\mathcal{P}=\{x_kx_lD_i\mid i\in Y_{0}, k,l \in
Y_{1}\}.$ We have the following

 \begin{lemma}

 Let $\phi\in \mathrm{Der}(\mathcal{W})$ with
 $\mathrm{zd}(\phi)=t\geq 0.$
 Suppose that $\phi(\mathcal{W}_{-1}\oplus \mathcal{W}_{0})=0.$
  Then $\phi$
 is a  scalar
 transformation of
 $\mathrm{span}_{\mathbb{F}}\mathcal{P}.$
 \end{lemma}

 \noindent \textit{Proof.}
 Since $\phi(\mathcal{W}_{-1})=\phi(\mathcal{W}_{0})=0,$
 by Lemma 2.1.1 we have

 $$\phi(x_{k}x_lD_i)\in \mathcal{G}_{t+1}\quad \mbox{for}\ i\in Y_{0}, k,l \in
Y_{1}. \eqno(3.1.3)$$
 In view of Lemma 3.1.1, we may assume that $t\equiv0 \pmod{p}.$
 Consider $\Gamma'=\sum_{r\in Y_1}x_rD_r\in \mathcal{W}_{0}.$
 Clearly, $[\Gamma',x_{k}x_lD_i]=2x_{k}x_lD_i.$
 Applying $\phi$
 to this equation, we obtain that
 $$[\Gamma',\phi(x_{k}x_lD_i)]=2\phi (x_{k}x_lD_i).\eqno(3.1.4)$$
 Note that $\mathrm{zd}(\phi(x_{k}x_lD_i))=t+1.$
 By (3.1.3) and Lemma 2.1.2(vi) we have
  $$[\Gamma',\phi(x_{k}x_lD_i)]
  =\left\{\begin{array}{cc}(t+2)\phi(x_{k}x_lD_i)
   ,&\quad\mathrm{zd}(\phi)\ \mbox{even;} \\
 (t+1)\phi(x_{k}x_lD_i),&\quad\mathrm{zd}(\phi)\ \mbox{odd}.\end{array}\right.$$
If $\mathrm{zd}(\phi)$ is odd, then we obtain from (3.1.4) and the
above equation that
  $\phi(x_{k}x_lD_i)=0,$ since $t\equiv 0\pmod{p}.$
That is, $\phi(\mathcal{P})=0.$
  Now consider the case that   $\mathrm{zd}(\phi)$
  is even. In view of (3.1.3), we may assume that
 $$\phi(x_{k}x_lD_i)=
 \sum_{r\in Y_0\atop |u|=t+2} c_{u,r}x^uD_r
 \quad\mbox{where}\ c_{u,r}\in \mathbb{F}.\eqno(3.1.5)$$
 Since $n>2,$
 choose $ s\in Y_1 \setminus \{k,l\}.$
 Then $[x_kx_lD_i,x_sD_s]=0.$
 Noting that $\phi(\mathcal{W}_0)=0,$
 we have $[\phi(x_kx_lD_i),x_sD_s]=0. $
 It follows that
 $$\Big[\sum_{r\in Y_0\atop|u|=t+2  }c _{u,r}x^uD_r,x_sD_s\Big]
 =0\quad \mbox{for all}\ s\in Y_1\setminus\{k,l\}.\eqno(3.1.6)$$
From (3.1.5) and (3.1.6), we have $c _{u,r}=0 $
  unless  $u=<k,l>.$
 Therefore we obtain from (3.1.5) that
 $$\phi(x_kx_lD_i)=\sum_{r\in Y_0}c_{<k,l>,r}x_kx_lD_r. \eqno(3.1.7)$$
If $\mathrm{zd}(\phi)>0,$ by (3.1.7)  we have $\phi(x_kx_lD_i)=0 $
and therefore $\phi(\mathcal{P})=0.$ Thus, it remains only  to
consider the case $\mathrm{zd}(\phi)=0.$
 Applying $\phi$
to the equation that $[x_kx_lD_i,x_iD_i]=x_kx_lD_i,$ we obtain
from (3.1.7) that
$$\phi(x_kx_lD_i)=c_{<k,l>,i}x_kx_lD_i.\eqno(3.1.8)$$
Let $\lambda :=c_{<k,l>,i}.$ For  arbitrary $ j\in Y_0,$ applying
$\phi$ to the equation that $[x_kx_lD_i,x_iD_j]=x_kx_lD_j,$ we
have
$$\phi(x_kx_lD_j)=\lambda(x_kx_lD_j).\eqno(3.1.9)$$
For arbitrary $r\in Y_1,$ applying $ \phi$ to the equation that
$[x_rD_k, x_kx_lD_j]=x_rx_lD_j,$ we obtain from (3.1.9) that
$$\phi(x_rx_lD_j)=\lambda x_rx_lD_j.\eqno(3.1.10) $$
(3.1.8)--(3.1.10) show that $\lambda=c_{<k,l>,i}$ is independent
of the choice of $k, l, i. $ Therefore,
 $$\phi(x_kx_lD_i)=\lambda x_kx_lD_i\quad\mbox{for all}\ k,l \in Y_1,\ i\in Y_0.$$
 The proof is complete.
\qed\newline

Recall that $\mathcal{N}=\{x_ix_kD_l\mid i\in Y_0, k,l\in Y_1\}.$
We have the following fact.

\begin{lemma}

 Let
 $\phi \in \mathrm{Der}(\mathcal{W}), \mathrm{zd}(\phi)=t\geq 0.$
 Suppose that $\phi(\mathcal{W}_{-1}\oplus \mathcal{W}_{0})=0. $
 Then $\phi(\mathcal{N})=0.$
  \end{lemma}

 \noindent  \textit{Proof.}  By Lemma 2.1.1,
 $\phi(x_ix_kD_l)\in \mathcal{G}_{t+1}  $ for $i\in Y_0, k,l\in
 Y_1.$
 In view of Lemma 3.1.1, it suffices  to consider the case that
 $t\equiv 0 \pmod{p}.$
 Recall $\Gamma':=\sum_{r\in Y_1}x_rD_r.$ Clearly,
 $$[\Gamma',x_ix_kD_l]=0\quad \mbox{for}\ i\in Y_0,k,l\in Y_1.$$
Therefore,
 $$[\Gamma',\phi (x_ix_kD_l)]=0\quad\mbox{for all }\ i\in Y_0,k,l\in Y_1.\eqno(3.1.11)$$
  Since $\phi(x_ix_kD_l)\in \mathcal{G}_{t+1}$, using Lemma 2.1.2 we
 obtain that
 $$[\Gamma',\phi (x_ix_kD_l)]= \left\{\begin{array}{cc}(t+2)
 \phi(x_{i}x_kD_l),&\mathrm{zd}(\phi)\ \mbox{even;} \\
 (t+1)\phi(x_{i}x_kD_l),&\mathrm{zd}(\phi)\ \mbox{odd}.\end{array}\right.\eqno(3.1.12)$$
Comparing (3.1.11) and (3.1.12) and noticing that $t\equiv
 0\pmod{p}$, we have
 $$\phi(x_ix_kD_l)=0\quad\mbox{for all}\ i\in Y_0,k,l\in Y_1.$$
The proof is complete.
 \qed\newline

Recall that $\mathcal{M}=\{x^{(q\varepsilon_i)}D_j\mid i,j\in
Y_{0}, 0\leq q\leq \pi_{i} \}.$ We have the following

\begin{lemma}

Let $\phi\in \mathrm{Der}(\mathcal{W}),
 \mathrm{zd}(\phi)=t\geq 0.$
 Suppose that $\phi(\mathcal{W}_{-1}\oplus \mathcal{W}_{0})=0.$
 Then $\phi(\mathcal{M})=0.$
\end{lemma}

 \noindent  \textit{Proof.}   We proceed by induction on $q\geq 2
 $ to prove that
$$\phi(x^{(q\varepsilon_i)}D_j)
=0\quad\mbox{for all}\ i,j\in Y_0,\ q\in \mathbb{N}_0.$$
 For $q=2,$
   Lemma 2.1.1 shows that $\phi(x^{(2\varepsilon_i)}D_j)\in \mathcal{G}_{t+1}.$
 Thus we may assume that
 $$\phi(x^{(2\varepsilon_i)}D_j)=
 \sum_{r\in Y}c_{u,r}x^uD_r \quad\mbox{where}\ c_{u,r}\in
 \mathbb{F}.\eqno(3.1.13)$$
 Note that $\mathrm{zd}(\phi(x^{(2\varepsilon_i)}D_j))\geq 1.$
For every pair $(u,r),$ there is
 $k\in  u \setminus r.$
 Clearly, $[x^{(2\varepsilon_i)}D_j,x_kD_k]=0.$
 Then
 $$[\phi (x^{(2\varepsilon_i)}D_j),x_kD_k]=0.\eqno(3.1.14)$$
 (3.1.13) and (3.1.14) imply that
 $c_{u,r}=0.$
 Since the pair $(u,r)$ is arbitrary,
 we have $\phi(x^{(2\varepsilon_i)}D_j)=0.$

 By the inductive hypothesis and Lemma 2.1.1, it is easily seen that
 $\phi(x^{(q\varepsilon_i)}D_j)\in \mathcal{G}.$
 Just as in the case $q=2,$
 one may check that $\phi(x^{(q\varepsilon_i)}D_j)=0.$
\qed\newline

Now we are able to prove the following fact.

 \begin{corollary}
 Let $\phi\in \mathrm{Der}_t(\mathcal{W}),t\geq 0.$
 Suppose that $\phi(\mathcal{W}_{-1}\oplus \mathcal{W}_{0})=0.$
Then $\phi \in \mathrm{ad}\mathcal{W}.$
 \end{corollary}

\noindent  \textit{Proof.}   By Lemma 3.1.2, there is $\lambda
\in \mathbb{F}$ such that
$$\phi(x_kx_lD_i)=\lambda x_kx_lD_i\quad\mbox{for all}\ k,l\in Y_1,i\in Y_0.$$
Put $\psi:=\phi-\frac{1}{2}\lambda\mathrm{ad}\Gamma'$ where
$\Gamma'=\sum_{r\in Y_1 }x_rD_r.$ Then
$\psi(\mathcal{W}_{-1}\oplus \mathcal{W}_{0})=0 $ and
 $\psi(\mathcal{P})=0.$ By Lemma 3.1.3 and
3.1.4, we have $\phi(\mathcal{N})=\phi(\mathcal{M})=0.$ By
Proposition 2.2.1, we obtain that $\psi=0;$ that is,
$\phi=\lambda\mathrm{ad}\Gamma'\in \mathrm{ad}\mathcal{W}.$ \qed

\begin{remark}
{\rm In view of Corollary 3.1.5, for determining the homogeneous
derivations of nonnegative $\mathbb{Z}$-degree of $\mathcal{W}$,
it suffices to reduce such a derivation to be vanishing on the
top. Just as in Proposition 2.1.6 and Corollary 3.1.5, the idea of
reduction will lead us throughout.}
\end{remark}

\subsection{Derivation algebra of $\mathcal{W}$}
\noindent In this subsection we shall determine the derivation
algebra
 of $\mathcal{W},$ the even part of the generalized Witt Lie superalgebra $W.$
 For the Lie superalgebra $W,$
 since the natural $\mathbb{Z}$-gradation is transitive, by a result analogous
  to Proposition 2.1.6 (see \cite[Lemma 5]{ZZ}), one may easily prove
 that the homogeneous superderivations of nonnegative $\mathbb{Z}$-degree
 of $ W $ are all inner.
 But, as mentioned above, the $\mathbb{Z}$-gradation of $\mathcal{W} $
   is not transitive. Thus we cannot obtain the
 corresponding  conclusion for $\mathcal{W}$ by using  Proposition 2.1.6
 directly.
 This observation leads us naturally  to  devote our
  attention  to the gradation component of $\mathbb{Z}$-degree zero. To do that, we  shall use
  a known result
   (see [15, Proposition 8.4, p. 193]).

 Let $V$ be a finite-dimensional vector space over $\mathbb{F}$. A
linear transformation $\varphi:V\rightarrow V$ is called
semisimple if the minimum polynomial of $\varphi$ has distinct
roots in some base field extension.

\begin{lemma} 
([15, Proposition 8.4, p. 193])
 Let $\frak{g}=\oplus_{i=-r}^{s}\frak{g}_{i}$ be a $\mathbb{Z}$-graded
 centerless Lie algebra and $\frak T\subset \frak{g}_{0}$   a torus of $ \frak{g}.$  If $\varphi\in \mathrm{Der}(\frak{g})$
 is homogeneous of $\mathbb{Z}$-degree $t$, then there is $e\in \frak{g}_t$ such that
 $(\varphi-\mathrm{ad}e)\mid _{\frak T}=0.$\qed
\end{lemma}

Recall that
$\mathcal{T}:=\mathrm{span}_{\mathbb{F}}\{\Gamma_{r}\mid r\in Y\}
$ is the canonical torus of $\mathcal{W}.$ Using the lemma above,
we can prove the following

 \begin{corollary} 

 Let $\phi\in
 \mathrm{Der}_{t}(\mathcal{W}), t\geq 0. $ Suppose that $\phi(\mathcal{W}_{-1})=0.$
 Then there exists $D\in \mathcal{G}_{t}$ such that
 $(\phi-\mathrm{ad}D)\mid _{\mathcal{T}}=0.$
\end{corollary}

 \noindent  \textit{Proof.}   By Lemma 2.1.2, all the  standard basis elements of
 $\mathcal{W}$ are  eigenvectors of $\mathrm{ad}\Gamma_{r}$ for each $r\in Y.$
  Note that every element of $\mathcal{T}$ is
 semisimple. By Lemma 3.2.1, there exists $E\in \mathcal{W}_{t}$
 such that$ (\phi-\mathrm{ad}E)\mid _{\mathcal{T}}=0.$ Note that
 $\phi(\mathcal{T})\subset \mathcal{G}_{t}.$ Then
 $$[E, \Gamma_{r}]=\phi(\Gamma_{r})\in \mathcal{G}_{t}\quad\mbox{for all}\
 r\in Y.$$
 By Lemma 2.1.2(vii), there is $D\in \mathcal{G}_{t}$ such that $[D,
 \Gamma_{r}]=[E,
 \Gamma_{r}] $ for all $r\in Y.$ Hence $(\phi-\mathrm{ad}D )\mid _{\mathcal{T}}=0.$
\qed\newline

 Now we can prove that every homogeneous derivation of nonnegative
 $\mathbb{Z}$-degree of
 $\mathcal{W}$ is inner. We first prove the following

 \begin{proposition}

  Let $t>0$. Then
  $\mathrm{Der}_t \mathcal{W}=\mathrm{ad}\mathcal{W}_t.$
\end{proposition}

 \noindent  \textit{Proof.}    It suffices to prove that
 $\mbox{Der}_t \mathcal{W}\subset\mathrm{ad}\mathcal{W}_t .$
 Let $\phi \in \mbox{Der}_t \mathcal{W},$
  $t>0.$
 By Proposition 2.1.6, we may assume that $\phi(\mathcal{W}_{-1})=0.$
Then Corollary  3.2.2 shows that there is $D\in \mathcal{G}_{t}$
such that $\psi( \mathcal{T} )=0, $
 where $\psi:=\phi-\mathrm{ad}D.$

 We treat separately the two  cases $t\geq 3$ and $0<t<3$
 in order to prove that
 $$\psi(x_rD_s)=0\quad
 \mbox{for }\ r,s\in Y\ \mbox{with}\
 \tau(r)=\tau(s).\eqno(3.2.1)$$

 \noindent \textit{Case} (i): \ $t\geq 3.$
 Pick $k\in Y_{1}\setminus \{r,s\}. $
 Then  $[\Gamma_k,x_rD_s]=0. $ Consequently,
  $$[\Gamma_k,\psi(x_rD_s)]=0 \quad \mbox{for all}\
   k \in Y_1\setminus\{r,s\}.\eqno(3.2.2)$$
 Note that $\psi(x_rD_s)\in \mathcal{G}$ by Lemma 2.1.1.
 Therefore,  by Lemma 2.1.2(vii), we obtain from (3.2.2)
 that $\mathrm{zd}(\psi(x_rD_s))\leq
 2.$ Since $\mathrm{zd}(\psi)=t\geq 3$ by our hypothesis, we
 obtain that $\psi(x_rD_s)=0;$ that is, (3.2.1) holds in this case.
\newline

 \noindent \textit{Case}(ii): $0<t<3.$ Then $t=1$ or $t=2.$
 Consider $\Gamma':=\sum_{r\in Y_1}\Gamma_r.$
 Obviously,
$[\Gamma',x_rD_s]=0.$
 Then
 $$[\Gamma',\psi(x_rD_s)]=0.\eqno(3.2.3)$$
On the other hand,  noticing that $\psi(x_rD_s)\in
\mathcal{G}$ (by Lemma 2.1.1),  we obtain by Lemma 2.1.2(vi) that
$$[\Gamma',\psi(x_rD_s)]=2\psi(x_rD_s).\eqno(3.2.4)$$
Since $\mbox{char}\mathbb{F}\neq 2, $ (3.2.3) and (3.2.4) imply
that $\psi(x_rD_s)=0$.

 So far, we have proved that $\psi(\mathcal{W}_{-1}+\mathcal{W}_0)=0.$
 Using Corollary 3.1.5, we have $\psi\in \mathrm{ad}\mathcal{W}_t$ and
 therefore,
 $\phi \in \mathrm{ad}\mathcal{W}_t.$
\qed\newline

 We  have also the following

 \begin{proposition}
  $\mathrm{Der}_0\mathcal{W}=\mathrm{ad}\mathcal{W}_0.$
  \end{proposition}

\noindent   \textit{Proof.}  Let $\phi \in  \mathrm{Der}
_0\mathcal{W}.$
 In view of Proposition 2.1.6 and Corollary  3.2.2, we may assume that
 $\phi(\mathcal{W}_{-1})=\phi(\mathcal{T})=0.$

 Let $i,j\in Y_0$
 with $i\neq j.$
 Noting that $\phi(x_iD_j)\in \mathcal{G}_0,$
 we have
 $$\phi(x_iD_j)=\phi([x_iD_i,x_iD_j])=[x_iD_i,\phi(x_iD_j)]=0.\eqno(3.2.5)$$
 We want to prove that
 $$\phi(x_kD_l)=c_{kl}x_kD_l
  \quad \mbox{for all }k,l\in Y_{1}\ \mbox{with}\ k\neq l, \eqno(3.2.6 )$$
  where $c_{kl}\in \mathbb{F}.$
Obviously, we may assume for fixed
$k, l\in Y_1$ with $k \neq l$ that
$$\phi(x_{k}D_{l})=\sum_{r,s\in Y_1}c_{rs}x_{r}D_{s}\quad
\mbox{where}\ c_{rs}\in \mathbb{F}.\eqno(3.2.7)$$ Then
$$\phi(x_{k}D_{l})=\phi([\Gamma_k, x_{k}D_{l}])=[\Gamma_k, \phi(x_{k}D_{l})].\eqno(3.2.8)$$
By (3.2.7) and (3.2.8), we can easily compute that
$$\sum_{r,s\in Y_1}c_{rs}x_{r}D_{s}=\sum_{s\in Y_1\setminus k}
c_{ks}x_{k}D_{s}+(c_{kk}x_{k}-\sum_{r\in Y_1}c_{rk}x_{r})D_k.
\eqno(3.2.9)$$ Comparing the coefficients of $D_k$ in (3.2.9), we
have
$$\sum_{r\in Y_1}c_{rk}x_{r}=c_{kk}x_{k}-\sum_{r\in Y_1}c_{rk}x_{r}
=-\sum_{r\in Y_1\setminus k}c_{rk}x_{r}.$$ Consequently, $
c_{kk}=0.$ Comparing the coefficients of $D_{s}$ for $s\in
Y_{1}\setminus k,$  one may get $c_{rr}=0$ for all $r\in
Y_{1}\setminus k.$ Thus we obtain from (3.2.7) that
$$\phi(x_{k}D_{l})=\sum_{r,s\in Y_1, r\neq s}c_{rs}x_{r}D_{s}\quad
\mbox{where}\ c_{rs}\in \mathbb{F}.\eqno(3.2.10)$$ For any $q\in
Y_1,$ we have $[\Gamma_q,
x_{k}D_{l}]=(\delta_{qk}-\delta_{ql})x_{k}D_{l}.$ Then
$$[\Gamma_q, \phi(x_{k}D_{l})]=(\delta_{qk}-\delta_{ql})\phi(x_{k}D_{l}).\eqno(3.2.11)$$
Using (3.2.10) we obtain from (3.2.11) that
$$\sum_{r,s\in Y_1, r\neq s}(\delta_{qr}-\delta_{qs})c_{rs}x_{r}D_{s}
=(\delta_{qk}-\delta_{ql})\sum_{r,s\in Y_1, r\neq
s}c_{rs}x_{r}D_{s}.
$$
It follows   that
$$(\delta_{qr}-\delta_{qs})c_{rs}=(\delta_{qk}-\delta_{ql})c_{rs}
\quad \mbox{for all}\ q\in Y_1.
$$
This implies that $c_{rs}=0$ unless $ (r,s)=(k,l). $ Hence (3.2.6)
holds.
 Applying $\phi$ to the following equation
 $$[x_rD_k,x_kD_l]=x_rD_l \quad \mbox{for}\ r,k,l\in Y_1, l\neq r,$$
 we obtain  that
 $$c_{rk}+c_{kl}=c_{rl}. \eqno(3.2.12)$$
 Similarly, we obtain from  the   equation
 $[x_lD_k,x_kD_l]=x_lD_l-x_kD_k$
  that
 $$c_{lk}+c_{kl}=0.\eqno(3.2.13)$$
 Clearly, the following  system of $n-1$ linear equations in $n$ unknowns
 $\lambda_{1'},\lambda_{2'},\ldots,\lambda_{n'}$ has solutions:
 $$\left.  \begin{array}{c}
    \lambda_{1'}-\lambda_{2'}=c_{1'2'}\\
  \lambda_{1'}-\lambda_{3'}=c_{1'3'}\\
  \cdots\ \cdots\ \cdots\ \cdots\\
  \lambda_{1'}-\lambda_{n'}=c_{1'n'}.\end{array}\right.$$
 Let $(\lambda_{1'},\lambda_{2'},\ldots,\lambda_{n'})$ be a solution.
Then
$\lambda_k-\lambda_l
=(\lambda_k-\lambda_1)+(\lambda_1-\lambda_l)=c_{k1}+c_{1l}=c_{kl}.$
Set
$$D:=\sum_{r\in Y_1}\lambda_r \Gamma_r,\quad \psi:=\phi-\mathrm{ad}D. $$
Then we have, for all $k,l\in Y_1$ with $k\neq l,$
\begin{eqnarray*}
\psi(x_kD_l)&=&\phi(x_kD_l)-[D,x_kD_l]\\
 &=&c_{kl}x_kD_l-(\lambda_k-\lambda_l)x_kD_l\\
 &=&0.
 \end{eqnarray*}
Using (3.2.5) we have also $\psi (x_iD_j)=0 $ for $i,j\in Y_0.$
Therefore, $\psi (\mathcal{W}_{0})=0.$ Clearly,
$\psi(\mathcal{W}_{-1})=0.$ Corollary 3.1.5 ensures that $\psi\in
\mathrm{ad}\mathcal{W} $ and therefore, $\phi\in
\mathrm{ad}\mathcal{W}_0.$ \qed\newline

Summarizing, we have the following

 \begin{proposition}

 The homogeneous derivations of nonnegative
 $\mathbb{Z}$-degree of $\mathcal{W}$ are all inner.
 \end{proposition}

\noindent  \textit{Proof.}   This is a direct consequence of
Propositions 3.2.3  and 3.2.4. \qed\newline

In view of Proposition 3.2.5, it remains only to determine the
homogeneous derivations of negative $\mathbb{Z}$-degree of
$\mathcal{W}.$ Our work is motivated by the corresponding results
and methods in [2, 23] and will depend heavily on the generator
set of $\mathcal{W}$ (Proposition 2.2.1). We first determine
$\mbox{Der}_{-1}\mathcal{W}.$ To do that, we need the following
lemma, which asserts that the derivations of $\mathbb{Z}$-degree
$-1$ are completely determined by $\mathcal{W}_{0}.$

 \begin{lemma} 

  Suppose that  $\varphi\in \mathrm{Der}_{-1}(\mathcal{W})$
 and $\varphi(\mathcal{W}_0)=0.$
 Then $\varphi=0.$
 \end{lemma}

 \noindent{\bf Proof.}  First, we  prove that $\varphi(\mathcal{N}\cup\mathcal{P})=0.$
 Observe that
 $$[x_kx_lD_i,x_iD_i]=x_kx_lD_i\quad\ \mbox{for all}\ k,l\in Y_1, i\in Y_0.\eqno(3.2.14)$$
 Note that $\varphi(\mathcal{W}_0)=0$
 and $\varphi (x_kx_lD_i)\in \mathcal{G}_0.$ Applying $\varphi$
to (3.2.14), one gets
 $$\varphi (x_kx_lD_i)=[\varphi (x_kx_lD_i),x_iD_i]=0.$$
 Similarly, one may check that $\varphi (x_ix_kD_l)=0 $ by applying $\varphi$
 to the equation that $[x_iD_i,x_ix_{k}D_l]=x_ix_kD_l.$
 Hence, $\varphi(\mathcal{N})=\varphi(\mathcal{P})=0.$

 Next, we  prove that $\varphi(\mathcal{M})=0.$ To do that we assert that
 $$\varphi(x^{(r\varepsilon_i)}D_i)=0 \quad \mbox{for all}\
   r\in \mathbb{N}, i\in Y_0.
 \eqno(3.2.15)$$
 For $r=1,$ it is clear that
  (3.2.15) holds.
 For $r=2,$ noticing that
 $\varphi(x^{(2\varepsilon_i)}D_i)\in \mathcal{G}_{0}\subset E(\mathcal{G}), $
 we obtain that
 $$\varphi(x^{(2\varepsilon_i)}D_i)=\varphi([x_iD_i, x^{(2\varepsilon_i)}D_i])
 =[x_iD_i,\varphi(x^{(2\varepsilon_i)}D_i)]=0.$$
 Now we proceed by induction on $r> 2 $ to prove (3.2.15).
 By inductive hypothesis, $\varphi(x^{(r\varepsilon_i)}D_i)\in \mathcal{G}_{r-2}.$
 Thus, we may assume that $\varphi(x^{(r\varepsilon_i)}D_i)=\sum_{u,q} c_{uq}x^uD_q
 $  where $c_{uq}\in \mathbb{F}.$ For any  given pair $(u,q)$,
 find $k\in u\setminus q$ since $r > 2.$
Then we may get $c_{uq}=0 $ from the  equation that
  $[x_kD_k,\varphi(x^{(r\varepsilon_i)}D_i)]=0. $
Thus,  (3.2.15) holds. Therefore,
  $$\varphi(x^{(r\varepsilon_i)}D_j)=\varphi( [x^{(r\varepsilon_i)}D_i,x_iD_j])
  =0 \quad\mbox{for all}\ j\in Y_0\setminus i.$$
Hence $\varphi(\mathcal{M})=0.$
 By Proposition 2.2.1, we conclude that $\varphi=0.$
\qed\newline

  Using Lemma 3.2.6, we can prove the following

 \begin{proposition}
  $\mathrm{Der}_{-1}(\mathcal{W})=\mathrm{ad}\mathcal{W}_{-1}.$
 \end{proposition}

 \noindent \textit{Proof. }  Let $\varphi \in \mbox{Der}_{-1}(\mathcal{W}).$
 For $i\in Y_0,k,l\in Y_1,$
 applying $\varphi$
 to the equation that  $[x_iD_i,x_kD_l]=0,$
one gets
 $$[\varphi(x_iD_i),x_kD_l]+[x_iD_i,\varphi(x_kD_l)]=0.$$
 Note that $[\varphi(x_iD_i),x_kD_l]=0,$ since $\varphi(x_iD_i)\in \mathcal{W}_{-1}.$
 It follows  that
 $$ [x_iD_i,\varphi(x_kD_l)]=0\quad \mbox{for}\ i\in Y_0, \ k,l \in Y_1.$$
 Since $\varphi(x_kD_l)\in \mathcal{W}_{-1},$
 this implies that $\varphi(x_kD_l)=0$ for all $k, l\in Y_1.$

 Given $i\in Y_0,$
 suppose that $\varphi(x_iD_i)=\sum_{r\in Y_0}c_{ir}D_r $ where $c_{ir}\in\mathbb{ F}.$
 Applying $\varphi$
 to the equation that $ [x_iD_i,x_jD_j]=0$ for $i,j\in Y_0$ with $i\neq j,$
 one may get $c_{ij}=0 $
 whenever $j\neq i.$
 Therefore, $\varphi(x_iD_i)=c_{ii}D_i.$
 Applying $\varphi$
 to the equation that $ [x_iD_i,x_iD_j]=x_iD_j,$
one gets
 $$c_{ii}D_j+[x_iD_i,\varphi(x_iD_j)]=\varphi(x_iD_j).$$
This implies that $\varphi(x_iD_j)=c_{ii}D_j $ for all $j\in Y_0,$
 since $\varphi(x_iD_j)\in \mathcal{W}_{-1}.$
Put $\psi :=\varphi -\sum_{r\in Y_0} c_{rr}(\mathrm{ad}D_r).$
 Then $\psi(\mathcal{W}_{0})=0.$
 By Lemma 3.2.6, we have $\psi=0; $ that is,
 $\varphi =\sum_{r\in Y_0} c_{rr}(\mathrm{ad}D_r)\in \mathrm{ad}
  \mathcal{W}_{-1}.$
\qed\newline

 To determine the derivations of $\mathbb{Z}$-degree less than
 $-1$, we first establish a technical lemma.

 \begin{lemma}

  Let $\phi \in \mathrm{Der}_{-q}(\mathcal{W}), q>1.$
 If $\phi (x^{(q\varepsilon_i)}D_i)=0$ for all $ i\in Y_0,$
 then $\phi=0.$
\end{lemma}

 \noindent  \textit{Proof.}  We first show that $\phi(\mathcal{M})=0.$
 If $t\leq q,$
  then $\phi (x^{(t\varepsilon_i)}D_i)=0$  for all $i\in Y_0.$
 Suppose that $t>q.$ We use induction on
  $t$
 to prove that $\phi (x^{(t\varepsilon_i)}D_i)=0.$
 By inductive hypothesis and Lemma 2.1.1,  $\phi (x^{(t\varepsilon_i)}D_i)
 \in \mathcal{G}_{t-q-1}.$
 Thus we may assume that
 $$\phi (x^{(t\varepsilon_i)}D_i)
 =\sum_{u\in \mathbb{ B}, r\in Y}c_{u,r} x^uD_r
 \quad\mbox{where}\
  c_{u,r} \in \mathbb{F}.\eqno(3.2.16)$$
 If $|u|\geq 2,$
 then for the pair $(u,r),$ there is $k\in  u \setminus r.$
 Note that
  $[x_kD_k,x^{(t\varepsilon_i)}D_i]=0.$
 Applying $\phi$
 to this equation  and using  (3.2.16), one may get $c_{u,r}=0.$
 Therefore, noting that $\mathrm{zd}(\phi (x^{(t\varepsilon_i)}D_i))=t-q-1
 \geq 0,$
 we have
 $$\phi (x^{(t\varepsilon_i)}D_i)=\sum_{u\in \mathbb{ B}_{1}, r\in Y_{1}}c_{u,r}
 x^uD_r=
 \sum_{l,r\in Y_1}c_{l,r}x_lD_r, \eqno(3.2.17)$$
 where $ c_{l,r}:=c_{u,r}$ if $u=\langle l\rangle.$ Applying $\phi$
 to the equation that $[x^{(t\varepsilon_i)}D_i,x_lD_r]=0$ for $ l,r\in Y_1,$
 we obtain that $c_{l,r}=0$ whenever $l\neq r$
  and that  $c_{l,l}=c_{r,r}$  for all $ l,r\in Y_1.$ Denote
  $\lambda:=c_{r,r}$ for all $r\in Y_{1}.$
It follows from (3.2.17) that
 $$\phi (x^{(t\varepsilon_i)}D_i)=
 \sum_{l\in Y_1}\lambda x_lD_l=\lambda \Gamma'.$$
Applying $\phi$ to the equation that
$[x^{(t\varepsilon_i)}D_i,x_kx_lD_j]=0 $ for $i,j\in Y_0$ with
$j\neq i $ and $k,l\in Y_1, $  we have
$$2\lambda x_kx_lD_j+[x^{(t\varepsilon_i)}D_i,\phi(x_kx_lD_j)]=0.  $$
Since $\mathrm{zd}(\phi(x_kx_lD_j))\leq -1,$  the equation above
yields $\lambda =0.$ For $j\in Y\setminus i,$ it follows that
$\phi
(x^{(t\varepsilon_i)}D_j)=\phi([x^{(t\varepsilon_i)}D_i,x_iD_j])=0.$
Therefore, $\phi(\mathcal{M})=0.$

It remains to show that $\phi(\mathcal{N}\cup\mathcal{P})=0.$
Since $\mathcal{N}\cup\mathcal{P}\subseteq \mathcal{W}_{1}$ and
$\mathrm{zd}(\phi)=-q\leq -2,$ it suffices to consider  only the
case $\mathrm{zd}(\phi)=-2.$  Let $k,l\in Y_1,i\in Y_0.$ Note that
$\phi (x_kx_lD_i)\in \mathcal{W}_{-1}.$  Applying $\phi $ to the
equation $[\Gamma',x_kx_lD_i]=2x_kx_lD_i,$ we have
$$\phi (x_kx_lD_i)=\frac{1}{2}\phi([\Gamma',x_kx_lD_i])
=\frac{1}{2}[\Gamma',\phi(x_kx_lD_i)]\in[\Gamma',
\mathcal{W}_{-1}]=0,
$$
proving that $\phi (\mathcal{P})=0.$ Similarly, we have
$$\phi (x_ix_kD_l)=\phi ([\Gamma'',x_ix_kD_l])= [\Gamma'',\phi(x_ix_kD_l)]
=- \phi(x_ix_kD_l).$$
Thus $\phi (x_ix_kD_l)=0,$ since $\mathrm{char} \mathbb{F}\neq 2.$

So far, we have showed that
$\phi(\mathcal{M})=\phi(\mathcal{N})=\phi(\mathcal{P})=0.$ By
Proposition 2.2.1, $\phi=0.$ \qed\newline

We are in the position to compute the homogeneous derivations of
negative $\mathbb{Z}$-degree. We treat two cases separately.
First, we give the following

\begin{proposition}

Suppose that $q>1 $ is not any $p$-power. Then
$\mathrm{Der}_{-q}(\mathcal{W})=0.$
\end{proposition}

 \noindent \textit{Proof.}   Let $\phi \in \mbox{Der}_{-q}(\mathcal{W}).$
 Clearly, $\phi(\mathcal{W}_{-1}+\mathcal{W}_{0})=0.$
 If $q \not \equiv  0\pmod{p},$
 then Lemma 3.1.1 shows that $\phi=0.$
 Thus we assume  that $q\equiv 0\pmod{p}.$
 Write $q$ to be the
  $p$-adic expression:
 $q=\sum^r_{k=1}c_kp^k,$
 where $0\leq c_k<p$ and $ c_r\neq 0.$
 We note that
 $
 {q\choose p^r} \not\equiv 0\pmod{p}$ and $\binom{q }{p^r-1}
 \equiv 0\pmod{p}.$
 According to Lemma 3.2.8, it suffices to show that
 $$\phi(x^{(q\varepsilon_i)}D_i)=0 \quad\mbox{for all}\ i\in Y_0.$$
Direct computation shows that
 $$
 \big[x^{((q-p^r+1)\varepsilon _i)}D_i,x^{(p^r\varepsilon
 _i)}D_i\big]
   = \left(  {q\choose p^r-1}
 -
  {q\choose p^r}
 \right)x^{(q\varepsilon _i)}D_i
= - {q\choose p^r}
 x^{(q\varepsilon _i)}D_i.
 $$
 Since $q-p^r+1<q$ and $p^r<q,$ from the equation above one easily gets
 $\phi(x^{(q\varepsilon_i)}D_i)=0.$
 The proof is complete.
\qed\newline

For the remaining case we have the following
 \begin{proposition} 

  Let $q$ be  $p$-power $p^r
 $ for some $r\in \mathbb{N}.$
  Then
$$\mathrm{Der}_{-q}(\mathcal{W})=
\mathrm{span}_{\mathbb{F}}\{(\mathrm{ad}D_i)^q| i\in Y_0\}.$$
\end{proposition}

 \noindent \textit{Proof.}  By Leibniz rule, it is easily seen
 that $(\mathrm{ad}D_i)^q $ is a derivation of $\mathcal{W} $  for all $i\in
 Y_{0}.$ It follows that one implication holds.
 To prove the converse implication, let $\phi\in \mbox{Der}_{-q}(\mathcal{W}).$
 One may  assume that
 $$\phi(x^{(q\varepsilon_i)}D_i)=\sum_{r\in Y_0}a_{ir}D_r \quad
 \mbox{where}\  a_{ir}\in \mathbb{F}.$$
 Applying $\phi$
 to the equation that
 $[x^{(q\varepsilon_i)}D_i,x_jD_j ]=0$  for   $j\in Y_0\backslash i,$
 yields that  $a_{ij}=0.$
 Therefore, $\phi(x^{(q\varepsilon_i)}D_i)=a_{ii}D_i$
 where $i\in Y_0.$
 Set $\psi :=\phi -\sum_{r\in Y_0}a_{rr}(\mathrm{ad}D_r)^q.$
 Then $\psi(x^{(q\varepsilon_i)}D_i)=0$ for all $i\in Y_0. $
 By Lemma 3.2.8 we have $\psi =0;$
 that is, $\phi=\sum_{r\in Y_0}a_{rr}(\mathrm{ad}D_r)^q.$
 The proof is complete.
\qed\newline

Assembling the main results obtained in this subsection we are
able to describe the derivation algebra of $\mathcal{W}.$ Recall
that $\mathcal{W}$ stands for the even part of
$W=W(m,n;\underline{t}),$ where
$\underline{t}=(t_{1},t_{2},\ldots, t_{m})\in \mathbb{N}^{m}.$

\begin{theorem} 

$\mathrm{Der} (\mathcal{W}) =\mathrm{ad}\mathcal{ W}\oplus\{(
\mathrm{ad} D_i)^{p^{r_{i}}}\mid i\in Y_{0}, 1\leq r_{i}<t_{i}
\}.$
\end{theorem}

\noindent\textit{Proof.} By Leibniz rule, it is sufficient to show
the inclusion`$ \subset $'. Note that when $r_{i}\geq t_{i},$ the
derivation $( \mathrm{ad} D_i)^{p^{r_{i}}}$ vanishes for $i\in
Y_{0}.$ Then the theorem follows from Propositions 3.2.5, 3.2.7,
3.2.9 and 3.2.10. \qed

\begin{remark} 

{\rm Now we can answer the question for $\mathcal{W}$ mentioned in
the introduction. By Theorem 3.2.11 and [23, Theorem 1], it is
easily seen that $ \mathrm{Der}
(W_{\overline{0}})=\mathrm{Der}_{\overline{0}}(W),$  where $W$
stands for the Lie superalgebra $W(m,n;\underline{t}).$ As
mentioned in Introduction, this implies that every derivation of
the Lie algebra $W_{\overline{0}}$ can extend  to be a
superderivation of the Lie superalgebra $W.$}
\end{remark}

\section{The derivation algebra of $\mathcal{S}$}

\noindent  In this section we shall determine the derivation
algebra  of $\mathcal S $ ($\mathcal S $ is the even part of the
Lie superalgebra of $S(m,n;\underline {t})$). In contrast to the
case of $\mathcal W,$
 we do not find a way to reduce directly a homogeneous derivation of nonnegative
  $\mathbb{Z}$-degree of $\mathcal S$
  to be
 a derivation vanishing on $ \mathcal {S}_{-1}\oplus\mathcal {S}_{0}.$
 Thus we encounter the difficulty that Lemma 3.2.1 is by no means  applicable  for
$\mathcal S.$  This observation forces us  to establish a
proposition similar to Proposition 2.1.6, but, where   $\mathcal
W_{-1}$ is replaced with the canonical torus $\mathcal
{T}_{\mathcal{S}}$ contained in $\mathcal{S}_{0}.$ This is why we
generalize Proposition 8.2 in [15, p. 192] to be Lemma 2.1.3 and
then give Corollary 2.1.5. As a result, we  establish  such a
proposition indeed, which is written to be Lemmas 4.2.4 and 4.2.5
in Section 4.2.

  Recall our convention that $m, n>2.$

  \subsection{ The top of $\mathcal{S}$}

\noindent  Just as in the case of $\mathcal W,$ in this subsection we
shall
   study the derivations of nonnegative $\mathbb Z$-degree of $\mathcal
   S$ to $\mathcal W,$
  which vanish on the top of $\mathcal S.$ As the final result in
  this subsection, it is proved that such a derivation must be
  inner  and
  determined by a scalar multiple  of $\Gamma'=\sum_{r\in Y_{1}}x_{r}D_{r}.$ The arguments will
  be based on the generator set of $\mathcal S$ (see Proposition
  2.2.3).

  First, we consider the generators of the form
  $D_{il}(x^{(2\varepsilon_{i})}x_{k})$ in $\mathcal{R},$ where $i\in Y_0$ and $  k,l\in Y_1.$

  \begin{lemma}
    Suppose that  $\phi \in \mathrm{Der} (\mathcal{S},\mathcal{W})$
    with
  $\mathrm{zd}(\phi)\geq 0$
  and that $\phi(\mathcal{S}_{-1}+\mathcal{S}_{0})=0.$
  Then
  $$\phi (D_{il}(x^{(2\varepsilon_{i})}x_{k}))=0 \quad\mbox{for all }
   \  i\in Y_0, \ k,l\in Y_1.\eqno(4.1.1)$$
   \end{lemma}

 \noindent   \textit{Proof.}  By the definition of $D_{il},$
  $$D_{il}(x^{(2\varepsilon_{i})}x_{k})=x_{i}x_{k}D_{l}
  +\delta_{kl}x^{(2\varepsilon_{i})}D_i
   \quad \mbox{for all}\ i\in Y_0, k,l\in Y_1.  \eqno(4.1.2)$$
  By our assumption
  and Lemma 2.1.1, one may assume that
  $$\phi (D_{il}(x^{(2\varepsilon_{i})}x_{k}))=
  \sum_{r\in Y, u\in \mathbb{ B}}c_{u,r}x^uD_r \quad \mbox{where}
  \ c_{u,r}\in \mathbb{F}. \eqno(4.1.3)$$
  We proceed in two steps to prove the equation (4.1.1).
\newline

 \noindent  \textit{Case} (i): $\mathrm{zd}(\phi)$
  is even. Then it follows from (4.1.3) that $c_{u,r}=0$
  for all
  $ \ r\in Y_1. $
  Thus
  $$\phi (D_{il}(x^{(2\varepsilon_{i})}x_{k}))
  =\sum_{r\in Y_0, u\in \mathbb{ B}}c_{u,r}x^uD_r. \eqno(4.1.4)$$

  First, consider  the case $k\neq l.$
  Then (4.1.2) yields
  $$D_{il}(x^{(2\varepsilon_{i})}x_{k})=x_{i}x_{k}D_{l}\quad\ \mbox{for}\
   k\neq l. \eqno(4.1.5)$$
  Clearly,
  $$[\Gamma_k-\Gamma_l,D_{il}(x^{(2\varepsilon_{i})}x_{k})]=
  2D_{il}(x^{(2\varepsilon_{i})}x_{k})\quad \mbox{ for } k \neq l. \eqno(4.1.6) $$
  Consequently,
  $$[\Gamma_k-\Gamma_l,\phi (D_{il}(x^{(2\varepsilon_{i})}x_{k}))]=
  2\phi (D_{il}(x^{(2\varepsilon_{i})}x_{k})) \quad \mbox{ for }\
   k\neq l. \eqno(4.1.7) $$
 Then we obtain from (4.1.4) and (4.1.7) that
  $$\sum_{r\in Y_0, u\in \mathbb{ B}}(\delta_{k\in u}-\delta_{l\in u})c_{u,r}x^uD_r
  =2 \sum_{r\in Y_0, u\in \mathbb{ B}}c_{u,r}x^uD_r. \eqno(4.1.8)$$
  It follows immediately from (4.1.8) that
  $$(\delta_{k\in u}-\delta_{l\in u})c_{u,r}=2c_{u,r}
  \quad \mbox{for} \ r\in Y_0,u\in \mathbb{ B}.$$
  This implies that $c_{u,r}=0$ for all $r\in Y_0, u\in \mathbb{ B}.$
  Therefore, (4.1.1) holds for the case $k\neq l.$

 Second,  we consider the case $k=l.$ By (4.1.2),
  $$D_{ik}(x^{(2\varepsilon_{i})}x_{k})
  =x_ix_kD_k+x^{(2\varepsilon_{i})}D_i.\eqno(4.1.9)$$
  For $q_1,q_2 \in Y_1,$ it is easy to see that
   $[\Gamma_{q_1}-\Gamma_{q_2},D_{ik}(x^{(2\varepsilon_{i})}x_{k})]=0.$
  Applying $\phi$, one gets
  $$ [\Gamma_{q_{1}}-\Gamma_{q_2},\phi (D_{ik}(x^{(2\varepsilon_{i})}x_{k}))]=0.
  \eqno(4.1.10)$$
  We then obtain from (4.1.4) and (4.1.10) that
  $$\sum_{r\in Y_0, u\in \mathbb{ B}}(\delta_{q_1\in u}-\delta_{q_2\in
  u})c_{u,r}x^uD_r=0
  \quad \mbox{for all} \ q_1,q_2\in Y_1,r\in Y_0.$$
 Consequently,
  $$(\delta_{q_1\in u}-\delta_{q_2\in u})c_{u,r}=0 \quad
   \mbox{for all}\ q_1,q_2\in Y_1,r\in Y_0. \eqno(4.1.11)$$
  Note that (4.1.11) implies that, if $u\neq  \omega $
  then $c_{u,r}=0 $ for all $r\in Y_0.$
  Therefore, it follows from (4.1.4) that
  $$\phi (D_{ik}(x^{(2\varepsilon_{i})}x_{k}))
  =\sum_{r\in Y_0}c_{\omega,r}x^{\omega}D_{r}.\eqno(4.1.12)$$
  Find $j\in Y_0 \setminus i, $ since $|Y_0|\geq 3.$
  Note that (Notice (4.1.9))
  $$[\Gamma_{i}-\Gamma_{j},D_{ik}(x^{(2\varepsilon_{i})}x_{k})]
  =D_{ik}(x^{(2\varepsilon_{i})}x_{k}).$$
    Applying $\phi $
   and substituting (4.1.12), we obtain that
   $$[\Gamma_{i}-\Gamma_{j},\sum_{r\in Y_0}c_{\omega,r}x^{\omega}D_{r}]
   =\sum_{r\in Y_0}c_{\omega,r}x^{\omega}D_{r}.$$
  It follows immediately that
  $$-c_{\omega,i}x^{\omega}D_{i}+c_{\omega,j}x^{\omega}D_{j}=\sum_{r\in Y_0}c_{\omega,r}x^{\omega}D_{r}.$$
  This implies that
  $$c_{\omega,r}=0,\quad \mbox{whenever}\ r\neq j.$$
  Then (4.1.12) yields that
  $$\phi (D_{ik}(x^{(2\varepsilon_{i})}x_{k}))=
  c_{\omega,j}x^{\omega}D_{j}
  \quad \mbox{whenever}\ j\in Y_0\setminus i.\eqno(4.1.13)$$
  The general assumption $m\geq 3 $
  and (4.1.13) show that (4.1.1) holds in the case $k=l.$
  \newline

  \noindent\textit{Case} (ii): $\mathrm{zd}(\phi)$
  is odd. Note that in this case it is easily seen in (4.1.3) that
   $c_{u,r}=0$ for all $r\in
  Y_0.$
  Thus
  $$\phi (D_{il}(x^{(2\varepsilon_{i})}x_{k}))
  =\sum_{r\in Y_1, u\in \mathbb{ B}}c_{u,r}x^uD_r.\eqno(4.1.14)
  $$
  Recall $\Gamma':=\sum_{r\in Y_1}x_{r}D_{r}.$
  Using the formulas (4.1.5) and (4.1.9)式, one obtains by direct
  computation that
  $$[\Gamma'+n\Gamma_i,D_{il}(x^{(2\varepsilon_{i})}x_{k})]
  =nD_{il}(x^{(2\varepsilon_{i})}x_{k}),\eqno(4.1.15)$$

  $$[\Gamma'+(n-1)\Gamma_i+\Gamma_j,D_{il}(x^{(2\varepsilon_{i})}x_{k})]
  =(n-1)D_{il}(x^{(2\varepsilon_{i})}x_{k}) \quad\mbox{for}\ j\in Y_0
  \setminus i.\eqno(4.1.16)$$
  Applying $\phi$
  to (4.1.15) and (4.1.16), respectively,  one
  gets from (4.1.14) that
  $$[\Gamma',\phi(D_{il}(x^{(2\varepsilon_{i})}x_{k}))]
  =[\Gamma'+n\Gamma_i,\phi(D_{il}(x^{(2\varepsilon_{i})}x_{k}))]
  =n\phi(D_{il}(x^{(2\varepsilon_{i})}x_{k})),\eqno(4.1.17)$$
   and
    $$[\Gamma',\phi(D_{il}(x^{(2\varepsilon_{i})}x_{k}))]=
    [\Gamma'+(n-1)\Gamma_i+\Gamma_j,\phi(D_{il}(x^{(2\varepsilon_{i})}x_{k}))]
    =(n-1)\phi(D_{il}(x^{(2\varepsilon_{i})}x_{k})).\eqno(4.1.18)$$
  Comparing (4.1.17) and (4.1.18), we have
   $\phi(D_{il}(x^{(2\varepsilon_{i})}x_{k}))=0.$
 Summarizing,  (4.1.1) holds.
  \qed\newline

  Next, we consider the generators of the form
  $D_{ij}(x_{i}x^v)$ in $\mathcal{R},$ where $i,j\in Y_0$ and  $v\in \mathbb{B}_{2}.$

  \begin{lemma}

   Let $\phi \in \mathrm{Der}_{t}
  (\mathcal{S},\mathcal{W})$ with
  $t\geq 0.$
  Suppose that $\phi(\mathcal{S}_{-1}+\mathcal{S}_{0})=0.$
  Then there is $\lambda \in\mathbb{F} $ such that
  $\phi (D_{ij}(x_{i}x^v))=\lambda D_{ij}(x_{i}x^v) $
  for all $i,j\in Y_0$ and  $v\in \mathbb{B}_{2}.$
 \end{lemma}

  \noindent\textit{ Proof. } Note that
  $$\phi(D_{ij}(x_{i}x^v)) =\phi([\Gamma_q-\Gamma_j, D_{ij}(x_{i}x^v) ])
  =[\Gamma_q-\Gamma_j, \phi(D_{ij}(x_{i}x^v))]\quad \mbox{for all}
  \ q\in Y_0\setminus j.\eqno(4.1.19)$$
  Assume that
  $$\phi(D_{ij}(x_{i}x^v))=\sum_{u\in \mathbb{B},r\in Y}c_{u,r}x^uD_r
  \quad\mbox{where}\ c_{u,r}\in \mathbb{F}.\eqno(4.1.20)$$

  \noindent\textit{Case }(i):  $\mathrm{zd}(\phi)$
  is even. Then
  $$\phi(D_{ij}(x_{i}x^v))=\sum_{u\in \mathbb{B},r\in Y_0}c_{u,r}x^uD_r.\eqno(4.1.21)$$
  If follows from (4.1.19) and (4.1.21) that
  $$\phi(D_{ij}(x_{i}x^v))=\sum_{u}c_{u,j}x^uD_j-\sum_{u}c_{u,q}x^uD_q\quad \mbox{for all}
  \ q\in Y_0\setminus j. \eqno(4.1.22)$$
  This implies that $c_{u,q}=0$ whenever $q\in Y_0\setminus j $ and therefore,
  $$\phi(D_{ij}(x_{i}x^v))=\sum_{u}c_{u,j}x^uD_j.\eqno(4.1.23)$$
   Find $k\in Y_1\setminus v, $ since
   $n\geq 3.$ Then
  $[\Gamma_q+\Gamma_k,D_{ij}(x_{i}x^v)]=0 $ for $q\in Y_0\setminus j.$
  Applying $\phi$
  to this equation and then substituting (4.1.23), one gets
  $$\sum_{u\in \mathbb{B} }\delta _{k\in u}c_{u,j}x^uD_j=0.$$
  For any fixed $u\neq v $ with $|u|\geq 2,$
  if $k\in u\setminus v,$
  then the equation above forces  $c_{u,j}=0.$
  It follows from (4.1.23) that
  $$\phi(D_{ij}(x_{i}x^v))=c_{v,j}x^vD_j.$$
  Note that $D_{ij}(x_ix^v)=x^vD_j.$
   Just as in the proof of Lemma 3.1.2, one may easily
  check that $\lambda:=c_{vj} $ is independent of  the choice of $v$ and $j.$
  Therefore, the lemma holds in this case.
  \newline

  \noindent\textit{Case} (ii): $\mathrm{zd}(\phi)$
  is odd.  Then
  $$\phi(D_{ij}(x_{i}x^v))=\sum_{u\in \mathbb{B},r\in Y_1}c_{u,r}x^uD_r.$$
  Find $q\in Y_0\setminus j.$ Then
  $$[\Gamma_q-\Gamma_j, D_{ij}(x_{i}x^v)]=D_{ij}(x_{i}x^v).$$
  Hence,
  $$\phi(D_{ij}(x_{i}x^v))=[\Gamma_q-\Gamma_j, \phi(x^vD_j)]=
  [\Gamma_q-\Gamma_j,  \sum_{u\in \mathbb{B},r\in Y_1} c_{u,r}x^uD_r]=0.$$
  The proof is complete.
  \qed\newline

  Finally, we consider the generators in $\mathcal{Q}.$

  \begin{lemma}

  Let
  $\phi \in \mathrm{Der}(\mathcal{S},\mathcal{W}), \mathrm{zd}(\phi)\geq 0.$
  Suppose that
  $\phi(\mathcal{S}_{-1}+\mathcal{S}_{0})=0.$
  Then
  $$\phi(D_{ij}(x^{(a\varepsilon_i)}))=0
  \quad \mbox{for all } \ i,j\in Y_0\ \mbox{and all}\ a\in \mathbb{N}.\eqno(4.1.24)$$
  \end{lemma}

  \noindent  \textit{Proof. }
  We proceed by induction on $a\geq 2$
  to prove (4.1.24).  Assume that (4.1.24) holds for $a-1.$
  Then Lemma 2.1.1 ensures that
  $$\phi(D_{ij}(x^{(a\varepsilon_i)}))=
  \sum_{u\in \mathbb{B},r\in Y} c_{u,r}x^uD_r
  \quad \mbox{where}\ c_{u,r}\in\mathbb{ F}.\eqno(4.1.25)$$
  Note that $D_{ij}(x^{(a\varepsilon_i)})=x^{((a-1)\varepsilon_i)}D_j.$
  Direct computation shows that
  $$[x_kD_k+x_jD_j,D_{ij}(x^{(a\varepsilon_i)})]=-D_{ij}(x^{(a\varepsilon_i)})
  \quad\mbox{for all}\ k\in Y_1.$$
 Applying $\phi$
  to the equation and using (4.1.25), one obtains that
  $$\sum_{u,r}(\delta_{k\in u}-\delta_{kr}) c_{u,r}x^uD_r
  -\sum_{u}c_{u,j}x^uD_j=- \sum_{u,r}c_{u,r}x^uD_r.\eqno(4.1.26)$$
  Comparing coefficients in (4.1.26), we have
  $$(\delta_{k\in u}-\delta_{kr}) c_{u,r}=-c_{u,r}
  \quad \mbox{for} \ r\neq j;\eqno(4.1.27)$$
  $$(\delta_{k\in u}-1) c_{u,j}=-c_{u,j} .\eqno(4.1.28)$$
  For arbitrary $u\in \mathbb{B}$ with $|u|\geq 2,$ find $k\in u$.
  Then (4.1.27) and (4.1.28) imply that $c_{ur}=0$ for $r\in Y,$
  since $\mbox{char}\mathbb{F} \neq 2.$
  This proves (4.1.24).
  \qed\newline

  Now we can conclude this subsection by the following main
  result.

 \begin{corollary} 

Suppose that
 $\phi \in \mathrm{Der}(\mathcal{S},\mathcal{W})$ with $\mathrm{zd}(\phi)\geq 0 $
  and  that $\phi(\mathcal{S}_{-1}+\mathcal{S}_{0})=0.$
 Then $\phi \in  \mathbb{F}\cdot \mathrm{ad}\Gamma'.$ In particular,  $\phi$ is inner.
 \end{corollary}

  \noindent \textit{ Proof.}  According to Lemma 4.1.2, there is $\lambda \in \mathbb{F}$
  such that
  $$\phi(D_{ij}(x_ix^v))=\lambda D_{ij}(x_ix^v)
  \quad \mbox{for all}\ i,j \in Y_0\ \mbox{and all }\ v\in \mathbb{B}_2.$$
  Set $\psi:=\phi -\frac{1}{2}\lambda \mathrm{ad}\Gamma'.$
   We then use Lemmas 4.1.1 and 4.1.3 in order to see that
  $\psi(\mathcal{Q}\cup \mathcal{R})= 0.$ Lemma 2.2.2 applies and
   $\psi=0.$
  Therefore, $\phi=\frac{1}{2}\lambda \mathrm{ad}\Gamma',$ completing the
  proof.

  \qed

  \begin{remark} 
  {\rm As noted in Remark 3.1.6, the central work in the
  sequel is to reduce the homogeneous derivation of nonnegative
  $\mathbb{Z}$-degree to be vanishing on the top of $\mathcal{S}.$
  In contrast to the case $\mathcal{W},$ just as remarked  in the
  introduction of this paper, we shall encounter the phenomenon that
  the reduction proposition [15, Proposition 8.4. p. 193] is not
  applicable in this case.}
  \end{remark}

  \subsection{Derivations of nonnegative $\mathbb{Z}$-degree}

  \noindent In this subsection, we shall determine
  $\mathrm{Der}_{t}(\mathcal{S}, \mathcal{W})$ and
  $\mathrm{Der}_{t}(\mathcal{S}) $ for $t\geq 0$ (Propositions 4.2.9 and 4.2.10).
  For  derivations of positive odd $\mathbb{Z}$-degree in
  $\mathrm{Der}(\mathcal{S}, \mathcal{W}),$ we establish a lemma
  (Lemma 4.2.4) analogous to Proposition 2.1.6, which contends that
  such a derivation vanishing on $\mathcal{S} _{-1}$ can be
  reduced to be vanishing on the toral elements $\Gamma _{k}-\Gamma
  _{1'},$ $k\in Y_{1}\setminus 1' $ (Lemma 4.2.4). For a derivation of
  nonnegative even $\mathbb{Z}$-degree in
  $\mathrm{Der}(\mathcal{S}, \mathcal{W}),$ we establish a
  corresponding lemma (Lemma 4.2.5). The reason that we treat
  those two cases separately is that the images of  elements in
  the canonical torus under a derivation are of different forms with respect to the
  standard $\mathbb{F}$-basis of $\mathcal{W}.$ In view of the results mentioned above,
  the remaining
  works in this subsection will be
  devoted to reducing the derivations
  of nonnegative $\mathbb{Z}$-degree to be vanishing on the top
  $\mathcal{S}_{-1}\oplus \mathcal{S}_{0}$. Thus we may conclude
  this subsection by using the result obtained in Section 4.1 (Corollary
  4.1.4).

 We first give two technical lemmas which will simplify our
 discussion.  Recall the notation
 $n=|\omega|=|Y_1|$.

 \begin{lemma}
  Suppose that $\phi \in \mathrm{Der}_{t}(\mathcal{S},\mathcal{W})$
  and $\phi(\mathcal{W}_{-1})=0.$

   $\mathrm{(i)}$ If $t=n-1$ is odd, then
    $$\phi(\Gamma_{1'}-\Gamma_{k})=0\quad \mbox{for all}\ k\in Y_1\setminus 1'.$$

    $\mathrm{(ii)}$ If $t=n-1$ is even, then there is $\lambda\in \mathbb{F}$
   such that
   $$(\phi -\lambda \mathrm{ad}(x^{\omega} D_{1'}))(\Gamma_{1'}-\Gamma_k)=
   0\quad \mbox{for all}\  k\in Y_1\setminus 1'.$$

   $\mathrm{(iii)}$ If $t>n-1,$ then $\phi=0.$
   \end{lemma}

   \noindent  \textit{Proof. } (i) We may assume that
   $$\phi(\Gamma_{1'}-\Gamma_{k})=\sum_{r\in Y_0} c_{r}x^{\omega} D_r
   \quad \mbox{where} \ c_r\in \mathbb{F}.\eqno(4.2.1)$$
   Applying $\phi$
   to $[\Gamma_1-\Gamma_i,\Gamma_{1'}-\Gamma_k]=0 $  for $i\in Y_0, k\in Y_1,$ we have
  $$[\Gamma_1-\Gamma_i,\phi(\Gamma_{1'}-\Gamma_k)]=
  [\Gamma_{1'}-\Gamma_k,\phi(\Gamma_1-\Gamma_i)]=0, \eqno(4.2.2)$$
  since $\phi(\Gamma_1-\Gamma_i)\in \mathrm{span}\{x^{\omega}D_r \ | \ r\in Y_0\}.$
  (4.2.1) and (4.2.2) then yield
   $$-c_1x^{\omega}D_1+c_ix^{\omega}D_i=0,\quad  i\in Y_0\setminus 1.$$
   Therefore $c_r=0 $ for all $r\in Y_0 $ and (i) holds.

   (ii)
   Applying $\phi$
   to the equation that $[\Gamma_{1'}-\Gamma_k,\Gamma_{1'}-\Gamma_l]=0,$ one gets
   $$[\Gamma_{1'}-\Gamma_k,\phi(\Gamma_{1'}-\Gamma_l)]
   =[\Gamma_{1'}-\Gamma_l,\phi(\Gamma_{1'}-\Gamma_k)]
   \quad \mbox{for all}\ k,l\in Y_1\setminus 1'.\eqno(4.2.3)$$
   Since $t$
   is even and $\phi(\mathcal{W}_{-1})=0,$
   by virtue of Lemma 2.1.1, we may assume that
   $$\phi(\Gamma_{1'}-\Gamma_k)=\sum_{r\in Y_1} c^{(k)}_{r}x^{\omega} D_r
   \quad \mbox{where }\  c_{ r}^{(k)} \in \mathbb{F}.\eqno(4.2.4) $$
   Combining (4.2.3) and (4.2.4), we have
   $$-c^{(l)}_{1'}x^{\omega} D_{1'}+c^{(l)}_{k}x^{\omega}D_{k}=
   -c^{(k)}_{1'}x^{\omega} D_{1'}+c^{(k)}_{l}x^{\omega} D_{l}.$$
   Consequently, $c^{(l)}_{k}=c^{(k)}_{l}=0,$ $c^{(l)}_{1'}=c^{(k)}_{1'}.$
   Let $\lambda:=c^{(2')}_{1'}=c^{(3')}_{1'}=\cdots=c^{(n')}_{1'}.$
   Then (4.2.4) shows that $\phi(\Gamma_{1'}-\Gamma_k)=\lambda x^{\omega} D_{1'} $
   for all
   $k\in Y_1 \setminus 1'.$
   Direct calculation shows that $\lambda\in \mathbb{F}$ is the
   desired scalar.

    (iii)
   By Lemma 2.1.1, $\phi(S_0)\subseteq \mathcal{G}.$ Note that
   $\mathcal{G} \subseteq  \bigoplus^{n-1}_{i=-1}\mathcal{W}_i. $ Since $t>
   n-1,$
   $$\phi(S_0)\subseteq \mathcal{W}_{t}\cap \mathcal{G}=0.$$
   Now, one may easily show  by induction on $r$ that
   $\phi(S_r)=0$ for all $r\in \mathbb{N}, $ proving
    $\phi=0.$
   \qed

   \begin{lemma} 
    Suppose that
   $\phi \in \mathrm{Der}(\mathcal{S},\mathcal{W})$
   and $\mathrm{zd}(\phi)>0$ is odd.

   $\mathrm{(i)}$ If $\mathrm{zd}(\phi)<n-1$ and
   $$\phi(\Gamma_{1'}-\Gamma_{2'})=\phi(\Gamma_{1'}-\Gamma_{3'})=\cdots
   =\phi(\Gamma_{1'}-\Gamma_{n'})=0 $$ then
   $$\phi(\Gamma_{1}-\Gamma_{2})=\phi(\Gamma_{1}-\Gamma_{3})=\cdots
   =\phi(\Gamma_{1}-\Gamma_{m})=0;\quad \phi(\Gamma_{1}+\Gamma_{1'})=0. $$

   $\mathrm{(ii)}$ If $\mathrm{zd}(\phi)=n-1$  then there are
   $\lambda_1,\ldots,\lambda_m\in \mathbb{ F}  $ such that
   $$\Big(\phi-\mathrm{ad}\Big(\sum _{i\in Y_{0}}\lambda_ix^{\omega} D_i\Big)\Big)
   (\Gamma_1-\Gamma_j)=0
   \quad \mbox{for all} \ j\in Y_0\setminus 1; $$
   $$\Big(\phi-\mathrm{ad}\Big(\sum _{i\in Y_0}\lambda_ix^{\omega} D_i\Big)\Big)
   (\Gamma_1+\Gamma_{1'})=0.
   \eqno$$
   \end{lemma}

  \noindent  \textit{Proof.}  Note that $\mathrm{zd}(\phi)$ is odd and
   $\phi(\mathcal{}S_{0})\subseteq \mathcal{G},$ by Lemma 2.1.1.
  We may assume that
  $$\phi(\Gamma_1-\Gamma_i)=\sum_{r\in Y_0,u\in \mathbb B}c^{(i)}_{u,r}x^uD_r
  \quad \mbox{where}\ c_{u,r}^{(i)}\in \mathbb{F}.\eqno(4.2.5)$$

  (i) Since $\mathrm{zd}(\phi)<n-1,$
  the coefficients $c^{(i)}_{\omega,r} $ in  (4.2.5)  vanish for all $r\in
  Y_{0}.$
  For any fixed $v\neq \omega,$
  find $k\in v$
  and $l\in Y_1\setminus v.$
  Applying $\phi $
  to the equation $[\Gamma_k-\Gamma_l,\Gamma_1-\Gamma_i]=0,$
  one may obtain by the hypothesis that
  $$[\Gamma_k-\Gamma_l,\phi(\Gamma_1-\Gamma_i)]=0
   \quad \mbox{for} \ i\in Y_0\setminus 1.
   $$
  Combining (4.2.5) with  the equation above, we have
  $$\sum_{r\in Y_0, u\in \mathbb B}(\delta_{k\in u}
  -\delta_{l\in u})c^{(i)}_{u,r}x^uD_r=0.\eqno(4.2.6)$$
  Since $\delta_{k\in v}-\delta_{l\in v}=1,$
   (4.2.6) implies that
  $$c^{(i)}_{v,r}=0\quad \mbox{for all } \ r\in Y_0.$$
  Then by (4.2.5),  $\phi (\Gamma_1-\Gamma_i)=0$ for $i\in Y_0\setminus 1.$
  One may show in the same way that $\phi (\Gamma_1+\Gamma_{1'})=0.$

  (ii) As $\mathrm{zd}(\phi)=n-1,$
  we know from (4.2.5) that
  $$\phi (\Gamma_1-\Gamma_i)
  =\sum_{r\in Y_0}c^{(i)}_{\omega,r}x^{\omega} D_r .\eqno(4.2.7)$$
  Applying $\phi$
  to the equation $[\Gamma_1-\Gamma_i, \Gamma_1-\Gamma_j]=0 $
  for $i,j\in Y_0\setminus 1 $ with $i\neq j,$
   and then combining with (4.2.7), we have
   $$\big[\sum_{r\in Y_0}c^{(i)}_{\omega,r}x^{\omega}D_r,\Gamma_1-\Gamma_j\big]=
  \big[ \sum_{r\in Y_0}c^{(j)}_{\omega,r}x^{\omega}D_r,\Gamma_1-\Gamma_i\big].$$
   Consequently,
   $$c^{(i)}_{\omega,1}x^{\omega} D_1-c^{(i)}_{\omega,j}x^{\omega} D_j=
   c^{(j)}_{\omega,1}x^{\omega} D_1-c^{(j)}_{\omega,i}x^{\omega} D_i.$$
   This implies that
   $$c^{(i)}_{\omega,1}=c^{(j)}_{\omega,1},\eqno(4.2.8)$$
   $$c^{(i)}_{\omega,j}=c^{(j)}_{\omega,i}=0
   \quad \mbox{for all }\ i,j\in Y_0\setminus 1\  \mbox{with }\ i\neq j.\eqno(4.2.9)$$
   By (4.2.8), let $\lambda_1:=c^{(i)}_{\omega,1} $ for $i\in Y_{0}\setminus 1.$
   Then (4.2.7) and (4.2.9) show that
   $$\phi(\Gamma_1-\Gamma_i)=\lambda_1x^{\omega} D_1
   +c^{(i)}_{\omega,i}x^{\omega} D_i.  \eqno(4.2.10)$$
   Let $\lambda_i:=-c^{(i)}_{\omega,i}$ for $ i=2,\ldots,m $ and
   $\varphi:=\phi-
   \sum _{i\in Y_{0}} \lambda_{i}\mathrm{ad}(x^{\omega} D_i). $
   Then (4.2.10) shows that
   $$\varphi(\Gamma_1-\Gamma_i)=0 \quad  \mbox{for all}
   \ i\in Y_0\setminus 1;\eqno(4.2.11)$$
   that is, the first equation in (ii)  holds.
   To show that $\varphi(\Gamma_1+\Gamma_{1'})=0,$
   assume that
   $$\varphi(\Gamma_1+\Gamma_{1'})
   =\sum_{r\in Y_0} d_{\omega,r}x^{\omega} D_r \quad\mbox{where}\  d_{\omega,r}
   \in \mathbb{F}.\eqno(4.2.12)$$
   For $i\in Y_0\setminus 1,$
   applying $\varphi$
   to the equation $[\Gamma_1-\Gamma_i,\Gamma_1+\Gamma_{1'}]=0,$
   and using (4.2.11) and (4.2.12), one may get
   $$d_{\omega,i}x^{\omega} D_i-d_{\omega,1}x^{\omega} D_1
   =[\Gamma_1-\Gamma_i,\sum_{r\in Y_0} d_{\omega,r}x^{\omega} D_r]=0.$$
   Consequently, $d_{\omega,i}=d_{\omega,1}=0$ for $ i\in Y_0\setminus 1. $
   This proves that
    $\varphi(\Gamma_1+\Gamma_{1'})=0.$
   Hence, the second equation in (ii) holds . The proof is complete.  \qed\newline

    Recall the canonical  torus of $\mathcal S$
    $$\mathcal{T}_{\mathcal{S}}:=\mbox{span}_{\mathbb{F}}
   \{\Gamma_1-\Gamma_2,\ldots,\Gamma_1-\Gamma_m,
   \Gamma_1+\Gamma_{1'} , \Gamma_{1'}-\Gamma_{2'},\ldots,\Gamma_{1'}-\Gamma_{n'}\}.$$
Summarizing, we have the following fact:

   \begin{corollary}

    Suppose that
   $\phi\in \mathrm{Der}(\mathcal{S},\mathcal{W})$ is homogeneous derivation
   of positive odd $\mathbb{Z}$-degree such that $\phi(\mathcal{S}_{-1})=0$ and
   $\phi(\Gamma_{k}-\Gamma_{1'})=0$ for all $k\in Y_1\setminus 1'.$ Then there is
   $E\in \mathcal{G}$
   such that
   $(\phi-\mathrm{ad}E) $ vanishes on the canonical torus $\mathcal T_{\mathcal S}.$
   \end{corollary}

   \noindent  \textit{Proof.}   This is a direct consequence of Lemma 4.2.1(iii)
    and Lemma 4.2.2. \qed\newline

    We now prove  two key lemmas in this subsection. First,
    consider the  derivations of  odd $\mathbb{Z}$-degree. We shall use Lemma 2.1.3.

   \begin{lemma} 

    Suppose that
   $\phi \in \mathrm{Der}_{t}(\mathcal{S},\mathcal{W}),$
   where $t>0$ is odd.  If $\phi(\mathcal{S}_{-1})=0,$
   then there is $D\in \mathcal{G}_t $
  such that
   $$(\phi - \mathrm{ad}D)(\Gamma_k-\Gamma_{1'})=0
   \quad \mbox{for all}\  k\in Y_1\setminus 1'.$$
   \end{lemma}

      \noindent  \textit{ Proof.}   If $t\geq n-1,$
      then Lemma  4.2.1(i) and (iii) show that
      $$\phi(\Gamma_k-\Gamma_{1'})=0 \quad \mbox{for all}\ k\in Y_1\setminus 1'.$$
      We therefore assume  that $t<n-1.$

      By Lemma 2.1.1, one may assume that
      $$\phi(\Gamma_k-\Gamma_{1'})=\sum_{r\in Y_0}f_{rk} D_r
       \quad\mbox{where}\ k\in Y_1\setminus 1';\ f_{rk}\in \Lambda(n).
       \eqno(4.2.13)$$
      Applying $\phi$ to the equation that
       $[\Gamma_k-\Gamma_{1'},\Gamma_l-\Gamma_{1'}]=0 $
      for $k,l\in Y_1\setminus 1',$
      and using (4.2.13), one gets
      $$[\sum_{r\in Y_0}f_{rk} D_r  ,\Gamma_l-\Gamma_{1'}]
      +[\Gamma_k-\Gamma_{1'},\sum_{r\in Y_0}f_{rl} D_r ]=0.$$
      Consequently,
      $$\sum_{r\in Y_0}(\Gamma_k-\Gamma_{1'})(f_{rl}) D_r=
      \sum_{r\in Y_0}(\Gamma_l-\Gamma_{1'})(f_{rk}) D_r. $$
      Since $\{D_{r}\mid r\in Y\} $ is a free basis of
      $\frak{A}$-module $ W ,$ we have
      $$(\Gamma_k-\Gamma_{1'})(f_{rl})= (\Gamma_l-\Gamma_{1'})(f_{rk})
       \quad \mbox{for all }
      \  r\in Y_0; \ k,l\in Y_1\setminus 1'.\eqno(4.2.14)$$
      Note that $f_{rk}\in \Lambda(n).$ Let
      $$f_{rk}=\sum_{|u|=t+1}c_{u,r,k}x^u\quad \mbox{where}\
       c_{u,r,k}\in \mathbb{F}.\eqno(4.2.15)$$
      Combining (4.2.15) and (4.2.14), we have
      $$\sum_{|u|=t+1}(\delta_{k\in u}-\delta_{1'\in u})c_{u,r,l}x^u=
      \sum_{|u|=t+1}(\delta_{l\in u}-\delta_{1'\in u})c_{u,r,k}x^u.$$
      Since $\{x^u\mid u\in \mathbb{B}\}$ is an $\mathbb{F}$-basis
      of $\Lambda(n),$ it follows
       that
      $$(\delta_{k\in u}-\delta_{1'\in u})c_{u,r,l}=
      (\delta_{l\in u}-\delta_{1'\in u})c_{u,r,k}
      \quad \mbox{for}\ r\in Y_0;\ k,l\in Y_1\setminus 1'.\eqno(4.2.16)$$

    Assume that $c_{u,r,k}$
    is an arbitrary fixed nonzero coefficient in (4.2.15), where $|u|=t+1<n,$ $r\in Y_0$ and
    $ k\in Y_1\setminus 1'.$
    If $1'\not\in u,$
    noting that $2\leq |u|,$
    one may find $l\in u \setminus k.$
    Then (4.2.16) shows that $\delta_{k\in u}=1;$
    that is, $k\in u.$
    If $1'\in  u, $
    noting that $|u|\leq n-1, $
    one may find $l\in Y_1\setminus u.$
    Then (4.2.16) shows that  $\delta_{k\in u}=0;$
    that is, $k\not\in u.$
    Summarizing, for any nonzero coefficient $c_{u,r,k}$ in (4.2.15),
    we have $\delta_{k\in u}+ \delta_{1'\in u}=1.$
    Then, we can rewrite (4.2.15) as follows
    $$f_{rk}=\sum_{ 1'\in u,k\not\in u} c_{u,r,k}x^u+
     \sum_{ 1'\not\in u, k\in u} c_{u,r,k}x^u.$$
    Direct computation shows that
    $$(\Gamma_k-\Gamma_{1'})(f_{rk})=-\sum_{1'\in u,k\not\in u} c_{u,r,k}x^u+
     \sum_{1'\not\in u, k\in u} c_{u,r,k}x^u.$$
   Therefore,
     $$(\Gamma_k-\Gamma_{1'})^2(f_{rk})=f_{rk}\quad\mbox{for all}\ r\in Y_0\
     \mbox{and}\ k\in Y_1\setminus 1'.
     \eqno(4.2.17)$$

      Put $V:=\Lambda(n), v_k:= f_{rk}$  and
    $ A_k = B_k:=\Gamma_k-\Gamma_{1'} $ for
     $ k\in Y_1\setminus 1'.$  Now, we are
    going to show that the conditions of Lemma 2.1.3 are fulfilled.
     (i) is automatic. (4.2.14) ensures that (ii) holds. (4.2.17) shows that
     (2.1.2) holds. (2.1.3) is clear. We check (2.1.1); that is,
    $$(\Gamma_k-\Gamma_{1'})^3
    =\Gamma_k-\Gamma_{1'}\quad \mbox{for} \ k \in Y_1\setminus 1'.\eqno(4.2.18)$$
    For every basis element $x^u$ of $V=\Lambda(n),$
    noticing that $\delta_{k\in u}-\delta_{1'\in u}=0,1$
    or $-1, $
    we obtain that
    $$(\Gamma_k-\Gamma_{1'})^3(x^u)=(\delta_{k\in u}-\delta_{1'\in u})^3x^u
    =(\delta_{k\in u}-\delta_{1'\in u})x^u= (\Gamma_k-\Gamma_{1'})(x^u).$$
    Thus (4.2.18) holds. By Lemma 2.1.3, there is $f_r\in \Lambda(n)$ for any fixed
     $r\in Y_0,$
    such that
    $$(\Gamma_k-\Gamma_{1'})(f_r)
    =f_{rk} \quad \mbox{for all} \ k \in Y_1\setminus 1'.\eqno(4.2.19)$$

    Set $D'=-\sum_{r\in Y_0}f_rD_r.$
    We  obtain by using (4.2.13) and (4.2.19) that,  for $k \in Y_1\setminus 1',$
    $$[D',\Gamma_k-\Gamma_{1'}] =
    \sum_{r\in Y_0}(\Gamma_k-\Gamma_{1'})(f_r)D_r=\sum_{r\in Y_0}f_{rk}D_r=
    \phi(\Gamma_k-\Gamma_{1'}).\eqno(4.2.20)
    $$
   Let $D:=D'_{t}.$ Then
    $D\in \mathcal{G}_t $ and
    (4.2.20) shows that
    $[D,\Gamma_k-\Gamma_{1'}]=\phi(\Gamma_k-\Gamma_{1'}) $
    for all $ k \in Y_1\setminus 1',$
   since $\mathrm{zd}(\phi)=t.$
   Hence
   $$(\phi-\mathrm{ad}D)(\Gamma_k-\Gamma_{1'})=0
   \quad \mbox{for all}\ k \in Y_1\setminus 1'.$$
   \qed\newline

   Let us consider the   case of even $\mathbb{Z}$-degree. Notice that, in
   contrast to the case of odd degree, we consider different elements
   in the canonical torus. This will  simplify the computation.

    \begin{lemma}
    Let $\phi \in \mathrm{Der}_{t}(\mathcal{S},\mathcal{W})
   $ where $t\geq 0$ is even. If $\phi(\mathcal{W}_{-1})=0,$
  then there is $D\in \mathcal{G}_t $
   such that
   $$(\phi - \mathrm{ad}D)(\Gamma_1+\Gamma_k)=0
   \quad \mbox{for all}\  k\in Y_1.\eqno(4.2.21)$$
  \end{lemma}

     \noindent  \textit{ Proof.}   Since $\phi$
      is of even $\mathbb{Z}$-degree, by Lemma 2.1.1, we may assume that for $k\in Y_1,$
      $$\phi(\Gamma_1+\Gamma_k)=\sum_{r\in Y_1}f_{rk}D_{r}\quad
      \mbox{where}\ f_{rk}\in \Lambda(n).\eqno(4.2.22)$$
      It is easily seen that
      $$[\Gamma_1+\Gamma_k,\phi(\Gamma_1+\Gamma_l)]=
      [\Gamma_1+\Gamma_l,\phi(\Gamma_1+\Gamma_k)]
      \quad\mbox{for all}\ k,l\in Y_1.\eqno(4.2.23)$$
      Thus we obtain from (4.2.22) and (4.2.23) that
$$\sum_{r\in Y_1}(\Gamma_k(f_{rl})-\Gamma_l(f_{rk}))D_r
+f_{lk}D_l-f_{kl}D_k=0.\eqno(4.2.24)$$
Comparing coefficients we have
$$\Gamma_k(f_{rl})=\Gamma_l(f_{rk})\quad \mbox{whenever }\ r\neq k,l;\eqno(4.2.25)$$
$$\Gamma_l(f_{kk})=\Gamma_k(f_{kl})-f_{kl}\quad \mbox{whenever}\ k\neq l.\eqno(4.2.26)$$
For $r,k\in Y_1,$ one may assume that
$$f_{rk}=\sum_{|u|=t+1}c_{urk}x^u \quad \mbox{where}\
  c_{urk}\in \mathbb{F}.\eqno (4.2.27)$$
We then obtain from (4.2.25) and (4.2.27) that
$$\delta _{l\in u}c_{urk} = \delta _{k\in u}c_{url} \quad \mbox{whenever}\
r\neq k,l.\eqno(4.2.28)$$ (4.2.28) implies that for $r\in
Y_1\setminus \{k,l\},$
$$ c_{urk} \neq 0\  \mbox{ and }\ l\in u \Longleftrightarrow
c_{url} \neq 0\  \mbox{and }\ k\in u.\eqno(4.2.29)
$$
Let $r, k\in Y_1$ with $ r\neq k.$ If $c_{urk} \neq 0,$ then
(4.2.29) ensures that   $k\in  u.$ Thus
$$\Gamma_k(f_{rk})=f_{rk}\quad \mbox{whenever}\ r\neq k.\eqno(4.2.30)$$
For any fixed $r\in Y_1,$ by Corollary 2.1.5, there is
$\overline{f}_r\in \Lambda(n) $ such that
$$\Gamma_k(\overline{f}_{r})=f_{rk}\quad\mbox{for all}\ k\in Y_1\setminus r. \eqno(4.2.31)$$
Assert that
$$\Gamma_k(f_{kk})=0 \quad\mbox{for all}\ k\in Y_1.\eqno(4.2.32)$$
We treat two cases separately.
\newline

\noindent\textit{Case }(i): $t\geq 2.$ Note that
$\Gamma^2_k=\Gamma_k $ for $k\in Y_{1},$  by Lemma 2.1.2(ii). We
obtain from (4.2.26) that
$$\Gamma_l\Gamma_k(f_{kk})=\Gamma_k\Gamma_l(f_{kk})
=\Gamma^2_k(f_{kl})-\Gamma_k(f_{kl})=0  \quad\mbox{for all}\
 l\in Y_1\setminus  k.\eqno(4.2.33)$$
Since $\mathrm{zd}(f_{kk})=t+1\geq 3 $ and
$\Gamma_k(x^u)=\delta_{k\in u}x^u,$ one may easily deduce (4.2.32)
from (4.2.33). \newline

 \noindent \textit{Case }(ii): $t=0.$
Applying $\phi$ to the equation that
$[x_lD_k,\Gamma_1+\Gamma_k]=x_lD_k $ for $k,l\in Y_1$ with $k\neq
l, $ we have
$$\phi(x_lD_k)+[\Gamma_1+\Gamma_k,\phi(x_lD_k)]
=[x_lD_k,\phi(\Gamma_1+\Gamma_k)].\eqno(4.2.34)$$ On the other
hand, noticing that $\mathrm{zd}(f_{rk})=1,$ it is easily seen from
(4.2.30) that
$$f_{rk}=c_{rk}x_{k} \quad\mbox{whenever }\ r\neq k,\eqno(4.2.35)$$
where $ c_{rk}\in \mathbb{F}.$  Then we obtain from (4.2.22) and
(4.2.35) that
$$\phi(\Gamma_1+\Gamma_k)=f_{kk}D_k
+\sum_{r\in Y_1\setminus k}c_{rk}x_kD_r.\eqno(4.2.36)$$ Combining
(4.2.34) and (4.2.36) and noticing that $D_k(f_{kk})\in
\mathbb{F},$  we have
$$\phi(x_lD_k)+[\Gamma_{1}+\Gamma_k,\phi(x_lD_k)]=(D_k(f_{kk})x_l-c_{lk}x_k)D_k+
\sum_{r\in Y_1 \setminus k}c_{rk}x_lD_r.\eqno(4.2.37)$$ Since
$\mathrm{zd}(\phi) =0,$ one may assume that for $k,l\in Y_1$
 with $k\neq l,$
 $$\phi(x_lD_k)=\sum_{s,r\in Y_1}\mu^{(l,k)}_{s,r}x_sD_r
\quad\mbox{where}\ \mu^{(l,k)}_{s,r}\in \mathbb{F}.\eqno(4.2.38)$$
 Clearly,
$$[\Gamma_1+\Gamma_k,x_sD_r]=(\delta_{ks}-\delta_{kr})x_{s}D_{r}
 \quad \mbox{for}\ k,l\in Y_{1}.\eqno(4.2.39)$$
It follows from (4.2.38) and (4.2.39) that
$$
 \phi(x_lD_k)+[\Gamma_1+\Gamma_k,\phi(x_lD_k)]=
 \sum_{s,r\in Y_1}\mu^{(l,k)}_{s,r}x_sD_r
 +\sum_{s,r\in Y_1}(\delta_{ks}-\delta_{kr})\mu^{(l,k)}_{s,r}x_sD_r.\eqno(4.2.40)
$$
 The coefficient of $D_k$ in the right hand side of (4.2.40) is
 $$\sum_{s\in Y_1}\mu^{(l,k)}_{s,k}x_s
 +\sum_{s\in
 Y_1}(\delta_{ks}-\delta_{kk})\mu^{(l,k)}_{s,k}x_s
 =\sum_{s\in Y_1}\delta_{ks}\mu^{(l,k)}_{s,k}x_s=\mu^{(l,k)}_{k,k}x_k.$$
 We then obtain from (4.2.37) that
 $$D_k(f_{kk})x_l-c_{lk}x_k=\mu^{(l,k)}_{k,k}x_k.$$
 It follows that $D_k(f_{kk})=0,$ proving (4.2.32).

 For $r\in Y_1,$
 put
 $$f_r:=-f_{rr}+(\Gamma_1+\Gamma_r)(\overline{f}_r).$$
 Clearly, $f_r\in \Lambda(n),$ since $\overline{f}_r, f_{rr}\in \Lambda(n).$
 Using (4.2.32), we have
 $$(\Gamma_1+\Gamma_r)(f_r)-f_r=f_{rr}.\eqno(4.2.41)$$
 For $k\in Y_1\setminus r,$ using  (4.2.26) and (4.2.31), we obtain that
 \begin{eqnarray*}
 (\Gamma_1+\Gamma_k)(f_r)&=&-(\Gamma_1+\Gamma_k)(f_{rr})
 +(\Gamma_1+\Gamma_k)(\Gamma_1+\Gamma_r)(\overline{f}_r)\\
 &=&-\Gamma_k(f_{rr})+\Gamma_r\Gamma_k(\overline{f}_r)\\
 &=&-(\Gamma_r(f_{rk})-f_{rk})+\Gamma_r(f_{rk})\\
 &=&f_{rk}.
 \end{eqnarray*}
 Putting $D'=-\sum_{r\in Y_1}f_rD_r, $
  we then obtain by using (4.2.41) that
\begin{eqnarray*}
 [D', \Gamma_1+\Gamma_k]&=&-\sum_{r\in Y_1}[f_rD_r,\Gamma_1+\Gamma_k]\\
 &=&\sum_{r\in Y_1}(\Gamma_1+\Gamma_k)(f_r)D_r-f_kD_k\\
 &=&\sum_{r\in Y_1\setminus k}(\Gamma_1+\Gamma_k)(f_r)D_r+
 ((\Gamma_1+\Gamma_k)(f_k)-f_k)D_k\\
 &=&\sum_{r\in Y_1\setminus k}f_{rk}D_r+f_{kk}D_k\\
 &=&\phi(\Gamma_1 +\Gamma_k).
 \end{eqnarray*}
Let $D:=D'_t.$ Then $D\in\ \mathcal{G}_{t}.$ Since
$\mathrm{zd}(\phi)=t,$ we have
$$[D,\Gamma_1+\Gamma_k]=\phi(\Gamma_1+\Gamma_k)\quad \mbox{for all}\ k\in Y_1.$$
Now it is easy to see that  (4.2.21) holds. \qed\newline

For our purpose, we need still the following three reduction
lemmas.

\begin{lemma}
  Suppose that $\phi \in
\mathrm{Der}_{t}(\mathcal{S},\mathcal{W}),
   $ and $\phi(\mathcal{S}_{-1})=0 $ where $t> 0$ is even.
   If $\phi(\Gamma_1+\Gamma_k)=0 $
   for all $k\in Y_1,$
   then $\phi(\mathcal{S}_{0})=0.$
   \end{lemma}

   \noindent  \textit{Proof.}    (i) First, we show that
   $\phi(\Gamma_1-\Gamma_i)=0$ for all $i\in Y_0\setminus 1.$
   Since $\mathrm{zd}(\phi)$
   is even, by Lemma 2.1.1  we may assume that
   $$\phi(\Gamma_1-\Gamma_i)=\sum_{r\in Y_{1},u\in \mathbb{B}}c^{(i)}_{u,r}x^uD_r
   \quad \mbox{where}\ c^{(i)}_{u,r}\in \mathbb{F}.\eqno(4.2.42)$$
   For any fixed coefficient $c^{(i)}_{v,s},$ the assumption that
    $\mathrm{zd}(\phi)=t>0 $ ensures that there is $k\in v\setminus s.$
   Applying $\phi$
   to the equation that $[\Gamma_1+\Gamma_k,\Gamma_1-\Gamma_i]=0 $
   and using (4.2.42), one gets
   $$\sum_{r\in Y_{1},u\in \mathbb{B}}\delta_{k\in u}c^{(i)}_{u,r}x^uD_r
   -\sum_uc^{(i)}_{u,k}x^uD_k=0.$$
   The choice of  $k$
   and the equation above imply that $c^{(i)}_{v,s}=0.$
  Thus $\phi(\Gamma_1-\Gamma_i)=0$
   for all $i\in Y_0\setminus 1.$

   (ii) Second, we show that $\phi(x_iD_j)=0 $
   for all $i,j\in Y_0$ with $i\neq j.$
   By our general assumption that $m\geq 3,$
   one may find $r\in Y_0\setminus \{i,j\}.$
   Then $[x_iD_j,\Gamma_r+\Gamma_k]=0$
   for all $k\in Y_1.$
   Notice  the assumption that $\phi(\Gamma_1+\Gamma_k)=0 $
   for all $k\in Y_1.$ By (i), it is easily seen that $\phi(\Gamma_r+\Gamma_k)=0. $
   Then
    $[\Gamma_r+\Gamma_k,\phi(x_iD_j)]=0.$
   Arguing as in  (i), one may see that the assertion holds.

   (iii) Finally, we assert that $\phi(x_kD_l)=0$ for $k,l\in Y_1$ with $k\neq l.$
   Just as in (i), one may assume that
   $$\phi(x_kD_l)=\sum_{r\in Y_1, u\in \mathbb{B} }c_{u,r}x^uD_r\quad \mbox{where}\
     c_{u,r}\in \mathbb{F}.\eqno(4.2.43)$$

   \noindent \textit{Case} 1: $\mathrm{zd}(\phi)=2.$
   Note that $n\Gamma_1+\Gamma'\in \mathcal{S }$ and $\phi(n\Gamma_1+\Gamma')=0 $ by
   the assumption.
   Noticing (4.2.43), we have
   $$2\phi(x_kD_l)=[\Gamma',\phi(x_kD_l)]= [n\Gamma_1+\Gamma',\phi(x_kD_l)]=0.$$
   Therefore, $\phi(x_kD_l)=0.$\newline

   \noindent\textit{Case} 2: $\mathrm{zd}(\phi)\geq 4.$
   For any given coefficient $c_{v,s} $ in (4.2.43), one may take
   $q\in v\setminus \{k,l,s\},$
   since $|v|=\mathrm{zd}(\phi)+1\geq 5.$
   Then
   $[\Gamma_1+\Gamma_q,x_kD_l]=0 $
   and therefore,
   $$[\Gamma_1+\Gamma_q,\sum _{r\in Y_1,u\in \mathbb{B}}c_{u,r}x^uD_r]=0.$$
  Furthermore,
   $$\sum_{r\in Y_1, u\in \mathbb{B}}\delta_{q\in u}c_{u,r}x^uD_r-\sum_uc_{u,q}x^uD_q=0. $$
     Consequently, $c_{v,s}=0 $
   and then
   $\phi(x_kD_l)=0 $ for all $k,l\in Y_1$ with $k\neq l.$
   Summarizing, we have
   $\phi(\mathcal{S}_{0})=0.$
   \qed

   \begin{lemma}
     Suppose that $\phi \in
   \mathrm{Der} (\mathcal{S},\mathcal{W})
   $ is a homogeneous derivation of positive odd $\mathbb{Z}$-degree
   and $\phi(\mathcal{S}_{-1}+\mathcal{T_{\mathcal{S}}})=0.$
   Then $\phi(\mathcal{S}_{0})=0.$
   \end{lemma}

\noindent  \textit{Proof.}  \ \  We first assert that
$$\phi(x_iD_j)=0\quad\mbox{ for all}\  i,j\in Y_0\
\mathrm{with}\ i\neq j.$$ Since $\mathrm{zd}(\phi)$ is odd, by Lemma
2.1.1 one may assume that
$$\phi(x_iD_j)=\sum_{r\in Y_0,u}c_{u,r}x^uD_r
 \quad \mbox{where}\ c_{u,r}\in \mathbb{F}.$$
Since $m\geq 3,$ one may take $s\in Y_0\setminus\{i,j\}.$ Clearly,
$$[\Gamma_s+\Gamma_k,x_iD_j]=0\quad \mbox{for all}\ k\in Y_1.$$
Applying $\phi$ to the equation above, we have
$$\sum_{r\in Y_0,u\in \mathbb{B}}\delta_{k\in u}c_{u,r}x^uD_r-\sum_uc_{u,s}x^uD_s=0.$$
Consequently,
$$\delta_{k\in u}c_{u,r}=0\quad \mbox{for}\ r\in Y_0\setminus s;\ k\in Y_1.
\eqno(4.2.44)$$
Take $k\in u.$ We  obtain from (4.2.44) that
$c_{u,r}=0 \ \mbox{for}\ r\in Y_0\setminus s.$ Thus
$$\phi(x_iD_j)=\sum_{u\in \mathbb{B}}c_{u,s}x^uD_s.\eqno(4.2.45)$$
Applying $\phi$ to the following equation
$$[x_jD_j+\Gamma_k,x_iD_j]=-x_iD_j \quad \mbox{for all }\ k\in
Y_1,$$ we obtain by using (4.2.45) that
$$[x_jD_j+\Gamma_k,\sum_u c_{u,s}x^uD_s]=-\sum_u c_{u,s}x^uD_s. $$
Then
$$\sum_u\delta_{k\in u} c_{u,s}x^uD_s=-\sum_u c_{u,s}x^uD_s.$$
A comparison of coefficients shows that
$$\delta_{k\in u} c_{u,s}=-c_{u,s}\quad \mbox{for all}\  k\in Y_1.\eqno(4.2.46)$$
Take  $k\in u.$ Then (4.2.46) yields $c_{u,s}=0$ for all $u\in \mathbb{B}.$ Then
we obtain from (4.2.45) that
$$\phi(x_iD_j)=0\quad \mbox{for}\ i,j\in Y_0 \ \mbox{with}\ i\neq j.$$
The assertion holds.

 It remains only to show that
$$\phi(x_kD_l)=0\quad \mbox{for}\ k,l\in Y_1\ \mbox{with}\ k\neq l.$$
As in the above, we may assume that
$$\phi(x_kD_l)= \sum_{r\in Y_0,u\in \mathbb{B}}c_{u,r}x^uD_r
   \quad \mbox{where}\  c_{u,r}\in \mathbb{F}.$$
Clearly,
$$[x_kD_l,\Gamma_i-\Gamma_j]=0 \quad \mbox{for }\  i,j\in Y_0\
\mbox{with }\  i\neq j.$$
Applying $\phi$  to this equation, we obtain that
$$\sum_{u}c_{u,i}x^uD_i-\sum_{u}c_{u,j}x^uD_j=0 $$
and therefore, $c_{u,i}=c_{u,j}=0.$ Hence $c_{u,r}=0 $ for all $r$
and all $u.$ The proof is complete. \qed\newline

\begin{lemma}
Suppose that $\phi \in \mathrm{Der}_{0}(\mathcal{S},\mathcal{W})
   $ and $\phi(\mathcal{W}_{-1})=0.$ If
   $$\phi(\Gamma_1+\Gamma_k)=0\quad \mbox{for all }\ k\in Y_1, $$
     then there is $D\in \mathcal{G}_0 $
     such that
      $(\phi-\mathrm{ad}D)(\mathcal{S}_{0})=0.$
 \end{lemma}

\noindent  \textit{ Proof.}   By Lemma 2.1.1,
 $\phi(x_iD_j)\in E(\mathcal{G})$ for all $i,j\in Y_0 $ with $i\neq j.
  $
Then
$$\phi(\Gamma_1-\Gamma_i)=[\phi(x_1D_i),x_iD_1]+[x_1D_i,\phi(x_iD_1)]=0
 \quad \mbox{for }\ i\in Y_0\setminus 1. $$
We next show that $\phi(x_iD_j)=0$ for $i,j\in Y_0 $ with $i\neq
j.$ Find $r\in Y_0\setminus \{i,j\}, $ since $m\geq 3.$ Then
$$\phi(x_iD_j)=\phi([\Gamma_i-\Gamma_r,x_iD_j ])=0, $$
since $\phi(\Gamma_i-\Gamma_r),\phi(x_iD_j)\in E(G).$ Finally,
consider $\phi(x_kD_l)$ for $k,l\in Y_1$ with $k\neq l.$
We may
assume that
$$\phi(x_kD_l)=\sum_{r\in Y_1}f_{r}D_r
\quad\mbox{where} \ f_r\in \Lambda (n)_{1}. $$
For $q\in
Y_1\setminus \{k,l\}, $ we have $[\Gamma_1+\Gamma_q, x_kD_l ]=0 $
and therefore,
$$\sum_{r}(\Gamma_1+\Gamma_q)(f_r)D_r-f_qD_q=0.$$
Consequently, $(\Gamma_1+\Gamma_q)(f_q)=f_q.$ This implies that
$f_q=c^{(q)}_{k,l}x_q$ for $  q\in Y_1\setminus \{k,l\},$ since
$f_q\in \Lambda (n)_{1}.$ Therefore,
$$\phi(x_kD_l)=\sum_{r\in Y_1\setminus\{ k,l\}}c_{r}x_rD_r+f_kD_k+f_lD_l.$$
Note that
 $[\Gamma_1+\Gamma_k, \phi(x_kD_l) ]=\phi(x_kD_l). $
Then
$$(\Gamma_1+\Gamma_k)(f_k)D_k-f_kD_k+\Gamma_k (f_l)D_l=
\sum_{r\in Y_1\setminus\{ k,l\}}c_{r}x_rD_r+f_kD_k+f_lD_l.$$
Comparing
coefficients we have
 $$c_r=0 \quad\mbox{for} \ r\in Y_1\setminus \{k,l\};$$
$$(\Gamma_1+\Gamma_k)(f_k)=2f_k; $$
$$\Gamma_k(f_l)=f_l. $$
 The last two equations  imply that $f_k=0 $  and
 there is  $\lambda _{kl}\in \mathbb{F} $
 such that
 $f_l=\lambda _{kl}x_k.$
Thus
$$\phi(x_kD_l) =\lambda _{kl}x_kD_l.$$
Just as in the proof of Proposition 3.2.4, one may show from the
equation above that there is $ D\in
\mathcal{G}_0 $ such that
$$(\phi-\mathrm{ad}D)(\mathcal{S}_{0})=0.$$
\qed\newline

Now we are able to characterize the homogeneous derivation space
of nonnegative $\mathbb {Z}$-degree.

\begin{proposition} 
$\mathrm{Der}
_{t}(\mathcal{S},\mathcal{W})=\mathrm{ad}\mathcal{W}_{t}$ for
$t\geq 0.$
\end{proposition}

\noindent\textit{Proof.} Clearly,
$\mathrm{ad}\mathcal{W}_{t}\subset\mathrm{Der}
_{t}(\mathcal{S},\mathcal{W}).$ To prove the converse inclusion,
let $\phi\in\mathrm{Der} _{t}(\mathcal{S},\mathcal{W}).$ In the
light of Proposition 2.1.6, we may assume that $\phi(
\mathcal{S}_{-1})=0.$ We treat two cases separately.

(i) Suppose that $\mathrm{zd} (\phi) $ is odd. By Lemma 4.2.4 and Corollary 4.2.3,
there is $ E\in \mathcal{G}$ such that $\phi -\mathrm{ad} E$ vanishes on
$\mathcal T_{\mathcal S}. $ Then Lemma 4.2.7 shows  that
$\phi -\mathrm{ad} E$ vanishes on
$\mathcal  \mathcal S_{0}.$ Clearly,
$(\phi -\mathrm{ad} E)(\mathcal{S}_{-1})=0.$ Now Corollary 4.1.4 ensures that
$\phi\in \mathrm {ad} W$ and therefore,  $\phi\in \mathrm {ad} W_{t}.$

(ii)   Suppose that $\mathrm{zd} (\phi) $ is even. By Lemma 4.2.5, there is
$D\in \mathcal{G}_{t}$ such that
 $(\phi-\mathrm{ad} D)(\Gamma_{1}+\Gamma _{k})=0 $  for all $k\in Y_{1}.$
 If $t>0,$ then  Lemma 4.2.6 implies that
 $(\phi-\mathrm{ad} D)(\mathcal{S}_{0})=0.$  If $t=0,$
 then   Lemma 4.2.8 shows that there is $E\in \mathcal{G}_{0}$ such that
 $(\phi-\mathrm{ad} D-\mathrm{ad} E)(\mathcal{S}_{0})=0.$ Now, Corollary 4.1.4
 ensures that $\phi\in \mathrm{Der} _{t}(\mathcal{S},\mathcal{W})$ is inner. The proof
 is complete.
 \qed\newline

 As an application of Proposition 4.2.9, we have:

 \begin{proposition}
$\mathrm{Der} _{t}(\mathcal{S} )
=\mathrm{ad}(\overline{\mathcal{S}}+\mathcal{T})_{t}$ for $t\geq
0.$
 \end{proposition}

\noindent\textit{Proof.} Just as in the case of Lie algebras, one may prove that
 $\mathcal S$ is an ideal of $\overline{\mathcal S}$
and $\mathcal{T}\subset \mathrm{Nor}_{\mathcal W}(\mathcal S)$
(cf. \cite{Z1,ZZ}). It follows that $
\mathrm{ad}(\overline{\mathcal{S}}+
\mathcal{T})_{t}\subset\mathrm{Der} _{t}(\mathcal{S} ).$  To prove
the converse
 inclusion, let $\phi\in\mathrm{Der} _{t}(\mathcal{S} ).$ One may identity $\phi$ with
 a derivation of $\mathrm{Der} _{t}(\mathcal{S},\mathcal{W}).$
 Then by Proposition 4.2.9, there is $D\in \mathcal{W}_t$ such that
 $\phi=\mathrm{ad}D$ as elements of  $\mathrm{Der} _{t}(\mathcal{S},\mathcal{W})$
 and therefore,
 $\phi=\mathrm{ad}D$ as elements of $\mathrm{Der} _{t}(\mathcal{S} ).$ Clearly,
 $D \in \mathrm{Nor}_\mathcal{W}(\mathcal{S})_{t}.$ Therefore, it is enough to show that
 $\mathrm{Nor}_\mathcal{W}(\mathcal{S})_{t}
 \subset (\overline{\mathcal{S}}+\mathcal{T})_{t}. $
  Let $E$ be an arbitrary element of
  $\mathrm{Nor}_\mathcal{W}(\mathcal{S})_{t}.$ Then
  $\mathrm{div} ([D_{i}, E] )=0 $ for all $i\in Y_{0}.$
  Since $\mathrm{div}$ is a derivation of $\mathcal{W}$ to $\frak{A} $ and
  $\mathrm{div}(D_{i})=0,$
  it follows that $D_i(\mathrm{div}(E))=0 $ for all $i\in Y_0.$ This implies that
  $\mathrm{div}(E)\in \Lambda(n)_{t}.$

  If $t=0 $ then $\mathrm{div}(E)\in \mathbb{F}.$ It follows that
   $\mathrm{div}(E- \mathrm{div}(E) \Gamma_1)=0.$ Hence
   $E\in(\overline{\mathcal{S}}+\mathcal{T})_{0}=
   \overline{\mathcal{S}}_{0}+\mathcal{T}.$

   Suppose that $t>0.$
   We contend that $\mathrm{Nor}_\mathcal W (\mathcal S)_t\subset
   \overline{\mathcal{S}}_{t}.$
   Assume that on the contrary that
   $\mathrm{Nor}_\mathcal W (\mathcal S)_t\not\subset
   \overline{\mathcal{S}}_{t}.$
   Then there is
   $E\in\mathrm{Nor}_\mathcal W (\mathcal S)_t\setminus
   \overline{\mathcal{S}}_{t}.$
   Since $\mathrm{div}(E)\in \Lambda(n),$
   one may write $E$ to be the following form
   $$ E=\sum_{i\in X_0, u\in \Omega_0}c_{ui}x^{u}x_{i}D_{i}+
        \sum_{k\in v, v\in \Omega_1}c_{vk}x^{v}D_{k}+G \eqno(4.2.47)
   $$
   such that
   $$\mathrm{div}(E)=\sum_{i\in X_0, u\in \Omega_0}c_{ui}x^{u} -
        \sum_{k\in v, v\in \Omega_1}c_{vk}D_{k}(x^{v})\eqno(4.2.48)
   $$
   and that no any cancellation occurs in the right hand side of  (4.2.48),
   where $X_0\subset Y_{0}, X_1\subset Y_1;$ $\Omega_{0}\subset \cup_{1\leq r\leq
   \frac{n}{2}}\mathbb {B}_{2r}, \Omega_{1}\subset \cup_{1\leq r\leq
   \frac{n-1}{2}}\mathbb {B}_{2r+1};$ $c_{ui}, c_{vk}\in \mathbb {F};$
   $G\in \overline{\mathcal{S}}.$
   Evidently, the assumption that
   $\mathrm{Nor}_\mathcal W (\mathcal S)_t\not\subset
   \overline{\mathcal{S}}_{t} $
   secures that
   $X_{0}\neq \emptyset$
   and
   $\Omega_0\neq \emptyset,$ or,
   $X_{1}\neq \emptyset$
   and
   $\Omega_1\neq \emptyset. $
   Assume that $X_{0}\neq \emptyset$
   and
   $\Omega_0\neq \emptyset.$ Given $w\in\Omega_0 $ and $j\in X_{0},$ choose $r\in  w.$
   It follows from (4.2.47) that
   $$
   \sum_{i\in X_0, u\in \Omega_0}
   \delta_{r\in u}c_{ui}x^{u}x_{i}D_{i}+
   \sum_{k\in v, v\in \Omega_1}(\delta_{r\in v}-\delta_{rk})c_{vk}x^{v}D_{k}
   =[\Gamma_r+\Gamma_1, E]-[\Gamma_r+\Gamma_1, G]\in \mathcal {S},
   $$
   since $\Gamma_r+\Gamma_1\in \mathcal{S}.$ Applying the divergence,
     we obtain from
   the equation above that
   $$
   \sum_{i\in X_0, u\in \Omega_0}\delta_{r\in u}c_{ui}x^{u} -
   \sum_{k\in v, v\in \Omega_1}(\delta_{r\in v}-\delta_{rk})
   c_{vk}D_{k}(x^{v})=0.\eqno(4.2.49)
   $$
   Note that $\delta_{r\in v}-\delta_{rk}=0$ or $1$. One should bear in mind,
    as remarked above, that no cancellation occurs in the right hand side of
    (4.2.48). Thus the nonzero summand  $c_{wj}x^{w}$ in the first sum in the
    left hand side of (4.2.49)
    cannot be cancelled, contradicting to that the right hand side is zero.

    It
    remains to discuss the case that  $X_{1}\neq \emptyset$
   and
   $\Omega_1\neq \emptyset.$ Given $w\in\Omega_1 $ and $l\in X_{1},$ one can
   choose $r\in  w\setminus l $ since $|w|\geq 3.$ Arguing as above one may find that
    $c_{wl}D_{l}(x^{w})$ is a nonzero summand in the second sum in the left hand side
    of (4.2.49), which cannot be cancelled by the same token, contradicting that the
    right hand side is zero.

    So far, we have proved that $\mathrm{Nor}_\mathcal W (\mathcal S)_t\subset
   \overline{\mathcal{S}}_{t} $ for $t>0.$
   Summarizing, the proof is complete. \qed

   \begin{remark}
   {\rm Our original idea is to study the derivation
   algebra $\mathrm{Der}(\mathcal{S}) $ rather than the derivation space
   $\mathrm{Der}(\mathcal{S}, \mathcal{W}), $ which contains the
   former in the obvious sense. But, in practice, it is convenient and
   effective first to study the derivation space
   $\mathrm{Der}(\mathcal{S}, \mathcal{W})  $ rather than the algebra
   $\mathrm{Der}(\mathcal{S}).  $ Our work had ever stopped for a time,
   since we observed that the natural $\mathbb{Z}$-gradation is not admissibly graded
   (see Remark 2.1.8). However,  when we considered the
   derivation space $\mathrm{Der}(\mathcal{S}, \mathcal{W})
   $ and determined it at last, almost all problems were solved at the last
   moment.}
   \end{remark}

\subsection{Derivation algebra of $\mathcal{S}$}

In this subsection we first determine the homogeneous derivations of negative
$\mathbb{Z}$-degree of $\mathcal{S}$ to $\mathcal{W}$.
This combining with the results obtained in Section 4.2 will give the structure of the
derivation space $\mathrm{Der}(\mathcal{S},\mathcal{W}).$  Using these results we are
able to
characterize the derivation algebra $\mathrm{Der}(\mathcal{S}).$

To compute the derivations of negative $\mathbb{Z}$-degree, recall
the generator set of $\mathcal{S} $ (see Proposition 2.2.3). We
still adopt the notations
$$\mathcal{R}:=\{D_{il}(x^{(2\varepsilon_i)}x_k) \mid i\in
Y_0,k,l\in Y_1\}\cup \{D_{ij}(x_ix_kx_l\mid  i,j\in Y_0,k,l\in
Y_1\} $$
and
$$\mathcal{Q}=\{ D_{ij}(x^{(a\varepsilon_j)})\mid i,j\in Y_0,i\neq
j,a\in \mathbb{N}\}.$$

The following lemma tells us that a derivation of
$\mathbb{Z}$-degree $-1$ of $\mathcal{S}$ to $\mathcal{W}$ is
completely determined by its action on $\mathcal{S}_{0}.$

\begin{lemma}

Suppose that $\varphi \in
\mathrm{Der}_{-1}(\mathcal{S},\mathcal{W}) $
       and $\varphi(\mathcal{S}_{0})=0.$
       Then $\varphi=0.$
\end{lemma}

\noindent \textit{Proof. }   First of all, we show that
$\varphi(\mathcal{R})=0.$
  We
shall use the following simple fact (by Lemma 2.1.1):
$$\varphi(\mathcal{S}_1)\subseteq \mathcal{G}_0\subseteq E(\mathcal{G}).\eqno(4.3.1)$$

Given $i\in Y_0, k,l\in Y_1,$ take  $j\in Y_0\setminus i.$ Then
$$[\Gamma_i-\Gamma_j,D_{il}(x^{(2\varepsilon_i)}x_k) ]
=D_{il}(x^{(2\varepsilon_i)}x_k).\eqno(4.3.2)$$
From (4.3.1) and (4.3.2), we obtain that
$$\varphi\big(D_{il}(x^{(2\varepsilon_i)}x_k)\big)=0
\quad\mbox{for}\ i\in Y_0,k,l\in Y_1.\eqno(4.3.3)$$

By a same argument, we can also obtain that
$$\varphi\big(D_{il}(x_ix_kx_l)\big)=0 \quad\mbox{for all}\ i\in Y_0,k,l\in Y_1.
\eqno(4.3.4)$$
It follows from (4.3.3) and (4.3.4) that
$\varphi(R)=0.$

We next show that $ \varphi(\mathcal{Q})=0.$
If $a=3,$
 it is easily showed  as above that
  $\varphi\big(D_{i,j}(x^{(a\varepsilon_j)})\big)=0.$
Now suppose that $a\geq 4.$ We proceed by induction on $a$ to show
that
$$\varphi\big(D_{ij}(x^{(a\varepsilon_j)})\big)=0,\quad a\geq 4.\eqno(4.3.5)$$
By inductive hypothesis,
$$\varphi\big(D_{i j}(x^{(a\varepsilon_j)})\big)\in \mathcal{G}_{a-3}.\eqno(4.3.6)$$
If $a$ is odd, then (4.3.6) ensures that
$$\varphi\big(D_{i j}(x^{(a\varepsilon_j)})\big)
= \sum_{r\in Y_1, u\in \mathbb{B}_{a-2}} c_{u,r}x^uD_r\quad\mbox{where}\
c_{u,r}\in\mathbb{ F}. \eqno(4.3.7)$$
Given any coefficient
$c_{v,s}$ in (4.3.7), one may take $k\in v\setminus r,$ since
$|u|=a-2 \geq 2.$ One may also take $q\in Y_0\setminus \{i,j\},$
since $m\geq 3.$ Then
$$[\Gamma_k+\Gamma_q, D_{ij}(x^{(a\varepsilon_j)})]=0.\eqno(4.3.8)$$
Applying $\varphi$ to (4.3.8) and using (4.3.7), one gets
$c_{v,s}=0.$ Thus (4.3.5) holds.

If $a$ is even, by (4.3.6),
$$\varphi\big(D_{ij}(x^{(a\varepsilon_j)})\big)
= \sum_{r\in Y_0, u\in \mathbb{B}_{a-2}} c_{u,r}x^uD_r.
\eqno(4.3.9)$$  Given any coefficient $c_{v,s}$ in (4.3.9), take
$k\in v, q\in Y_0\setminus \{s,j\}.$ Then
$$[\Gamma_k+\Gamma_q, D_{i j}(x^{(a\varepsilon_j)})]
=-\delta_{qj}D_{ij}(x^{(a\varepsilon_j)}).\eqno(4.3.10)$$ From
(4.3.9) and (4.3.10), one may compute that
 $c_{v,s}=-\delta_{qj}c_{v,s}.$
Therefore, $c_{v,s}=0 $ and (4.3.5) holds.
 This proves $\varphi(\mathcal{Q})=0.$
\qed\newline

 Using Lemma 4.3.1 we can  determine the derivations of $\mathbb{Z}$-degree
 $-1;$ in particular, they
 are all inner.

\begin{proposition}
$\mathrm{Der}_{-1}(\mathcal{S},\mathcal{W})=\mathrm{ad}
\mathcal{W}_{-1}.$ In particular, $\mathrm{Der}_{-1}(\mathcal{S}
)=\mathrm{ad} \mathcal{S}_{-1}.$
\end{proposition}

  \noindent \textit{Proof.}   Let $\phi\in \mbox{Der}_{-1}(\mathcal{S},\mathcal{W}). $
  For $i\in Y_0,$  suppose that
  $$\phi(\Gamma_{1'}+\Gamma_i)=\sum_{r\in Y_0}c_{i,r}D_r \quad\mbox{where}\
   c_{i,r}\in \mathbb{F}.\eqno(4.3.11)$$
  Let $j\in Y_0\setminus i.$
  Then $[\Gamma_{1'}+\Gamma_i,\Gamma_{1'}+\Gamma_j]=0$
  and therefore,
  $$[\Gamma_{1'}+\Gamma_i,\phi(\Gamma_{1'}+\Gamma_j)]
  =[\Gamma_{1'}+\Gamma_j,\phi(\Gamma_{1'}+\Gamma_i)].$$
    Then by (4.3.11),
    $$c_{i,j}=0 \quad\mbox{whenever}\  i,j\in Y_0\ \mbox{with}\ i\neq j.$$
    Thus, by (4.3.11), one gets
    $$\phi(\Gamma_{1'}+\Gamma_i)=c_{i,i}D_i
    \quad \mbox{where}\ c_{i,i}\in \mathbb{F}.\eqno(4.3.12)$$
    Obviously,
    $$[\Gamma_{1'}+\Gamma_i,x_iD_j]=x_iD_j
    \quad\mbox{for}\ i,j\in Y_0 \quad\mbox{with}\ i\neq j.$$
    Using (4.3.12), we then have
    $$c_{i,i}D_j+[\Gamma_{1'}+\Gamma_i,\phi(x_iD_j)]=\phi(x_iD_j).$$
    Noticing that $\phi(x_iD_j)\in \mathcal{S}_{-1}=\mathcal{W}_{-1},$
    we obtain from the equation above  that $\phi(x_iD_j)=c_{i,i}D_j.$
    Now let $\psi:=\phi- \sum_{r\in Y_0} c_{r,r}\mathrm{ad} D_r.$ Then
    $$\psi(\Gamma_{1'}+\Gamma_i)=\psi(x_iD_j)= 0
    \quad\mbox{for}\ i,j\in Y_{0} ,\ i\neq j.\eqno(4.3.13)$$
    We want to prove that
    $$\psi(x_kD_l)=0\quad\mbox{for}\ k,l\in Y_1 \ \mbox{with}\
    k\neq l.\eqno(4.3.14) $$
    To do that, we choose $q\in Y_1\setminus \{k,l\}.$
    Then
    $$[\Gamma_i+\Gamma_q,x_kD_l]=0\quad\mbox{for all}\ i\in Y_0. $$
    It follows that
    $$[\psi(\Gamma_i+\Gamma_q) ,x_kD_l]
    +[\Gamma_i+\Gamma_q,\psi(x_kD_l)]=0.\eqno(4.3.15)$$
   Since $\psi(\Gamma_i+\Gamma_q)\in \mathcal{W}_{-1},$
   we have $[\psi(\Gamma_i+\Gamma_q) ,x_kD_l]=0.$
   Then (4.3.15) implies that $[\Gamma_i,\psi(x_kD_l)]=0 $
   for all $i\in Y_0.$
   Hence $\psi(x_kD_l)=0,$
   since $\psi(x_kD_l)\in \mathcal{W}_{-1}.$

   Similarly, we may check that
   $$\psi(\Gamma_{1'}-\Gamma_k)=0 \quad\mbox{for all}\ k\in Y_1\setminus 1'.\eqno(4.3.16)$$
  By (4.3.13), (4.3.14) and (4.3.16),  $\psi(\mathcal{S}_{0})=0.$
   By Lemma 4.3.1, we obtain that  $\psi=0$
   and $\phi \in \mathrm{ad}\mathcal{W}_{-1}.$
   \qed\newline

   To compute the derivation of $\mathbb{Z}$-degree  less than $-1 $ of $\mathcal{S}$
    to $\mathcal{W},$ we establish  the following lemma.

  \begin{lemma}
    Let $\phi \in \mathrm{Der}_{-t}(\mathcal{S},\mathcal{W}) ,\ t>1.$
      Suppose that
       $\phi( D_{ij}(x^{((t+1)\varepsilon_i)}))=0 $ for all $i,j\in Y_{0}.$
      Then $\phi=0.$
   \end{lemma}

  \noindent  \textit{Proof. }  First claim that $\phi(\mathcal{Q})=0. $
  To that aim, we proceed by induction on $q$ to show that
  $$\phi( D_{ij}(x^{(q\varepsilon_i)}))=0
  \quad\mbox{for all}\ i,j\in Y_{0}.\eqno(4.3.17)$$
  If $q\leq t+1,$ then (4.3.17) holds.  Suppose that $q> t+1 $ in the following.
  By inductive hypothesis and Lemma 2.1.1,
  $\phi( D_{ij}(x^{(q\varepsilon_i)}))\in \mathcal{G}_{q-t-2}.$
  Then one may assume that
  $$\phi( D_{ij}(x^{(q\varepsilon_i)}))=\sum_{r\in Y,|u|=q-t-1}
  c_{u,r}x^uD_r \quad\mbox{where}\ c_{u,r} \in \mathbb{F}.\eqno(4.3.18)$$
  We treat two cases separately.\newline

 \noindent \textit{Case} (i): $q-t\geq 3.$
  For any fixed coefficient $c_{u_0,r_0} $ in (4.3.18),
  choose
   $k\in u_0\setminus r_0 $ since $|u_{0}|\geq 2.$
  Choose also $s\in Y_0\setminus \{r_0,i\}.$
  Then
  $$[\Gamma_s+\Gamma_k,D_{ij}(x^{(q\varepsilon_i)}) ]
  =-\delta_{sj}(D_{ij}(x^{(q \varepsilon_i)})).$$
  Applying $\phi$ to the equation above and then combining that with (4.3.18),
  one may obtain by a comparison of the coefficients of  $x^{u_{0}}D_{r_{0}}$
  that
  $$c_{u_{0},r_{0}}=-\delta_{sj}c_{u_{0},r_{0}}.$$
  Consequently, $c_{u_{0},r_{0}}=0 $
  and $\phi(D_{ij}(x^{(q\varepsilon_i)}))=0.$\newline

 \noindent \textit{Case} (ii):  $q-t<3.$
 Note that $q> t+1 $ and then
  $q-t-1=1.$
  Then rewrite (4.3.18) as
  $$\phi(D_{ij}(x^{(q\varepsilon_i)}))
  =\sum_{s\in Y_1,r\in Y_1}c_{s,r}x_sD_r\quad \mbox{where}\
   c_{s,r}\in \mathbb{F}. $$
  Arguing as in Case (i),  one may easily obtain that $c_{s,r}=0 $
  whenever $s\neq r.$
  Then
  $$\phi(D_{ij}(x^{(q\varepsilon_i)}))=\sum_{r\in Y_1}c_{r,r}x_rD_r. $$
  Applying $\phi$ to the identity $[x_sD_l,D_{ij}(x^{(q\varepsilon_i)}) ]=0$ for
  $ s,l\in Y_1,$
  we have  $c_{s,s}=c_{l,l}.$
  Let $\lambda:=c_{r,r},\ r\in Y_1.$

  Then
   $\phi(D_{ij}(x^{(q\varepsilon_i)}))=\lambda \Gamma'.$
  Clearly, $x_lx_sD_r\in\mathcal{ S}$ for
   $l,s\in Y_1,r\in Y_0\setminus i.$ Note that
   $[x_lx_sD_r,D_{ij}(x^{(q\varepsilon_i)})]=0.$
 Applying $\phi$  to this equation, one has
  $$[x_lx_sD_r,\lambda \Gamma']+
  [\phi(x_lx_sD_r), D_{ij}(x^{(q\varepsilon_i)})]=0.\eqno(4.3.19)$$
  Note that  $[x_lx_sD_r,\lambda \Gamma']
  =-2\lambda x_lx_sD_r$ and $ \phi(x_lx_sD_r)\in \mathcal{S}_{-1}.$
  Then by (4.3.19), $\lambda =0.$
  Thus, (4.3.17) holds for all $q $
  and therefore, $\phi(\mathcal{Q})=0.$

  We next prove that $\phi(\mathcal{R})=0.$
  Since  $\mathcal{R}\subseteq \mathcal{S}_1,$ $\mathrm{zd}(\phi)\leq -2,$
  it suffices to consider the case that $\mathrm{zd}(\phi)= -2.$
  Note that  $\phi(\mathcal{S}_1)\subset \mathcal{S}_{-1}.$

For $i\in Y_0, k,l\in Y_1,$ choose $q \in Y_1\setminus \{k,l\}.$
Then $m\Gamma_q+\Gamma''\in \mathcal{S}_0$ and
$$[m\Gamma_q+\Gamma'',D_{il}(x^{(2\varepsilon_i)}x_k)]=D_{il}
(x^{(2\varepsilon_i)}x_k). $$
Applying $\phi,$ one gets
$$\phi(D_{il}(x^{(2\varepsilon_i)}x_k))
=[m\Gamma_q+\Gamma'',\phi(D_{il}(x^{(2\varepsilon_i)}x_k))]
=-\phi(D_{il}(x^{(2\varepsilon_i)}x_k)).$$
Since
$\phi(D_{il}(x^{(2\varepsilon_i)}x_k))\in \mathcal{W}_{-1},$
$$\phi(D_{il}(x^{(2\varepsilon_i)}x_k))=0
\quad \mbox{for all} \ i\in Y_0, k,l \in Y_1. $$ Note that
$n\Gamma_q+\Gamma'\in S_0$ for $q\in Y_0 $ and that
$$[n\Gamma_q+\Gamma',D_{ji}(x_jx_kx_l)]
=(2-n\delta_{qi})D_{ji}(x_jx_kD_l)\quad\mbox{for}\ i,j\in
Y_0,k,l\in Y_1. \eqno(4.3.20)$$
 Assume that  $$\phi(D_{ji}(x_jx_kx_l))=\sum
_{r\in Y_0}a_rD_r \quad \mbox{where} \ a_r\in \mathbb{F}. $$
Then
we obtain from (4.3.20) that
$$-a_qD_q=(2-n\delta_{qi})\sum_{r\in Y_0}a_rD_r.
$$
Since $p\neq 3,$ it follows that $a_r=0 $ for all $r\in Y_0.$

This proves that
$$\phi(D_{ji}(x_jx_kx_l))=0\quad \mbox{for all} \ i,j\in Y_0,k,l\in Y_1. $$
By Proposition 2.2.3, $\phi=0.$ \qed\newline

We  are  in  the position to determine the homogeneous derivations
of $\mathbb{Z}$-degree $<-1$ of $\mathcal{S}$ to $\mathcal{W}.$ We
first give the following fact.

\begin{proposition}

Suppose that $t>1$ is not any $p$-power. Then
 $\mathrm{Der}_{-t}(\mathcal{S},\mathcal{W})=0.$
 In particular,
 $\mathrm{Der}_{-t}(\mathcal{S} )=0.$
\end{proposition}

\noindent \textit{Proof.} Let $\phi\in
\mbox{Der}_{-t}(\mathcal{S},\mathcal{W}).$ In view of Lemma 4.3.3,
it is sufficient to show that
$$\phi(D_{ij}(x^{((t+1)\varepsilon_i)}))
=0\quad \mbox{for all} \ i,j\in Y_0.\eqno(4.3.21)$$ We treat two
 cases separately.\newline

\noindent \textit{Case} (i):  $t\not\equiv 0 \pmod{p}.$ Recall
$\Gamma''=\sum _{r\in Y_0}\Gamma_r.$ Clearly,
$m\Gamma_k+\Gamma''\in \mathcal{S}_0 $ for any $k\in Y_1 .$ By
Lemma 2.1.2(viii),
$$[m\Gamma_k+\Gamma'',D_{ij}(x^{(t+1)\varepsilon_i})]
=(t-1)D_{ij}(x^{(t+1)\varepsilon_i}). $$   Applying
 $\phi$ to the equation above, we have
$$[m\Gamma_k+\Gamma'',\phi(D_{ij}(x^{(t+1)\varepsilon_i}))]
=(t-1)\phi(D_{ij}(x^{(t+1)\varepsilon_i})).\eqno(4.3.22)$$
On the
other hand, since $\phi(D_{ij}(x^{(t+1)\varepsilon_i}))\in
\mathcal{S}_{-1}, $
$$[m\Gamma_k+\Gamma'',\phi(D_{ij}(x^{(t+1)\varepsilon_i}))]
=-\phi(D_{ij}(x^{(t+1)\varepsilon_i})).\eqno(4.3.23)$$ A
comparison of (4.3.22) and (4.3.23) shows that (4.3.11) holds,
since $t\not\equiv 0\pmod{p}.$\newline

\noindent\textit{Case }(ii): $t\equiv 0 \pmod{p}.$ Write $t$ to be the
$p $-adic expression
$$t=\sum^{r}_{s=1}a_sp^s \quad\mbox{where}\ 0\leq a_s<p\ \mbox{and}\ a_r\not=0.$$
Note that
$$\mathrm{zd}(D_{ij}(x^{((t-p^r+2)\varepsilon_j)}))=t-p^r<t-2;$$
$$\mathrm{zd}(D_{ij}(x^{(p^r\varepsilon_j)}x_i))=p^r-1<t-2, $$  since $t$ is not any
$p$-power.
Then
$$\phi(D_{ij}(x^{((t-p^r+2)\varepsilon_j)}))
=\phi(D_{ij}(x^{(p^r\varepsilon_j)}x_i))=0.\eqno(4.3.24)$$ Direct
computation shows that (by using the fact that $\binom{t}{p^{r}-1
}\equiv 0 \pmod{p}$ )
$$[D_{ij}(x^{((t-p^r+2)\varepsilon_j)}),
D_{ij}(x^{(p^r\varepsilon_j)}x_i)]
=-\binom{t+1}{p^r}D_{ij}(x^{((t+1)\varepsilon_j)}).\eqno(4.3.25)$$
 Note that $\binom{t+1}{p^r}\not\equiv 0 \pmod{p}.$
 Applying $\phi$ to (4.3.25),  we obtain (4.3.21) by
 using
 (4.3.24).
\qed\newline

Now we give the following

\begin{proposition}

  Let $t=p^r,\ r>0.$ Then
$
\mathrm{Der}_{-t}(\mathcal{S},\mathcal{W})=\mathrm{span}_{\mathbb{F}}
\{(\mathrm{ad}D_i)^t\mid i\in Y_0\}.$ In particular, $
\mathrm{Der}_{-t}(\mathcal{S} )=\mathrm{span}_{\mathbb{F}}
\{(\mathrm{ad}D_i)^t\mid i\in Y_0\}.$
\end{proposition}

\noindent \textit{ Proof.}   Clearly, $(\mathrm{ad}D_i)^{p^r}  $
is a derivation of $\mathbb{Z}$-degree $-p^{r}$ for any $i\in
Y_{0} $ and $r\in \mathbb{N}.$ Let $\phi\in
\mbox{Der}_{-t}(\mathcal{S},\mathcal{W}).$ Consider the action of
$\phi$ on the element $D_{ij}(x^{((t+1)\varepsilon_i)})$ for
$i,j\in Y_0.$ Note that $\mathrm{zd}
(\phi(D_{ij}(x^{((t+1)\varepsilon_i)})))=-1.$ Suppose that
$$\phi(D_{ij}(x^{((t+1)\varepsilon_i)}))
=\sum_{r\in Y_0}a_{ijr}D_r \quad \mbox{where}\ a_{ijr}\in
\mathbb{F}.\eqno(4.3.26)$$ For any $s\in Y_0\setminus j,$ by Lemma
2.2.2, $\Gamma_k+\Gamma_s\in \mathcal{S}_0 $ for $k\in Y_1.$
Moreover,
$$[\Gamma_k+\Gamma_s,D_{ij}(x^{((t+1)\varepsilon_i)})]
=\delta_{si} \binom{t}{1 }x^{(t\varepsilon_i)}D_j=0,  $$
since $\binom{t}{1}\equiv 0 \pmod{p}.$ Applying $\phi$
to the equation above, we have
$$a_{ijs}=0 \quad\mbox{for}\ s\in Y_0\setminus j.$$
Therefore, we obtain from that (4.3.26) that
$$\phi(D_{ij}(x^{((t+1)\varepsilon_i)}))=a_{i jj}D_j.\eqno(4.3.27)$$
Observe that
$$[D_{ij}(x^{((t+1)\varepsilon_i)}),x_jD_r]=
 D_{ir}(x^{((t+1)\varepsilon_i)})\quad \mbox{for}\
 i,j,r\in Y_0 \ \mbox{with}\ r\neq i,j.$$
 Applying $\phi,$
 we obtain from (4.3.27) that
  $[a_{ijj}D_j,x_jD_r]=a_{irr}D_r.$
 Consequently,
 $$a_{ijj}=a_{irr},\ i\neq j,i\neq r.$$
 Write $a_i:=a_{irr}$ for $r\in Y_0\setminus i.$
 Put
 $$\psi:=\phi-\sum_{r\in Y_0}a_r(\mathrm{ad}D_r)^t.$$
Then $\psi\in \mbox{Der}_{-t}(\mathcal{S},\mathcal{W})$ and for all
 $ i,j\in Y_0$ with $i\neq j,$ we have
 \begin{eqnarray*}
  \psi(D_{ij}(x^{((t+1)\varepsilon_i)}) )&=& \phi(D_{ij}(x^{((t+1)\varepsilon_i)})) -
  \sum_{r\in Y_0}a_r(\mathrm{ad}D_r)^t(D_{ij}(x^{((t+1)\varepsilon_i)}) )\\
  &=&a_{ijj}D_j-a_iD_j\\
  &=&a_{i}D_j-a_iD_j\\
  &=&0.
 \end{eqnarray*}
 Lemma 4.3.3 ensures that $\psi=0;$
 that is,
 $\phi=\sum_{r\in Y_0}a_r(\mathrm{ad}D_r)^t.$ The proof is complete.
\qed\newline

Now we can describe the derivation space
$\mathrm{Der}(\mathcal{S},\mathcal{W} )$ and the derivation
algebra $\mathrm{Der}(\mathcal{S} )$.

\begin{theorem}
 $\mathrm{Der}(\mathcal{S},\mathcal{W} )
 =\mathrm{ad}(\mathcal {W})\oplus
 \mathrm{span}_{\mathbb {F}}\{ (\mathrm{ad} D_{i})^{p^{r_{i}}}
 \mid i\in Y_{0}, 1\leq r_{i} < t_{i}\}.
 $
\end{theorem}

 \noindent\textit{Proof.} This is  a direct consequence of
 Propositions 4.2.9, 4.3.2, 4.3.4 and 4.3.5. The proof is complete.\qed
 \newline

\begin{theorem}
 $\mathrm{Der}(\mathcal{S} )
 =\mathrm{ad}(\overline{\mathcal{S}}
 +\mathcal {T})\oplus
 \mathrm{span}_{\mathbb {F}}\{ (\mathrm{ad} D_{i})^{p^{r_{i}}}
 \mid i\in Y_{0}, 1\leq r_{i} < t_{i}\}.
 $
  \end{theorem}

 \noindent\textit{Proof.} This is a direct consequence of
 Propositions 4.2.10, 4.3.2, 4.3.4 and 4.3.5.
\qed

 \section{Outer derivation algebras}

 Let $\frak{g}$ be a Lie algebra. Denote by
 $\mathrm {Der_{out}}(\frak{g}):
 =\mathrm {Der} (\frak{g}) /\mathrm{ad}\frak{g}$
the outer derivation algebra of $\frak{g}$. Using the results
 obtained in Sections 3 and 4, we shall determine the outer derivation
 algebras of $\mathcal {W}$ and $\mathcal {S}.$

 Recall the our
 notations $\frak{A}:=\frak{A}(m,n;\underline{t}),$
 $W:=W(m,n;\underline{t}), $ and
$S:=S(m,n;\underline{t}),$ where $m,n$ are integers at least 3 and
$\underline{t}:=(t_{1},t_{2},\ldots,t_{m})\in \mathbb{N}^{m}$ is
an $m$-tuple of positive integers. We use still the symbols
$\pi_{i}:=p^{t_{i}}-1$, $\pi:=(\pi_{1},\pi_{2},\ldots, \pi
_{m})\in \mathbb{N}^{m}.$ Recall the canonical torus of $\mathcal
{W},$ denoted by $\mathcal
{T}:=\mathrm{span}_{\mathbb{F}}\{\Gamma_{r}\mid r\in Y\},  $ where
$\Gamma_{r}:=x_{r}D_{r}.$

 Since the full superderivation  algebra $\mathrm{Der}(\frak{A})$ is
a restricted Lie superalgebra (see [6] and [12]), one see easily
that $ D_{i}^{p^{r }}\in \mathrm{Der}(\frak{A}) $ and
$\mathrm{ad}D_{i}^{p^{r}}=(\mathrm{ad}D_{i})^{p^{r}} $ for any
$i\in Y_{0}$ and $r\in \mathbb{N}_{0},$ where ad is the adjoint
representation of the Lie superalgebra $\mathrm{Der}(\frak{A}).$
Put $\mathcal{J} :=\mathrm{span}_{\mathbb{F}}\{
D_{i}^{p^{r_{i}}}\mid i\in Y_{0}, 1\leq r_{i}< t_{i}\}.$ Note that
both $\mathcal{T}$ and $\mathcal{J}$ are   contained in the even
part of the Lie superalgebra $\mathrm{Der}(\frak{A}).$
\newline

\subsection{The outer derivation algebra of $\mathcal{W}$}

\begin{lemma} 
  $\mathcal{J}$ is a $(\sum_{i\in
Y_{0}}t_{i}-m)$-dimensional abelian Lie-subalgebra of
$\mathrm{Der}(\frak{A})$ and $[\mathcal{J},\mathcal{T}]=0.$
\end{lemma}

\noindent\textit{Proof.} Since $[ D_{i},D_{j}]=0,$ $\mathcal {J}$
is abelian. By the formula (1.2.1), $[\mathcal{J},\mathcal{T}]=0.$
Applying $D_{i}^{p^{r_{i}}}$ to all basis elements of $\frak {A} $ of the form
$x^{(\pi_{j} \varepsilon_{j})}, $    one may show that
$\{ D_{i}^{p^{r_{i}}}\mid i\in Y_{0}, 1\leq r_{i}< t_{i}\}$ is
$\mathbb{F}$-linear independent. The proof is complete.
\qed\newline

Now we can characterize the outer derivation algebra of
$\mathcal{W}.$

\begin{theorem}
$\mathrm {Der_{out}}(\mathcal{W})$ is an abelian Lie algebra of
dimension $ \sum_{i\in Y_{0}}t_{i}-m.$
\end{theorem}

\noindent\textit{Proof.} This is  a direct consequence of Theorem
3.2.11 and Lemma 5.1.1. \qed\newline

As a direct consequence of Theorem 5.1.2, we give the dimension
formula of the derivation algebra of $\mathcal{W}.$

\begin{corollary}
$\dim_{\mathbb{F}}(\mathrm{Der}(\mathcal{W}))=(m+n)2^{n-1}
p^{\sum_{i\in Y_{0}}t_{i}} +\sum_{i\in Y_{0}}t_{i}-m.$
\end{corollary}

\noindent\textit{Proof.} Note that $\dim_{\mathbb{F}} \mathcal{W}
=(m+n)2^{n-1} p^{\sum_{i\in Y_{0}}t_{i}}.$ Since $\mathcal{W} $ is
centerless, the dimension formula is a direct consequence of
Theorem 5.1.2.\qed

\subsection{The outer derivation algebra of $\mathcal{S}$}

To study the outer derivation algebra of $\mathcal {S},$ we first
establish the following lemma.

\begin{lemma}

$\mathrm{(i)}$ $\mathcal{T}\subset
\mathrm{Nor}_{\mathcal{W}}(\overline{\mathcal{S}}); $
$\mathrm{(ii)}$ $\mathcal{J}\subset
\mathrm{Nor}_{\mathrm{Der}(\frak{A} )}(\overline{\mathcal{S}});$
$\mathrm{(iii)}$ $\mathcal{T}\subset \mathrm{Nor}_{\mathcal{W}}(
\mathcal{S}); $ $\mathrm{(iv)}$ $\mathcal{J}\subset
\mathrm{Nor}_{\mathrm{Der}(\frak{A} )}( \mathcal{S} ).$
\end{lemma}

\noindent\textit{Proof.}  (i)  Since the divergence  is a
derivation of $\mathcal{W}$ to $\frak {A},$ we have
$$
\mathrm{div}([\Gamma_{r},E])=\Gamma_{r}( \mathrm{div}(E))
\quad \mbox{for all }\ r\in Y,\  E\in \mathcal{W}.
$$
The assertion follows.

(ii) Note that $D_{i}^{p^{r_{i}}}$ and $D_{s} $ commute.  Using
the formula (1.2.1), we  have
$$ [D_{i}^{p^{r_{i}}}, fD_{s}]
  = D_{i}^{p^{r_{i}}} (f)D_{s}
  \quad \mbox{where}\ i\in Y_{0}, s\in Y, f\in \frak{A}.
$$
Furthermore,
 $ \mathrm{div}([D_{i}^{p^{r_{i}}},fD_{s}])
  =D_{i}^{p^{r_{i}}} (\mathrm{div} (fD_{s}) ).
$
Since $\mathrm{div} $ is linear, (ii) holds.

 (iv) Just as in (ii), one may obtain that
$[D_{i}^{p^{r_{i}}}, D_{rs}(f)]=D_{rs}(D_{i}^{p^{r_{i}}}(f)) $
for $r,s \in Y,$  $f\in \frak{A}.$ Therefore, (iv) holds.

(iii) The proof is analogous to the one of (iv). \qed\newline

Let
$\widetilde{\mathcal{S}}:=\overline{\mathcal{S}}+\mathcal{T}+\mathcal{J}.
$
  By Lemma 5.1.1 and 5.2.1(i), (ii),
   $\widetilde{\mathcal{S}}$ is a subalgebra of $\mathrm{Der}(\frak{A}).$
   We need the following lemma.

\begin{lemma}

   $\mathcal{S}$ is an ideal of $\overline{\mathcal{S}}$; in particular,
   $\mathcal{S}$ is an ideal of $\widetilde{\mathcal{S}}.$
 \end{lemma}

\noindent\textit{Proof.}  The first assertion is  direct, since
$S$ is an ideal of $\overline{S}.$  Then the second follows from
Lemma 5.2.1(iii) and (iv). \qed

\begin{proposition} 
$\mathrm {Der_{out}}(\mathcal{S})\cong \widetilde{\mathcal{S}}/
\mathcal{S}.$
\end{proposition}

\noindent\textit{Proof.} By Lemma 5.2.2, the   mapping $
\widetilde{\mathcal{S}}\rightarrow \mathrm{Der}(\mathcal{S}), $
$E\mapsto (\mathrm{ad}E)\mid _{\mathcal{S}} $ is well defined. By
Theorem 4.3.7, it is an epimorphism of Lie algebras. Assert that
it is also injective. It suffices to show that   the centralizer
$\mathrm{C}_{ \widetilde{\mathcal{ S  }}}( \mathcal{S} )$ of
$\widetilde{\mathcal{S}}$ in $\mathcal{S} $ is trivial.

Clearly, $\mathrm{C}_ { \mathcal{W}}  (\mathcal{S}) \subset
\mathcal{G},$ since $\mathcal{S}_{-1}=\mathcal{W}_{-1}.$ Then
$\mathrm{C}_{  \mathcal{W} } (\mathcal{S}) \subset \mathrm{C}_{
\mathcal{G} }(\mathcal{S})\subset \mathrm{C}_{ \mathcal{G} }(
\mathcal{S} _{0}).$ Using Lemma 2.1.2,  one may verify by an
elementary computation that $\mathrm{C}_{ \mathcal{G} }(
\mathcal{S} _{0})=0.$  Consequently, $\mathrm{C}_ { \mathcal{W}}
(\mathcal{S})=0 $ and therefore, $\mathrm{C}_{
\overline{\mathcal{S}}+\mathcal{T} }( \mathcal{S})=0.$ On the
other hand, using the identity that $[D_{i}^{p^{r_{i}}},
D_{rs}(f)]=D_{rs}(D_{i}^{p^{r_{i}}}(f)) $ for $r,s \in Y,$ $f\in
\frak{A}(m),$ one may easily show that $\mathrm{C}_{ \mathcal{J}}(
\mathcal{S} )=0. $

Combing $\mathrm{C}_{ \overline{\mathcal{S}}+\mathcal{T} }(
\mathcal{S})=0 $ with $\mathrm{C}_{ \mathcal{\mathcal{J}} }(
\mathcal{S} )=0,$ one may obtain that $\mathrm{C}_{
\widetilde{\mathcal{ S  }}}( \mathcal{S} )=0,$  since
$\overline{\mathcal{S}}$ is spanned by certain homogeneous
elements of $\mathbb{Z}$-degree at least $-1$, but,
$\mathcal{\mathcal{J}}$  is spanned by certain elements of
$\mathbb{Z}$-degree less than $-1.$ The proof is complete. \qed\\

We have the following dimension formula for the derivation algebra of $ \mathcal{S}.$
\begin{corollary}

$$\dim_{\mathbb{F}}( \mathrm{Der}(\mathcal{S}) )
 =\left\{
\begin{array}{l}
  (m+n-1)2^{n-1}p^{\sum_{i\in Y_{0}} t_{i}}+\sum_{i\in Y_{0}}t_{i}-m+2
   \quad\quad n\ \mbox{even};
  \\
(m+n-1)2^{n-1}p^{\sum_{i\in Y_{0}} t_{i}}+\sum_{i\in
Y_{0}}t_{i}-m+1 \quad\quad n\ \mbox{odd}.
\end{array}
\right.
 $$
 \end{corollary}

\noindent\textit{Proof.} By the proof of Proposition 5.2.3,
$\mathrm{Der}(\mathcal{S})\cong  \widetilde{\mathcal{S}}.$ Note
that
 $\widetilde{\mathcal{S}}
 =\overline{\mathcal{S}}\oplus \mathbb {F}\cdot \Gamma_{1'}\oplus \mathcal{J}.$
 By Lemma
 5.1.1, it suffices
 to determine the dimension of $\overline{\mathcal{S}}.$ By the definition of
 the divergence,
 it is easily seen that
 $$\mathrm{div}(\mathcal{W})
 =\mathrm{span}_{\mathbb{F}}\{x^{\alpha}x^{u}\mid \alpha\in \mathbb {A},
 u\in \mathbb{B}, |u|\ \mbox{even}; |\alpha|+|u|<|\pi|+n\}.
 $$
 Therefore,  $\dim_{\mathbb{F}}(\mathrm{div}(\mathcal{W}))
 =2^{n-1}p^{\sum_{i\in Y_{0}} t_{i}}-1, $ if $n$ is even;
 $\dim_{\mathbb{F}}(\mathrm{div}(\mathcal{W}))
 =2^{n-1}p^{\sum_{i\in Y_{0}} t_{i}}, $ if $n$ is odd.
 Note that $\overline{\mathcal{S}}=\ker(\mathrm{div}\mid_{\mathcal {W}}) $
 and that
 $\dim_{\mathbb{F}}(\mathcal {W})=(m+n)2^{n-1}p^{\sum_{i\in Y_{0}} t_{i}}.$
 Then
 $$\dim_{\mathbb{F}}(\overline{\mathcal{S}} )
 =\dim_{\mathbb{F}}(\mathcal{W})-\dim_{\mathbb{F}} (\mathrm{div}(\mathcal{W}))
 =\left\{
\begin{array}{l}
  (m+n-1)2^{n-1}p^{\sum_{i\in Y_{0}} t_{i}}+1  \quad n\ \mbox{even};
  \\
(m+n-1)2^{n-1}p^{\sum_{i\in Y_{0}} t_{i}}  \quad\quad\quad n\
\mbox{odd}.
\end{array}
\right.
 $$
 The proof is complete.
\qed\newline

We shall describe explicitly   the structure of the outer
derivation algebra of $\mathcal{S}.$ If $n(=|\omega|)$ is even,
put $\frak{P}:=\mathrm{span}\{x^{(\pi-\pi_{i}\varepsilon_{i})}
x^{\omega}D_{i}\mid i\in Y_0\}.$
 Then $\frak{P}$ is an abelian subalgebra of $\mathcal{W}.$

 Let $V$ be an abelian Lie algebra of dimension $\sum_{i\in Y_{0}}t_{i}.$
 Let $A\in \mathrm{End}_{\mathbb{F}}(V)$ be semisimple  with eigenvalues
  $0$ and $1$ such
 that $\dim_{\mathbb{F}} V_{0}(A)= \sum_{i\in Y_{0}}t_{i}-m$
 and $\dim_{\mathbb{F}} V_{1}(A)=  m.$ Denote the semidirect product by
 $\mathcal{L}_{\mathcal{S}}:=\mathbb{F}\cdot A \ltimes V.$
 Then
 $\mathcal{L}_{\mathcal{S}}$
 is a metabelian Lie algebra of dimension
 $1+\sum_{i\in Y_{0}}t_{i}.$

\begin{theorem}

$\mathrm{(i)}$ If $n$ is odd then $\mathrm{Der
_{out}}(\mathcal{S})$ is an abelian Lie algebra of dimension
 $1+\sum_{i\in Y_{0}}t_{i}-m.$

 $\mathrm{(ii)}$ If $n$ is even then $\mathrm{Der _{out}}(\mathcal{S})$
 is isomorphic to the metabelian Lie algebra $\mathcal{L}_{\mathcal{S}} $ of dimension
 $1+\sum_{i\in Y_{0}}t_{i}. $
 \end{theorem}

 \noindent\textit{Proof.}
 According to [8, Proposition 2.8], we know that
 $\overline{\mathcal{S}}=\mathcal{S}\oplus \frak{P}$ if $n$ is even;
 $\overline{\mathcal{S}}=\mathcal{S}$ if $n$ is odd. Consequently,
 $\dim_{\mathbb{F}}(\mathcal{S})=\dim_{\mathbb{F}} (\overline{\mathcal{S}})-m$,
 if $n$ is even;
 $\dim_{\mathbb{F}}(\mathcal{S})=\dim_{\mathbb{F}} (\overline{\mathcal{S}}) $,
 if $n$ is odd. Then by Proposition 5.2.3, we have
 $$
 \dim_{\mathbb{F}} (\mathrm{Der_{out}}(\mathcal{S}))
=\left\{
\begin{array}{l}
1+\sum_{i\in Y_{0}}t_{i},  \quad\quad\quad  n\  \mbox{even};
\\
 1+\sum_{i\in Y_{0}}t_{i}-m,  \quad n\ \mbox{odd}.
\end{array}
\right.
$$
Again by Proposition 5.2.3, if $n$ is odd then $\mathrm{Der
_{out}}(\mathcal{S})\cong \mathbb {F}\cdot \Gamma_{1'}\oplus
\frak{J} $ is an abelian Lie algebra, since we have noted that
$\widetilde{\mathcal{S}}
 =\overline{\mathcal{S}}\oplus \mathbb {F}\cdot \Gamma_{1'}\oplus \frak{J}.$

 It remains to consider the case that $n$ is even. By Proposition 5.2.3, we obtain
 in this case that
 $$\mathrm{Der _{out}}(\mathcal{S})
 \cong (\mathcal{S}+\mathbb {F}
 \cdot \Gamma_{1'}+\mathcal{J}+\frak{P})/ \mathcal{S}
 =  (\mathcal{S}+\mathbb {F}
 \cdot \Gamma_{1'})/ \mathcal{S}\oplus(\mathcal{S}+\mathcal{J})/ \mathcal{S}
 \oplus(\mathcal{S}+\frak{P})/ \mathcal{S}.$$
 Note that $[\Gamma_{1'}, x^{(\pi-\pi_{i}\varepsilon_{i})} x^{\omega}D_{i}]
 =x^{(\pi-\pi_{i}\varepsilon_{i})} x^{\omega}D_{i},$
 $[\Gamma_{1'}, \frak{J}]=0,$ $[\mathcal{J}, \frak{P}]\subset \mathcal{S}.$
 Now, a straightforward verification shows that
 $\mathrm{Der _{out}}(\mathcal{S})\cong \mathcal{L}_{\mathcal{S}}.$
 The proof is complete.
\qed

\begin{remark}

{\rm According to  the known results on the outer superderivation
algebra of Lie superalgebra $W$ and $S $ (see [8,  Theorems 2.4
and 2.12]), we conclude from Theorem 5.1.2 that, for the
generalized Witt Lie superalgebra $W(m,n;\underline{t}),$ the
outer superderivation algebra coincides with the outer derivation
algebra of the even part $ W_{\overline{0}}.$ We can also conclude
from Theorem 5.2.5 that, for the special Lie superalgebra
$S(m,n;\underline{t}),$ the same conclusion holds (that is,
$\mathrm{Der _{out}}(S)\cong\mathrm{Der _{out}}(S_{\overline{0}} )
$) if and only if $n$ is even.}
\end{remark}

\noindent\textbf{Acknowledgments}\\

\noindent The authors are partially supported by  NSF  grant (10271076) of China
and NSF grant of
 Heilongjiang Province, China.

\end{document}